\documentclass[nonblindrev]{informs_diy}

\RequirePackage{tgtermes}
\RequirePackage{newtxtext}
\RequirePackage{newtxmath}
\RequirePackage{bm}
\RequirePackage{endnotes}

\OneAndAHalfSpacedXI 

\usepackage{algorithm}
\usepackage{algpseudocode}
\usepackage{tikz}

\usepackage{natbib}
 \bibpunct[, ]{(}{)}{,}{a}{}{,}%
 \def\bibfont{\small}%

\usepackage{fix-cm}
\usepackage{multirow}
\usepackage{rotating}
\usepackage{fancyvrb}
\usepackage{array}
\usepackage{booktabs}
\usepackage{accents}
\usepackage{array}
\usepackage[font=scriptsize,skip=0pt]{subcaption}
\usepackage{bm}
\usepackage{threeparttable}
\usepackage{adjustbox}
\usepackage[labelfont=bf, labelsep=period, skip=0pt]{caption}

\newcommand{\revision}[1]{{\color{black}#1}}
\newcommand{\revisiontwo}[1]{{\color{black}#1}}

\newcommand{\ubar}[1]{\underaccent{\bar}{#1}}

\EquationsNumberedThrough    

\TheoremsNumberedThrough     
\ECRepeatTheorems  %

\MANUSCRIPTNO{MSOM-2024-1328.R2}

\begin{document}

\RUNAUTHOR{Ye et al.}


\RUNTITLE{Contextual Stochastic Optimization for Order Fulfillment Optimization}

\TITLE{Contextual Stochastic Optimization for \\ Omnichannel Multi-Courier Order Fulfillment Under Delivery Time Uncertainty}

\ARTICLEAUTHORS{%
\AUTHOR{Tinghan Ye, Sikai Cheng, Amira Hijazi, Pascal Van Hentenryck}

\AFF{H. Milton Stewart School of Industrial \& Systems Engineering, Georgia Institute of Technology, \EMAIL{\{joe.ye, sikaicheng, ahijazi6, pvh\}@gatech.edu}}




%
}

\ABSTRACT{\textbf{\textit{Problem definition}}: The paper studies a
  large-scale order fulfillment problem for a leading e-commerce
  company in the United States. The challenge involves selecting
  fulfillment centers and shipping carriers with observational data
  only to efficiently process orders from a vast network of physical
  stores and warehouses. The company's current practice relies on
  heuristic rules that choose the cheapest fulfillment and shipping
  options for each unit, without considering opportunities for
  batching items or the reliability of carriers in meeting expected
  delivery dates.  \textbf{\textit{Methodology / results}}: The paper
  develops a data-driven Contextual Stochastic Optimization (CSO)
  framework that integrates distributional forecasts of delivery time
  deviations with stochastic and robust order fulfillment optimization
  models. The framework optimizes the selection of fulfillment centers
  and carriers, accounting for item consolidation and delivery time
  uncertainty. Validated on a real-world data set containing tens of
  thousands of products, each with hundreds to thousands of fulfillment options,
  the proposed CSO framework significantly enhances the accuracy of
  meeting customer-expected delivery dates compared to current
  practices. It provides a flexible balance between reducing
  fulfillment costs and managing delivery time deviation risks, emphasizing the
  importance of contextual information and distributional forecasts in
  order fulfillment.  \textbf{\textit{Managerial implications}}: This
  is the first study of omnichannel multi-courier order
  fulfillment problem with delivery time uncertainty through the lens
  of contextual optimization, fusing machine learning and
  optimization. The results offer actionable guidance for retailers to enhance service quality and customer satisfaction while balancing cost efficiency and risk, supporting higher retention and profitability.
}
\FUNDING{This research was partly supported by the NSF AI Institute for Advances in Optimization (Award 2112533). }




\maketitle


\section{Introduction}

Over the past decade, online shopping has surged dramatically. The
e-commerce share of total retail sales in the United States has risen
from 6\% in 2015 to 16\% in 2024, exceeding 1.1 trillion dollars
\citep{us-department-of-commerce-2024}. Amid intense competition in
pricing, service, and marketing, driven by the lucrative e-commerce
market, logistical performance has become a key element to
success. Last-mile fulfillment, in particular, is critical as it
accounts for 55\% of the total transportation costs in the fulfillment
process \citep{viswanathan2023elevating}. This evolving market
landscape poses unprecedented challenges for e-commerce companies
striving to deliver satisfactory customer services. Fortunately, the
growing volume of data also offers opportunities for data-driven
operational improvements, enabling companies to boost revenue and cut
costs \citep{fisher-2019}.

In addition, after the COVID-19 pandemic, the demand for timely home
delivery has also increased significantly. As a result, omnichannel
fulfillment and the use of multiple couriers, including
crowd-shipping, have become essential components of the fulfillment
landscape, drawing significant attention both in industry and academia
\citep{dethlefs-2022, mohri-2023, das-2023}. The omnichannel approach
combines distribution centers (DCs) and brick-and-mortar stores to
fulfill online orders. Unlike the traditional method that uses DCs for
online orders and stores for in-store purchases separately, the
omnichannel approach leverages the proximity of physical stores to
delivery addresses. This can potentially offer more cost-effective and
faster delivery options for online customers
\citep{acimovic2019fulfillment}. Meanwhile, in the crowd-shipping
model, retailers recruit non-dedicated drivers, such as gig couriers
like Roadie, to handle pick-up and delivery tasks for online
customers. This approach can alleviate the workload of employed
drivers during peak hours and offer faster delivery to customers,
resulting in higher service satisfaction \citep{zehtabian-no-date,
  behrendt2023prescriptive}.  {\em The problem considered in this paper
involves fulfilling orders by choosing from various types of
fulfillment centers, such as warehouses and stores, and different
couriers, including traditional shipping companies and gig couriers.}

Given that transportation and shipping are significant cost drivers in
e-commerce \citep{dethlefs-2022, kuhn-2013, hubner-2013}, delivery
consolidation has emerged as an effective strategy for fulfilling
orders with multiple items. This approach has been shown to enhance
consumer satisfaction \citep{wagner-2023}, while also reducing
shipping costs for online retailers \citep{chen-2024}. In addition to
cutting parcel expenses, item consolidation enhances the customer
experience by eliminating the need to receive multiple packages for a
single order. This approach not only improves convenience but also
promotes a more sustainable and environmentally friendly fulfillment
process by decreasing the number of packages shipped and the number of
delivery trips needed. This reduction in shipments and trips leads to
lower carbon emissions and less packaging waste, contributing to
greener logistics practices \citep{ulku2012dare}. {\em This paper also
investigates multi-item order consolidation in its omnichannel
fulfillment problem.}

Most existing studies on order fulfillment optimization have not
explicitly considered the uncertainty in delivery timeliness. However,
\cite{freedman-2019} highlight the strong desire of customers for
accurate delivery time promises and faster service. The survey
indicates that around 44\% of online shoppers have abandoned their
shopping carts because the items they wanted would not arrive on time,
while 20\% abandon orders with unclear delivery dates. Furthermore,
28\% of shoppers are willing to pay for expensive, expedited shipping
to receive products at their desired time. The statistics clearly show
that customers prioritize service quality, sometimes even more than
the product value, when making purchase decisions. Moreover,
\cite{salari2022real} point out the importance of considering delivery
time uncertainty in order fulfillment to better meet customer
expectations and improve satisfaction. Their proposed promised
delivery time policy demonstrates a 6.1\% increase in sales volume in
numerical experiments.  Notably, overpromising delivery times can lead
to higher return rates, lower repurchase rates, and dissatisfaction
with delayed service. Additionally, customers can have mixed
sentiments over underpromised delivery times, as early deliveries may
not align with their availability. Therefore, missing the exact
promised delivery date can negatively impact short-term and long-term
sales, both implicitly and explicitly \citep{salari2022real, cui2024sooner}. {\em Given the crucial role of delivery timeliness, the
paper incorporates delivery time uncertainty into the omnichannel
multi-courier order fulfillment optimization problem through a
contextual stochastic optimization framework.}

One of the key challenges of Omnichannel Multi-Courier Order
Fulfillment Optimization comes from the nature of historical data:
\revision{only observational data---i.e., delivery time deviations for the fulfillment options actually chosen---are available. Counterfactual deviations for unchosen options are never observed, similar to the bandit feedback setting \citep{lattimore2020bandit}.} {\em The availability of observational data only, raises fundamental
challenges that are typically not addressed in prior work.}

To address these challenges, {\em this paper proposes a Contextual
Stochastic Optimization (CSO) framework, encompassing both
risk-neutral and robust approaches.}  Unlike traditional stochastic
programming methods that depend on static, predefined probability
distributions for uncertainties, CSO combines machine learning
techniques with mathematical optimization to dynamically predict
uncertainties based on contextual information. Thanks to the
availability of contextual information, CSO has the flexibility of
generating more tailored and effective decisions, leveraging the most
current and relevant information. Such flexibility is particularly
beneficial in the dynamic environment of online retail, where
conditions and customer behaviors can change rapidly, requiring a more
responsive approach to decision making. The effectiveness of CSO has been
validated in the existing operations research and management science
literature \citep{mivsic2020data, sadana2024survey}. This approach has
proven successful in various applications, including inventory
control \citep{bertsimas-2016, meller-2018}, price
and revenue management \citep{ito-2016, perakis2023robust}, and, more recently, school redistricting \citep{guan2024contextual}.

\subsection{Contributions}

{\em The key novelty in this paper is a CSO framework for an
  omnichannel multi-courier order fulfillment optimization with
  observational data only.}   
  Existing studies on CSO typically assume
the availability of a dataset $\mathcal{D} = \{(\mathbf{s}_1,
\mathbf{c}_1), (\mathbf{s}_2, \mathbf{c}_2), \ldots \}$, pairing
covariates with fully observed uncertain parameter
vectors. However, real-world scenarios often involve high-dimensional
uncertain parameters where individual data points capture only partial
information, particularly when uncertainty measurement depends on the
chosen solution, which is the case in the application considered in
this paper.

To address the challenge of observational data, this paper introduces
a generalized CSO framework addressing these real-world
complexities. {\em The novel CSO framework is applied to a
  sophisticated stochastic omnichannel order fulfillment optimization
  problem with multiple carrier options and delivery time
  uncertainties.} In collaborating with a major US online retailer,
the framework is rigorously evaluated on an industrial order
fulfillment dataset featuring a vast, intricate fulfillment
network. This application demonstrates the framework capabilities and
its potential impact on large-scale, real-world optimization
challenges. \revisiontwo{Based on the significant performance improvements demonstrated in this study, the industrial partner is actively considering a pilot program to integrate key elements of the proposed framework into their production fulfillment optimization engine.} The contributions of the paper can be summarized as
follows.
\revision{
\begin{enumerate}
  \item \revisiontwo{ {\em A modular and scalable CSO framework for observational data}:  
    The paper presents a novel, distribution‑agnostic contextual stochastic optimization (CSO) framework that converts problems with missing‑counterfactual observational data into tractable problems. A contextual distribution oracle is learned from partial observations, and the resulting problem is solved via two scalable solution paradigms---Contextual Sample Average Approximation (C‑SAA) and Contextual Robust Optimization (C‑RO). This features a tractable formulation for assignment-type problems with piecewise linear objectives.}
  \item {\em Tailored contextual distribution learning for e-commerce delivery time deviations}: 
    The paper adapts calibrated probabilistic multi‐class classification and tree‐based quantile regression to exploit the discrete, ordinal nature of delivery time deviations, yielding well‑calibrated contextual distributions for the CSO framework.
  \item {\em Consolidation‑ and timeliness‑aware MILP for omnichannel, multi‑courier fulfillment}:
    The paper proposes the first data‑driven mixed‑integer linear program (MILP) that jointly selects fulfillment centers and carriers and captures item‑consolidation discounts while embedding learned contextual deviation distributions in a tractable way.
  \item {\em Industry‐scale real‐world validation}:
    The paper includes an extensive case study at a major e‐commerce retailer to (1) validate the efficacy of the proposed approaches in reducing expected fulfillment costs and (2) show significantly improvements in on‐time delivery rates compared to standard heuristics and baseline models. These results demonstrate the potential for substantial real‑world impact.
\end{enumerate}
}

\revision{
\subsection{Structure of the Paper}
The remainder of this paper is organized as follows. Section \ref{sec: lit-review} reviews related literature. Section \ref{sec:application} formulates the omnichannel multi‑courier order‑fulfillment problem under delivery‑time deviation uncertainty. Section \ref{sec:methodology} develops a generic contextual stochastic optimization framework for observational data and introduces two tractable solution paradigms---Contextual Sample Average Approximation and Contextual Robust Optimization. Section \ref{sec: learning-app} describes the machine‑learning methods (multi‑class classification and tree‑based quantile regression) used to implement the contextual distribution oracle. Section \ref{sec: experiments} presents the case study, computational results, and actionable insights from a real dataset. Finally, Section \ref{sec:conclusion} concludes and outlines directions for future research.
}

\section{Related Literature} \label{sec: lit-review}
There are several streams of literature related to this paper in terms of methodology and problem context.

\subsection{Data-Driven Decision Making}

This paper adds to the growing body of literature on data-driven
decision-making, particularly at the intersection of machine learning
and optimization, which has led to the development of
CSO. \cite{sadana2024survey} put forth a comprehensive survey on
contextual optimization. Tutorials on various methodologies and
applications of CSO models in logistics and operations management are
available in the works of \cite{mivsic2020data},
\cite{qi2022integrating}, and \cite{tian2023tutorial}. Among the body
of research on CSO, the most relevant to this paper are those that
focus on learning contextual distributions of uncertain parameters
and incorporating these forecasts into optimization
models. Data-driven approaches typically build these forecasts from
empirical residual errors \citep{deng2022predictive,
  kannan2023residuals, perakis2023robust} or weighted empirical
distributions based on proximity to the training data
\citep{bertsimas2020predictive, notz2022prescriptive}.  However, these
studies assume full observation of uncertain parameters for each data
point, which differs from the context of this paper. Recent research
has also explored integrating downstream optimization effects into ML
model training. \cite{elmachtoub2022smart} introduce the Smart Predict
then Optimize (SPO) loss and its convex surrogate for point prediction
models, with subsequent extensions to contextual distributional
estimation \citep{qi2021integrated, kallus2023stochastic}. However,
the SPO loss framework is not directly applicable to this paper's
context, as the dataset analyzed here is observational. In such
datasets, outcomes are only partially observed based on decisions
made, presenting unique challenges not addressed by these existing
methods.

This paper thus contributes to the field of data-driven
decision-making that handles observational data, a domain explored by
prior research, e.g., in \citet{bertsimas2019optimal} and \citet{jo2021learning},
using optimal prescriptive trees for prescriptive analytics in
personalized treatment. While these methods effectively learn policies
mapping features to actions, they cannot be directly applied to the
fulfillment problem considered in this paper due to the presence of
instance-specific hard constraints on decisions. The proposed
framework addresses these limitations, offering a novel approach to
CSO problems with observational data and constrained decision spaces.

Unlike earlier CSO applications that focus on risk-neutral objectives,
a growing stream of research is concerned with risk-averse
optimization models that incorporate contextual information
\citep{bayram2022optimal, pervsak2023contextual, patel2024conformal,
  sun2024predict}.  Building on these studies, this paper also adapts
traditional robust optimization models to tackle a unique CSO problem
that differs from those explored in other research.

\subsection{Order Fulfillment Problem}

The novelty in this paper is the application of CSO to a unique
instance of the order fulfillment problem. The order fulfillment
problem has been extensively studied, and
\cite{acimovic2019fulfillment} provide a tutorial on related
algorithms developed for addressing the fulfillment optimization
problem. In terms of problem context, the most relevant works are
those related to omnichannel and multi-item order fulfillment
\citep{xu2009benefits, acimovic2015making, jasin2015lp, zhao-2020,
  wei2021shipping, ma2023order, chen-2024}. However, these works do
not address fulfillment involving a variety of carrier shipping
options.

Among the research on order fulfillment, demand uncertainty is the
most commonly studied source of randomness \citep{jasin2015lp,
  zhao-2020, das-2023, devalve2023understanding, ma2023order}. In this
paper, demand uncertainty is not considered because the focus is on a
single-period order fulfillment problem, rather than a forward-looking
approach that spans multiple periods.

A few studies do address delivery time uncertainty at various stages
of the order fulfillment process or throughout the entire process. For
instance, \cite{raj2024stochastic} propose an integrated queuing model
that accounts for delivery time uncertainty in different parts of the
order life cycle.  They incorporate these uncertainties into an
optimization model that minimizes costs while adhering to a delivery
probability constraint to ensure on-time delivery.

In contrast to using queuing models, this study is particularly
relevant to data-driven, ML-based methods that address for delivery
time uncertainty in order fulfillment. Specifically,
\cite{liu2021time} integrate point predictions of travel time with a
last-mile delivery order assignment optimization problem. They propose
both a SAA model and a distributionally robust optimization (DRO)
model, with the latter constructing a moment-based ambiguity
set. Similarly, \cite{kandula2021prescriptive} explore a last-mile
delivery problem, focusing on delivery success as the
uncertainty. They map success probabilities to delivery times, which
are then integrated into a Vehicle Routing Problem with Time
Windows. \cite{bayram2022optimal} study a robust order batching
optimization problem in warehouse picking and packing. They use random
forest to predict order processing time to form an uncertainty set
within a robust optimization model. This paper builds on similar
concepts, developing tractable approximations and reformulations of
integer programming-based optimization models. More precisely, the
uncertainty in the order fulfillment problem considered in this paper
is the potential order delivery deviations, which can occur at various
stages of the fulfillment process. A similar uncertainty addressed
using data-driven approach can be found in \cite{salari2022real}, who
develop a quantile regression forest-based method to generate
distributional forecasts of order delivery times. Their approach then
applies a decision rule, designed to balance the asymmetric costs of
early and late deliveries, to produce expected arrival times for
customers.  This paper differs in two key aspects. First, the
distributional forecasts are directly incorporated into contextual
optimization models, eliminating the need for a separate decision rule
to generate point predictions. Second, the CSO framework proposed in
this paper is more comprehensive, accommodating multiple carrier
options, whereas \cite{salari2022real} focus on a single-carrier
scenario.

\section{\revision{The Omnichannel Multi-Courier Order Fulfillment Problem}} \label{sec:application}

This section introduces a real-world order fulfillment problem encountered by the industrial partner. 

\subsection{High-Level Problem Description}

\begin{figure}[!ht]
\centering
    \includegraphics[width=0.8\textwidth]{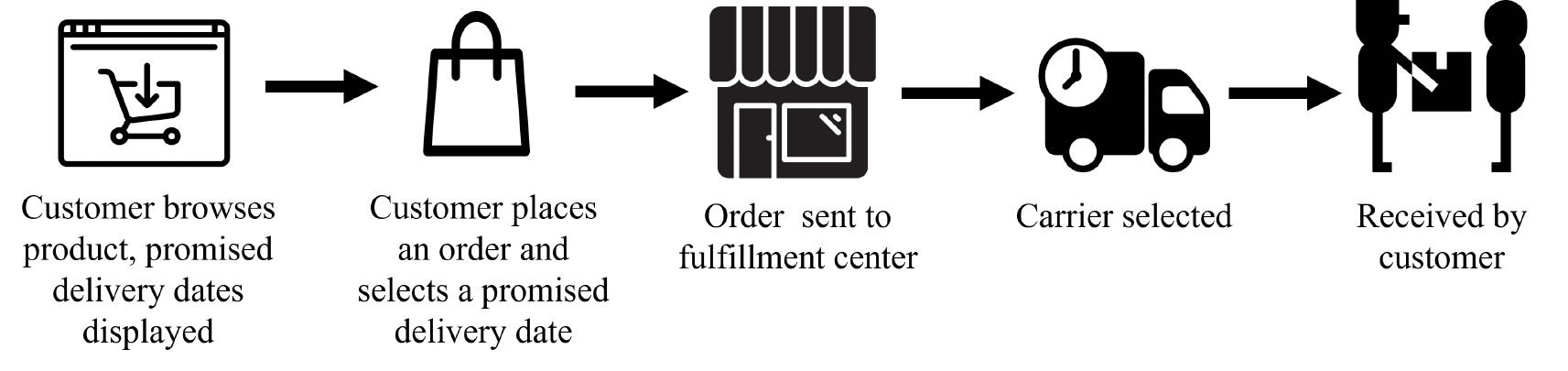} 
    \caption{The Life Cycle of the Order Fulfillment Process.}
    \label{fig:flowchart}
\end{figure}

This section provides a high-level overview of the online order
fulfillment challenges faced by the industrial partner. Similar to
other major e-commerce companies, tens of thousands of orders are
placed continuously throughout the day on the website and await
processing. Each incoming order is organized by stock-keeping units
(SKUs), resulting in one or more order lines. For each order line,
there are hundreds of fulfillment centers across the company network
eligible to source the SKUs. These centers include primarily physical
stores, as well as regional distribution centers and metro e-commerce
centers. Each location is equipped with multiple eligible carriers
capable of providing different levels of service, such as
Ground, Same-Day, and Overnight shipping. Throughout this paper, ``carrier" refers to one of these specific carrier-service pairs. The set of available
carriers spans from major shipping companies to crowdsourced gig
carriers. Figure \ref{fig:flowchart} illustrates the order fulfillment
life cycle within the operations of the industrial partner.

Before making an order fulfillment decision, the company has access to
both the transactional details of the order and the real-time status
of the fulfillment network. Specifically, the order information
includes the date and time the order was placed, the type, dimensions,
weight, and quantities of the SKUs, the shipping destination, and the
desired delivery date. The network status data provides details on the
on-hand inventory levels of the requested SKUs and the remaining open
capacities at various fulfillment centers. In addition, each carrier
provides an estimated table for shipping costs and planned
time-in-transit.

Given this information, the problem is to decide which fulfillment
center to source each order line from and which carrier to use for
shipping from that location. The goal is to identify the most
cost-effective fulfillment option (a location-carrier pair) while
ensuring the customer's desired delivery date is met. The uncertainty
of this problem lies in the delivery timeliness performance of each
fulfillment option.

\subsubsection{Overview of the Current Practice}

Currently, the industrial partner uses a fulfillment optimization engine that employs
a straightforward greedy algorithm. For each unit in an order, the
fulfillment optimization engine sorts the eligible location-carrier pairs by the
estimated shipping costs and transit times.

There are several limitations with this approach. First, since the
engine makes decisions for each unit independently, it misses the
opportunity to consolidate multiple units of the same or different
products within an order. This can lead to higher shipping costs, as
the potential to reduce the number of packages---and thereby lower
shipping expenses---is not fully realized.

Second, the method relies solely on static transit time estimates
provided by the carriers based on their service levels to select the
best carrier for meeting the promised delivery date. This approach
assumes that these estimates are always accurate. However, static
transit times may not accurately reflect the actual delivery
performance of carriers. Historical data from the industrial partner
reveals that over 10\% of orders experience varying degrees of early
or late deliveries, which are called {\em deviations} in this
paper. Given the complexity of the entire fulfillment process, these
deviations can be attributed to numerous factors, such as warehouse
pick-and-pack performance, late carrier arrivals, no-shows, and
external factors like extreme weather or other unexpected
circumstances. Such deviations from customer expectations can
negatively impact their experience, potentially harming the company in
the long run.

\subsection{The Nominal Optimization Problem}

This section introduces the nominal MILP model developed for
optimizing the omnichannel multi-courier order fulfillment problem. A comprehensive list of the notations
presented in this section is provided in Appendix A
for reference.

\subsubsection{Main Decision Variables}

Let $\mathcal{I}$
be the set of SKUs and $\mathcal{V}$ be the set of fulfillment
options. Here, $\mathcal{V}$ is defined as the Cartesian product $
\mathcal{K} \times \mathcal{L}$, where $\mathcal{K}$
denotes the set of carriers, and $\mathcal{L}$
denotes the set of locations.

Each order needs to solve an independent
optimization problem to decide the selection of carrier-location
pairs for each SKU. The primary decision variables are denoted by the vector
$\mathbf{z}= (z_{ik\ell})_{i \in \mathcal{I}, k \in \mathcal{K},
  \ell \in \mathcal{L}}$, where each $z_{ik\ell} \in \mathbb{Z}$
indicates the quantity of SKU $i$ sourced from location $\ell$ and shipped by carrier $k$.

\subsubsection{Objective Function}

The objective of the optimization problem for an order is to minimize
the overall fulfillment costs while ensuring timely delivery. The
fulfillment costs for an order include the fixed fulfillment (non-parcel)
costs and the shipping (parcel) costs, with consolidation discounts
applied when multiple units are assigned to the same carrier-location
pair. Delivery timeliness is measured by deviations from the
desired delivery date, whether early or late. To encourage
timeliness, the objective imposes penalties for both early and late
deliveries.

The objective function is formally defined as follows. Let $c_{ik \ell}^{ship} > 0$ be the per-unit shipping cost when using location
$\ell$ and carrier $k$ to fulfill SKU $i$. For each
location $\ell$, let $c_{\ell}^{fixed} > 0$ be the per-unit fixed
fulfillment cost. For each carrier $k$, let $\beta_{k }\in (0, 1)$ be the
discount percentage applied to shipping costs when multiple units are
shipped using carrier $k$ from the same fulfillment center. Let
$\gamma^+ \geq 0$ and $\gamma^- \geq 0$ be constants that convert late
and early delivery penalties into per-unit costs, respectively.

In an ideal setting where all carrier-location deviations are known at the
time of fulfillment, let the vector $\mathbf{d} = (d_{k \ell})_{k
  \in \mathcal{K}, \ell \in \mathcal{L}}$ represent the realized
uncertain delivery deviations, where each $d_{k \ell}$ denotes the number of days that the carrier-location pair
$(k, \ell)$ deviates from desired delivery date. A positive deviation
indicates late delivery, a zero deviation indicates on-time delivery, and
a negative deviation indicates early delivery.

The weighted objective function for an order can be constructed as follows:
\begin{equation} \label{eq:nonlinear-obj}
    \begin{aligned}
        g(\mathbf{z}, \mathbf{d}) = 
        \begin{cases}
        \sum_{i \in \mathcal{I}} \sum_{k \in \mathcal{K}}  \sum_{\ell \in \mathcal{L}} [c_{\ell}^{fixed} + (1- \beta_k)c_{ik \ell}^{ship} +  \gamma^+ d_{k \ell}^+ + \gamma^- d_{k \ell}^-]  z_{ik \ell}  & \text{if } \sum_{i \in \mathcal{I}} z_{ik \ell} \geq 2,\\
        \sum_{i \in \mathcal{I}} \sum_{k \in \mathcal{K}}  \sum_{\ell \in \mathcal{L}} (c_{\ell}^{fixed} + c_{ik \ell}^{ship} +  \gamma^+ d_{k \ell}^+ + \gamma^- d_{k \ell}^-) z_{ik \ell} & \text{otherwise},
        \end{cases}
    \end{aligned}
\end{equation}
where $d_{k \ell}^+ = \max \{0, d_{k \ell} \}$ and $d_{k \ell }^-
= \max \{0, -d_{k \ell} \}$.  This objective function encourages
item consolidation within an order by applying a carrier-specific
discount if at least two units are sourced from the same location and shipped by the same carrier. In addition, it leverages an asymmetric
penalty cost function for delivery deviations, where $\gamma^+ \gg
\gamma^-$. This asymmetry strongly penalizes late deliveries to discourage delays while also applying a smaller penalty to disincentivize early deliveries. Specifically, the penalties can be interpreted as the opportunity cost associated with early deliveries and the lost-sale cost associated with late deliveries. The objective
function is piecewise linear and can be linearized using standard
rewritings.

\subsubsection{Operational Constraints}
Several operational constraints ensure that the selection of fulfillment center locations and carriers is feasible in practice. For location $\ell \in \mathcal{L}$, let $\operatorname{inv}_{i\ell} \geq 0$ be the inventory level of SKU $i$ at location $\ell$ and $\operatorname{cap}_{\ell} \geq 0$ be the available capacity at location $\ell$ when the order is placed. Let $q_{i} \geq 0$ be the quantity of SKU $i$ in the order, and let $e_{ik \ell} \in \{0, 1\}$ be a binary indicator of whether the location-carrier pair $(\ell, k)$ is an eligible fulfillment option for SKU $i$. An assignment $\mathbf{z}$ in the feasible set $\mathcal{Z}$ must obey the following constraints.

\begin{itemize}
\item Each unit of an SKU must be sourced from exactly one location and shipped by exactly one carrier.
        \begin{equation}
            \sum_{k \in \mathcal{K}} \sum_{\ell \in \mathcal{L}}  z_{ik \ell}= q_{i}, \quad \forall i \in \mathcal{I}
        \end{equation}
        
    \item Each SKU can only be fulfilled by eligible candidate location-carrier pairs, which may vary depending on the carrier's coverage and package requirements.
    \begin{equation} \label{eq:eligibility}
         z_{ik \ell} \leq q_{i}e_{ik\ell}, \quad \forall i \in \mathcal{I}, k \in \mathcal{K}, \ell \in \mathcal{L}
    \end{equation}
        
\item Only locations with available inventory of the requested SKUs can be chosen.
        \begin{equation}
             \sum_{k \in \mathcal{K}} z_{ik \ell} \leq \operatorname{inv}_{i\ell}, \quad \forall i \in \mathcal{I}, \ell \in \mathcal{L}
        \end{equation}

\item Each location has a daily capacity limit on the total number of units it can process.
        \begin{equation}
            \sum_{i \in \mathcal{I}} \sum_{k \in \mathcal{K}}   z_{ik \ell} \leq \operatorname{cap}_{\ell}, \quad \forall \ell \in \mathcal{L}
        \end{equation}

\end{itemize}

\subsection{The Contextual Stochastic Order Fulfillment Problem}

In practice, the delivery deviations of each carrier-location pair cannot
be perfectly estimated in advance. To account for this uncertainty for
each order, this section models the vector of uncertain delivery
deviations $\tilde{\mathbf{d}}= (\tilde{d}_{k \ell})_{k \in
  \mathcal{K}, \ell \in \mathcal{L}} \sim \mathbb{P}$ as a
$(|\mathcal{K}||\mathcal{L}|)$-dimensional random vector governed by a joint probability
distribution $\mathbb{P}$. \revision{Since deviations measured in days only take on discrete values in the
dataset, each $\tilde{d}_{k\ell}$ inherently follows a discrete probability distribution. 
It is further assumed that every $\tilde{d}_{k \ell}$ shares the same finite, ordered support. To formalize this, let $C$ be the total number of distinct deviation values and denote them by $\xi_1 < \xi_2 < \ldots < \xi_{C}$. By definition, any realization $d_{k \ell}$ of $\tilde{d}_{k\ell}$ must lie in this set: $d_{k \ell} \in \{\xi_1, \xi_2, \ldots, \xi_{C}\}.$}

In addition, let $\mathbf{s} =
(\mathbf{s}_{k \ell})_{k \in \mathcal{K}, \ell \in \mathcal{L}}$ denote the state vector observable to the order before the fulfillment decision is made, where each $\mathbf{s}_{k
  \ell}$ is a vector of covariates
encompassing information about
order details and fulfillment network status related to carrier-location pair $(k, \ell)$.  The
contextual distribution of the deviation vector, upon observing
$\mathbf{s}$, is denoted by $\mathbb{P}(\mathbf{d} \mid
\mathbf{s})$, and the individual deviation contextual distribution is
$\mathbb{P}_{k \ell}(d_{k \ell} \mid \mathbf{s}_{k \ell})$ for each
$k\in \mathcal{K}$ and $\ell \in \mathcal{L}$.

This leads to the following contextual stochastic order fulfillment problem (CSOFP):
\begin{equation} \label{eq: CSOFP}
\begin{aligned}
\text{(CSOFP)} \quad & \min_{\mathbf{z} \in \mathcal{Z}} \mathbb{E}_{\tilde{\mathbf{d}} \sim \mathbb{P}(\mathbf{d} \mid \mathbf{s})}[g(\mathbf{z},  \tilde{\mathbf{d}})] .
\end{aligned}
\end{equation}

Solving Problem (\ref{eq: CSOFP}) in practice requires leveraging data from past orders. Unlike classical CSO formulations where the entire random cost vector is observed, CSOFP is based on an observational dataset of past decisions and their realized costs. Let $\mathcal{O}$ be
the set of historical orders. For each $o \in \mathcal{O},$ define the decision vector $\mathbf{z}_o = (z_{i k \ell, o})_{i \in \mathcal{I}, k \in
\mathcal{K}, \ell \in \mathcal{L}}.$ Since deviations were observed
only for chosen $(i, k, \ell)$ triples, the full deviation vector $\mathbf{d}_o=(d_{k\ell, o})_{k \in \mathcal{K}, \ell \in \mathcal{L}}$ contains many unobserved entries. The subset of realized deviations for order $o$ is: 
\[
\mathbf{D}_o : = \{ d_{k\ell, o} \mid \sum_{i \in \mathcal{I}}z_{ik \ell, o} > 0, k \in \mathcal{K}, \ell \in \mathcal{L} \}.
\]
The resulting dataset available for estimating the contextual distribution $\mathbb{P}(\mathbf{d} \mid \mathbf{s})$ is $\mathbb{D} =
\{(\mathbf{s}_o, \mathbf{D}_o) \}_{o = 1}^{|\mathcal{O}|}$.

\section{\revision{Methodology: A CSO Framework for Observational Data}} \label{sec:methodology}

\begin{figure}
    \centering
    \includegraphics[width=\linewidth]{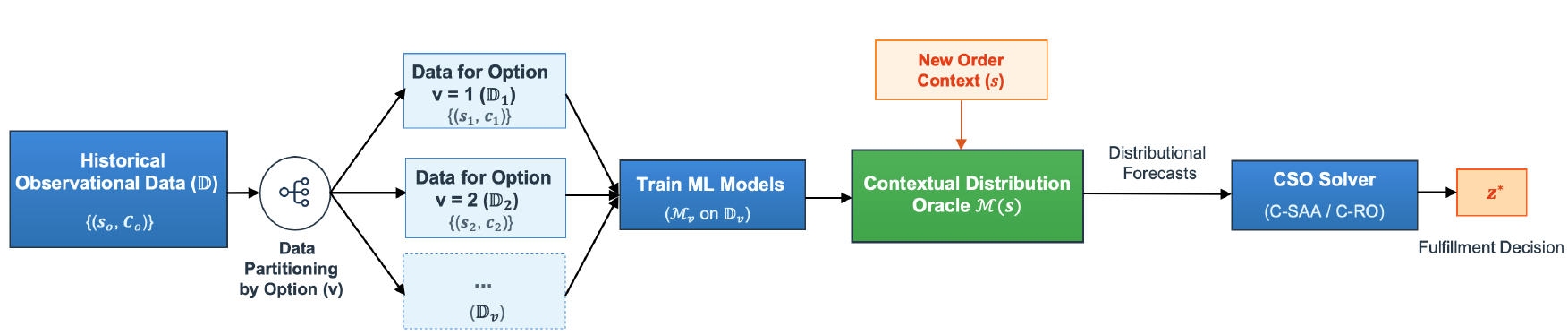}
    \caption{\revisiontwo{A Schematic Overview of the CSO Framework for Observational Data.}}
    \label{fig:schematic}
\end{figure}

\revision{
Solving CSOFP under observational data requires:
(1) learning a contextual distributional oracle from partial realizations, and
(2) optimizing with that oracle in a tractable way.  

This section lays out a contextual stochastic optimization (CSO) framework that does exactly that. In Section \ref{sec:generic_cso}, a generic CSO problem (denoted by ${\cal P}$) is defined. Section \ref{sec:learn_from_observational} shows how to construct a contextual distribution oracle  \(\mathcal M(\mathbf s)\) from partial observations. Sections \ref{sec:c-saa} and \ref{sec:c-ro} introduce two solution paradigms---Contextual Sample Average Approximation and Contextual Robust Optimization---which use \(\mathcal M(\mathbf s)\) as input. Finally, Section \ref{sec:specialize} maps all of the above back to the omnichannel multi-courier CSOFP introduced in Section \ref{sec:application}. }
\revisiontwo{Figure~\ref{fig:schematic} provides a schematic overview of the proposed framework, illustrating the flow from historical observational data through the learning and optimization stages to produce a final decision.
}

\subsection{Generic Contextual Stochastic Optimization (CSO) Problem} \label{sec:generic_cso}

Problem ${\cal P}$ consists of fulfilling a set of incoming requests
sequentially, each of which specifies demands for a subset of products. Let \(\mathcal I\) be the set of products, \(\mathcal V\) the set of selection options (e.g., fulfillment or transport modes), and \(\mathcal Z\subseteq\mathbb Z_+^{|\mathcal I|\times|\mathcal V|}\) the feasible assignment set. \revision{Each request arrives with known demand quantities \(\{q_i : i\in\mathcal I\}\) at decision time.} The decision variable
\(
  \mathbf z = (z_{iv})_{i\in\mathcal I,\;v\in\mathcal V}\in\mathcal Z
\)
assigns quantity \(z_{iv}\) of product \(i\) to option \(v\).  

\revision{
Let \(\tilde{\mathbf c}=(\tilde c_{iv})_{i\in\mathcal I,\;v\in\mathcal V}\sim\mathbb P\) denote the random cost vector. Assume that \(\tilde{\mathbf c}\) is element-wise positive and the cost function \(g(\mathbf z,\tilde{\mathbf c})\) is piecewise linear in \(\mathbf z\).  Introduce, for \(j=0,1,\dots,n\), the breakpoint vectors
\[
  \mathbf b^j = \bigl(b^j_{iv}\bigr)_{i\in\mathcal I,\;v\in\mathcal V},
  ~
  \mathbf b^0 < \mathbf b^1 < \cdots < \mathbf b^n
  \quad(\text{component-wise}),
\]
so that each segment is defined as
\[[\mathbf b^{j-1},\,\mathbf b^j]
  = \bigl\{\mathbf z : b^{\,j-1}_{iv}\le z_{iv}\le b^j_{iv}
  \;\forall\,i\in\mathcal I,\;v\in\mathcal V\bigr\}.
\]
Let \(\mathbf a^j=(a^j_{iv})_{i\in\mathcal I,\;v\in\mathcal V}\) be the marginal‐cost adjustment on segment \(j\).  Then
\begin{equation} \label{eq:piecewise-linear}
  g(\mathbf z,\tilde{\mathbf c})
  = (\tilde{\mathbf c} + \mathbf a^j)^\mathsf{T}\mathbf z
  \quad\text{if}\quad
  \mathbf z \in [\mathbf b^{j-1},\mathbf b^j],\;j=1,\dots,n,
\end{equation}
thereby capturing volume‐dependent pricing behaviors.

For each request, the fulfillment decision solves the risk‐neutral stochastic optimization problem
\begin{equation} \label{eq: so}
  \min_{\mathbf z\in\mathcal Z}
  \;\mathbb E_{\tilde{\mathbf c}\sim\mathbb P}
  \bigl[g(\mathbf z,\tilde{\mathbf c})\bigr].
\end{equation}
The piecewise‐linear cost in \eqref{eq:piecewise-linear} partitions the decision space into shipment‐size tiers determined by the breakpoints \(\{\mathbf b^j\}\).  Within tier \(j\), the marginal cost is \(\tilde{\mathbf c}+\mathbf a^j\), which shifts to \(\tilde{\mathbf c}+\mathbf a^{j+1}\) upon entering tier \(j+1\).  This formulation accommodates economies of scale (declining unit costs at higher volumes), quantity discounts (step‐down pricing beyond thresholds), capacity surcharges (unit‐cost premiums when exceeding lower‐cost bands), and other tiered‐pricing schemes.
}

In Problem \eqref{eq: so}, the true distribution $\mathbb{P}$ is not
readily available. However, side information often exists for each product-option pair.
Let $\mathbf{s}_{iv}$ denote the covariate vector for product $i$ under selection option $v$ (e.g., network state, weather, current demand). All these covariate vectors are aggregated into the
state vector $\mathbf{s} = (\mathbf{s}_{iv})_{i \in \mathcal{I}, v
  \in \mathcal{V}}$, which encompasses all contextual information
available at the time of decision. The distribution $\mathbb{P}$ is then approximated by the conditional
distribution $\mathbb{P}(\mathbf{c} \mid \mathbf{s})$, yielding the Contextual Stochastic Optimization (CSO) problem
\begin{equation} \label{eq: cso}
\begin{aligned}
\min_{\mathbf{z} \in \mathcal{Z}} \mathbb{E}_{\tilde{\mathbf{c}} \sim \mathbb{P}(\mathbf{c} \mid \mathbf{s})}[g(\mathbf{z}, \mathbf{\tilde{c}})].
\end{aligned}
\end{equation}

Let
$\tilde{\mathbf{c}}_{v} = (\tilde{c}_{iv})_{i \in \mathcal{I}}$ and $\mathbf{s}_{v} = (\mathbf{s}_{iv})_{i \in \mathcal{I}}$ for each $v \in \mathcal{V}$. Under the assumption that, conditional on its own covariates, each component \(c_v\) is independent across options \(v\in\mathcal V\), one has
\[
\mathbb{P}(\mathbf{c} \mid \mathbf{s}) = \Pi_{v \in \mathcal{V}} \mathbb{P}_v(\mathbf{c}_{v} \mid \mathbf{s}_{v}),
\]
where $\mathbb{P}_v(\cdot \mid \mathbf{s}_{v})$ is the
contextual distribution for option $v$. This assumption is reasonable when each option’s cost depends primarily on its own covariates. Shared features, such as weather, create similarity but do not induce direct dependence. Consequently, Problem (\ref{eq: cso}) can be written as
\begin{equation} \label{eq: cso-de}
\begin{aligned}
\min_{\mathbf{z} \in \mathcal{Z}} \mathbb{E}_{\tilde{\mathbf{c}} \sim \Pi_{v \in \mathcal{V}} \mathbb{P}_v(\mathbf{c}_{v} \mid \mathbf{s}_{v})}[g(\mathbf{z}, \mathbf{\tilde{c}})].
\end{aligned}
\end{equation}

\subsection{Constructing a Contextual Distribution Oracle from Observational Data} \label{sec:learn_from_observational}

The contextual distribution $\mathbb{P}(\mathbf{c} \mid
\mathbf{s})$ is also generally unknown but can often be inferred from historical observations. In many real-world COS applications (including CSOFP) however, data are observational, i.e., only the costs corresponding to decisions actually taken were recorded. As a result, the dataset available to
estimate $\mathbb{P}(\mathbf{c} \mid \mathbf{s})$ takes the form $\mathbb{D} =
\{(\mathbf{s}_o, \mathbf{C}_o) \}_{o = 1}^{|\mathbb{D}|}$, 
where
\(
  \mathbf{s}_o = (s_{iv,o})_{i\in\mathcal I,\;v\in\mathcal V}
  ~\text{and} ~
  \mathbf{C}_o = (\mathbf{c}_{v,o})_{v\in\mathcal V}
\) collect the context and realized costs, with \(\mathbf{c}_{v,o}\) observed only if option \(v\) was selected in request \(o\).  

The
learning task consists in finding a full‑vector forecasting oracle  $\mathcal{M}$ that
approximates $\mathbb{P}(\mathbf{c} \mid \mathbf{s})$ in a
supervised manner using the dataset $\mathbb{D}$: $\mathcal{M}$
receives a context $\mathbf{s}$ as input and returns a probability
distribution $\mathcal{M}(\mathbf{s})$ on the cost vector $\tilde{\mathbf{c}}$. Under the conditional independence assumption from the previous section, learning
$\mathbb{P}(\mathbf{c} \mid \mathbf{s})$ amounts to learning
$\mathbb{P}_v(\mathbf{c}_{v} \mid \mathbf{s}_{v})$ for each selection option
$v$ separately. Partition the full dataset $\mathbb{D}$ into subsets
$\mathbb{D}_v = \{(\mathbf{s}_{v, 1}, \mathbf{c}_{v, 1}), (\mathbf{s}_{v, 2}, \mathbf{c}_{v, 2}), \ldots \},$ representing the dataset for option $v$. A machine learning model $\mathcal{M}_v$ is then trained on $\mathbb{D}_v$ to approximate each $\mathbb{P}_v(\mathbf{c}_{v} \mid \mathbf{s}_{v})$. Assemble the oracle as the product
\(
  \mathcal{M}(\mathbf{s})
  = \bigl(\mathcal{M}_v(\mathbf{s}_v)\bigr)_{v\in\mathcal V},
\)
which generates joint forecasts by treating each component independently.
\revision{
Note that if every \(\mathbf{s}_v\) shares the same structure and dimension, a single learning architecture can serve as $\mathcal{M}_v$ for all $v \in \mathcal{V}$, trained on the union $\cup_{v \in \mathcal{V}}\mathbb{D}_v$.

By reducing the high‑dimensional joint learning task to \(|\mathcal V|\) manageable subproblems, this decomposition turns observational data into a black‑box contextual distribution oracle \(\mathcal M(\mathbf s)\), which will serve as the input for the solution methods developed in Section \ref{sec:c-saa} (C-SAA) and Section \ref{sec:c-ro} (C-RO).}

\subsection{Contextual Sample Average Approximation (C-SAA)}
\label{sec:c-saa}

Problem \eqref{eq: cso-de} is approximated by generalizing the Sample
Average Approximation (SAA) method \citep{kleywegt2002sample} to
contextual stochastic optimization. Given a request with covariates
$\mathbf{s}$, the Contextual Sample
Average Approximation (C-SAA) proceeds in two stages:

\begin{enumerate}
  \item \textbf{Candidate generation.}  
    \begin{enumerate}
      \item Draw \(N_1\) independent cost vector samples 
      \(\{\boldsymbol\xi^n\}_{n=1}^{N_1}\sim\mathcal M(\mathbf s)\), 
      where \(\boldsymbol\xi^n=(\xi^n_{iv})_{i\in\mathcal I,\;v\in\mathcal V}\).  
      \item Solve
      \[
        \mathbf{z}_i \;=\;\arg\min_{\mathbf z\in\mathcal Z}
        \;\frac{1}{N_1}\sum_{n=1}^{N_1} g(\mathbf z,\boldsymbol\xi^n),
        \quad
        i=1,\dots,Q,
      \]
      independently \(Q\) times to obtain candidate solutions 
      \(\{\mathbf z_1,\dots,\mathbf z_Q\}\) with sample‐average objectives 
      \(\{f_1,\dots,f_Q\}\).  
      \item Compute the average objective 
      \(\ubar{f} = \tfrac1Q\sum_{i=1}^Q f_i\), which is a statistical lower bound on the true optimal value \(f^*\) because $\ubar{f}$ is an unbiased estimator of $\mathbb{E}_{{\tilde{\mathbf{c}} \sim
       \mathbb{P}(\mathbf{c} \mid \mathbf{s})}}[\ubar{f}]$ and \(\mathbb{E}_{{\tilde{\mathbf{c}} \sim
       \mathbb{P}(\mathbf{c} \mid \mathbf{s})}}[\ubar{f}]\le f^*\).
    \end{enumerate}

  \item \textbf{Candidate evaluation.}  
    \begin{enumerate}
      \item For each \(\mathbf z_i\), draw a larger set of evaluation samples
      \(\{\boldsymbol\xi^n\}_{n=1}^{N_2}\sim\mathcal M(\mathbf s)\), with \(N_2\gg N_1\).
      \item Compute the re‐estimated objectives 
      \(\bar f(\mathbf z_i)=\tfrac1{N_2}\sum_{n=1}^{N_2}g(\mathbf z_i,\boldsymbol\xi^n)\).  
      \item Each \(\bar f(\mathbf z_i)\) is an upper bound on \(f^*\) since each candidate solution is feasible; the difference \(\bar f(\mathbf z_i)- \ubar{f}\) approximates the optimality gap.  
      \item Select the candidate \(\mathbf z^*\) with the smallest estimated gap.
    \end{enumerate}
\end{enumerate}

\smallskip
\noindent\textbf{Discussion.}  Unlike the reweighting approach of \citet{bertsimas2020predictive}, which leverages complete historical samples for each fully-observed cost vector, C–SAA relies solely on observational data (i.e., only the costs of the selected options are available).  Attempting to form all possible cost vector combinations from partial observations would lead to a combinatorial explosion as \(|\mathcal V|\) and \(|\mathbb{D}|\) grow, making such method intractable in large, high‐dimensional selection spaces.

\subsection{Contextual Robust Optimization (C-RO)} \label{sec:c-ro}
The performance of C-SAA largely hinges on the accuracy of the
estimated contextual distribution. Robust optimization provides a hedge against misspecification of the contextual distribution by guarding against worst–case scenarios, making it a viable approach for
contextual stochastic optimization.  Let $\mathcal{U}(\mathbf{s}) \subset \mathbb R^{|\mathcal I|\times|\mathcal V|}$
represent a contextual uncertainty set for the cost vector
$\tilde{\mathbf{c}}$. The Contextual Robust Optimization (C-RO) problem is
\begin{equation} \label{eq: C-RO}
    \begin{aligned}
        \min_{\mathbf{z} \in \mathcal{Z}}  \max_{\bm{\xi} \in \mathcal{U}(\mathbf{s})} g(\mathbf{z}, \bm{\xi}).
    \end{aligned}
\end{equation}

Two constructions of \(\mathcal U(\mathbf s)\) are considered: a discrete sampling–based set \(\mathcal U_d(\mathbf s)\) and a budgeted interval set \(\mathcal U_b(\mathbf s, B)\).

\subsubsection{Contextual Discrete Uncertainty Set}

Problem
(\ref{eq: C-RO}) can be reformulated into a single-level problem by
enumerating every possible realization of the random vector $\tilde{\mathbf{c}}$, a standard technique in robust
optimization. In principle, a discrete distribution can be generated as mentioned before, by combining all combinations of the
observed data points. However, this will make the problem intractable as discussed.

\revision{
To retain tractability, a sampling-based discrete uncertainty set in the same spirit of C-SAA is introduced. Let $\mathcal{U}_d(\mathbf{s}) = \{\bm{\xi}^1, \ldots, \bm{\xi}^N\}$ denote a collection of \(N\) deviation scenarios drawn from \(\mathcal M(\mathbf s)\). The robust counterpart of \eqref{eq: C-RO} becomes
\begin{equation} \label{eq: C-RO-D}
    \begin{aligned}
         \text{(C-RO-D)} \quad \min_{\mathbf{z} \in \mathcal{Z},\;t\in\mathbb R} \quad &  t\\
          \text{s.t.} \quad & t \geq g(\mathbf{z}, \bm{\xi}^n),
    \quad n=1,\dots,N.
    \end{aligned}
\end{equation}
}
\subsubsection{Contextual Budgeted Interval Uncertainty Set}

An alternative uncertainty set is the budgeted interval set
$\mathcal{U}_b(\mathbf{s}, B)$ \citep{bertsimas2004price}, which is defined by a robust budget $B$. 
This budget limits the extent to
which limits how many components of the uncertain parameter $\tilde{\mathbf{c}}$ can simultaneously deviate to their upper bounds. 

Rather than estimating the full joint distribution, only the conditional lower and upper quantiles at levels \(q_{\rm low}\) and \(q_{\rm up}\) are needed to define a \((q_{\rm low},q_{\rm up})\) prediction interval for each option.

For each option \(v\in\mathcal V\), let
\(
\underline{\boldsymbol\xi}_v 
= \bigl(\underline\xi_{iv}\bigr)_{i\in\mathcal I}
\) and
\(
\overline{\boldsymbol\xi}_v 
= \bigl(\overline\xi_{iv}\bigr)_{i\in\mathcal I}
\)
denote the estimated conditional \(q_{\rm low}\)- and \(q_{\rm up}\)-quantiles extracted from \(\mathcal M_v(\mathbf s_v)\):
\[
\underline{\bm \xi}_{v}
= \inf\{\mathbf  u:\mathbb P_v(\mathbf c_{v}\le \mathbf u\mid\mathbf s_v)\ge q_{\rm low}\}, 
\quad
\overline{\bm \xi}_{v}
= \inf\{\mathbf u:\mathbb P_v(\mathbf c_{v}\le \mathbf u\mid\mathbf s_v)\ge q_{\rm up}\}.
\]
Given a budget \(B\), the budgeted interval uncertainty set is
\[
\mathcal U_b(\mathbf s,B)
= \Bigl\{\boldsymbol\xi:\;
\xi_{iv}
= \underline\xi_{iv}
+ \delta_{iv}\,(\overline\xi_{iv}-\underline\xi_{iv}),\;
\sum_{i,v}\delta_{iv}\le B,\;
\delta_{iv}\in[0,1]\Bigr\}.
\]



With a fixed $\mathbf{z}$, the inner maximization problem in \eqref{eq: C-RO} is 
    \begin{align}
        \max_{\delta_{iv}} \quad &  \sum_{i, v} (\overline\xi_{iv}-\underline\xi_{iv})  z_{iv} \delta_{iv} \\
        \text{s.t.} \quad &  \sum_{i,v} \delta_{iv} \leq B \label{eq: imp:c1}\\
        & \delta_{iv} \in [0, 1], \quad \forall i, v \label{eq: imp:c2}.
    \end{align}

Let \(y_j \in \{0,1\}\) and \(w_{iv}^j \ge 0\) be auxiliary variables.  Denote by 
\(\pi \ge 0\) and \(\lambda_{iv}\ge 0\) the dual variables corresponding to Constraints (\ref{eq: imp:c1}) and (\ref{eq: imp:c2}), respectively.  Finally, let \(M\) be a sufficiently large constant. \revision{A valid choice for the ``big‑M” constant is
\(
M \;=\;\max_{i, v, j}
\bigl\{\,b_{iv}^n-b_{iv}^1,\;b_{iv}^{\,n-1},\;(\underline{\xi}_{iv}+a_{iv}^j)\,b_{iv}^n\bigr\}.
\)}
By taking the dual of the inner maximization problem and introducing \(y_j\) and \(w_{iv}^j\), a robust counterpart of (\ref{eq: C-RO}) is obtained as
\begin{equation}\label{eq: RCP}
\begin{aligned} 
    \text{(C-RO-B)} \quad \min_{\mathbf{z} \in \mathcal{Z}, y,w,\pi,\lambda}  \quad & (\sum_{i, v, j} w_{iv}^j  + \pi B + \sum_{i, v} \lambda_{iv})\\
    \text{s.t.} \quad & \pi + \lambda_{iv} \geq (\overline\xi_{iv}-\underline\xi_{iv}) z_{iv},  \quad \forall i, v \\
    & z_{iv} \leq b^j_{iv} + M (1-y_j),  \quad \forall i, v, j \\
    & z_{iv} \geq b^{j-1}_{iv} - M (1-y_j), \quad \forall i, v, j \\
    & w_{iv}^j \geq (\underline{\xi}_{iv}+a_{iv}^j)z_{iv} - M(1-y_j), \quad \forall i, v, j \\
    & y_j \in \{0, 1\}, ~w_{iv}^j \geq 0, ~\pi \geq 0, ~\lambda_{iv} \geq 0. \\
\end{aligned}
\end{equation}

\revision{
\subsection{Performance Guarantees on C-SAA \& C-RO}
C–SAA replaces the unknown cost distribution with Monte-Carlo samples generated by a learned contextual oracle and minimizes the resulting empirical objective.  
Under standard regularity conditions, the SAA objective is statistically consistent: its optimal value converges exponentially fast to the true optimum as $N_{1}$ grows \citep{kleywegt2002sample,shapiro2021lectures}.  
When the oracle is well-calibrated and $N_{1}$ is large, C–SAA thus attains the smallest \emph{expected} cost.  
However, it remains sensitive to distributional misspecification and finite-sample error, which can expose the solution to rare but severe tail events.

C–RO mitigates this risk by optimizing the worst-case cost over a prescribed uncertainty set. Let a feasible solution $\mathbf z$ satisfy
\(
\max_{\bm\xi \in \mathcal U(\mathbf s)} g(\mathbf z,\bm\xi) \le t.
\)
For a discrete uncertainty set built from $N$ samples (C-RO-D), one obtains with probability $1-\delta$ that
\[
\Pr\bigl(g(\mathbf z,\tilde{\mathbf{c}})>t\bigr)\;\le\;\epsilon
\]
provided
\[
N \;\ge\; 2|\mathcal V||\mathcal I| 
      \;+\;\frac{2|\mathcal V||\mathcal I|}{\epsilon}\ln\!\Bigl(\frac{2}{\epsilon}\Bigr)
      \;+\;\frac{2}{\epsilon}\ln\!\Bigl(\frac{1}{\delta}\Bigr)
\]
samples are drawn \citep{calafiore2006scenario, campi2008exact, bertsimas2021probabilistic}.  
Although this bound grows linearly in the dimension of the uncertainty $|\mathcal V||\mathcal I|$, it too can degrade under distributional misspecification.  
Under a budgeted-interval uncertainty set (C-RO-B), the violation probability instead satisfies
\[
\Pr\bigl(g(\mathbf z,\tilde{\mathbf{c}})>t\bigr)\;\le\;\exp\!\Bigl(-\frac{B^{2}}{2\,|\mathcal V||\mathcal I|}\Bigr),
\]
so that the budget parameter $B$ directly controls the trade-off between average-case performance and robustness \citep{bertsimas2004price}.  
This tunable flexibility makes C–RO-B an attractive middle ground between C–SAA and C–RO-D.

Extending these guarantees to bound worst-case regret or post-decision surprise under partial-feedback settings---where counterfactual costs are unobserved---remains an important direction for future work.
}

\subsection{\revision{Application to the Omnichannel Multi-Courier CSOFP}}\label{sec:specialize}

The CSO framework applies directly to CSOFP via the following identifications:
\[
\begin{aligned}
  &\mathcal I \;\mapsto\;\{\text{SKUs}\}, 
  &&\mathcal V \;\mapsto\;\mathcal K\times\mathcal L,\\
  &z_{i v}\;\mapsto\;z_{ik\ell}, 
  &&\mathbf s_v\;\mapsto\;\mathbf s_{k\ell},\\
  &\tilde{\mathbf c}=(\tilde c_v)_{v\in\mathcal V}
  \;\mapsto\;\tilde{\mathbf d}=(\tilde d_{k\ell})_{k\in\mathcal K,\ell\in\mathcal L}.
  && \mathbb{P}_v \;\mapsto\; \mathbb{P}_{k \ell}.
\end{aligned}
\]

Under the conditional independence assumption, \[
  \mathbb P({\mathbf d}\mid \mathbf s)
  = \prod_{k\in\mathcal K}\prod_{\ell\in\mathcal L}
    \mathbb P_{k\ell}(d_{k\ell}\mid \mathbf s_{k\ell}).
\]
Hence CSOFP (Problem \ref{eq: CSOFP}) can be rewritten as
\begin{equation} \label{eq: CSOFP2}
\begin{aligned}
\min_{\mathbf z\in\mathcal Z} \mathbb{E}_{\tilde{\mathbf{d}} \sim \Pi_{k, \ell} \mathbb{P}_{k \ell}(d_{k \ell} \mid \mathbf{s}_{k \ell})}[g(\mathbf{z},  \tilde{\mathbf{d}})].
\end{aligned}
\end{equation}

\revision{
To learn each marginal distribution \(\mathbb P_{k\ell}(d_{k\ell}\mid\mathbf s_{k\ell})\), the historical dataset is first partitioned into
\(
\mathbb D_{k\ell}
=\bigl\{(\mathbf s_{k\ell,m},\,d_{k\ell,m})\bigr\}_{m=1}^{M_{k\ell}}.
\)
This is followed by training a per‑carrier model \(\mathcal M_k\) on the dataset
\(\mathbb D_k = \bigcup_{\ell\in\mathcal L}\mathbb D_{k\ell}\).  This per-carrier approach is a practical adaptation of the generic per-option framework from Section~\ref{sec:learn_from_observational}. Rather than training a distinct model for each individual location-carrier pair (which can be inefficient and suffer from data sparsity), the data is pooled to train a single, more robust model for each carrier. This decomposition simplifies learning, allows the model to generalize to new locations, and was found to outperform a single global model in preliminary experiments.  The full contextual distribution oracle is then assembled as
\(
\mathcal M(\mathbf s)
= \bigl(\mathcal M_{k}(\mathbf s_{k\ell})\bigr)_{k\in\mathcal K,\,\ell\in\mathcal L}.
\)
}

\subsubsection{Applying C-SAA and C-RO}
\revision{
Computing the expectation in
\(\mathbb{E}_{\tilde{\mathbf d}\sim\mathcal M(\mathbf s)}\bigl[g(\mathbf z,\tilde{\mathbf d})\bigr]\)
directly would require either enumerating all \(C^{|\mathcal K|\times|\mathcal L|}\) deviation vectors or numerically integrating a continuous CDF—both computationally intractable.  Instead, the two solution paradigms use the contextual distribution oracle \(\mathcal M(\mathbf s)\) as follows:
\begin{itemize}

    \item  C-SAA replaces the expectation by a sample average: draw \(N_1\) full‑vector scenarios \(\{\tilde{\mathbf d}^n\}\sim\mathcal M(\mathbf s)\) and solve  
  \[
    \min_{\mathbf z\in\mathcal Z}\;\frac{1}{N_1}\sum_{n=1}^{N_1} g\bigl(\mathbf z,\tilde{\mathbf d}^n\bigr).
  \]
    \item  C-RO builds an uncertainty set from oracle outputs—either the discrete scenario set \(\{\tilde{\mathbf d}^n\}_{n \in [N]}\) or a budgeted‑interval set from marginal quantiles—and then solves  
  \[
    \min_{\mathbf z\in\mathcal Z}\;\max_{\mathbf d\in\mathcal U(\mathbf s)}g(\mathbf z,\mathbf d).
  \]
\end{itemize}

Both approaches avoid the exponential support growth while remaining agnostic to the internal structure of \(\mathcal M\).  A full C‑RO‑B reformulation is provided in the Appendix F.
}

\section{Learning the Contextual Distribution of Delivery Deviations for CSOFP} \label{sec: learning-app}
Having defined in Section \ref{sec:methodology} a generic contextual distribution oracle \(\mathcal M(\mathbf s)\), this section describes its concrete implementation for CSOFP.  It begins by selecting features to construct the
covariate vector. Two classes of machine learning
methods are then proposed to estimate the contextual distribution of deviations: (1) Probabilistic multi-class classification, which
leverages the discrete support of deviations; and (2) \revision{Tree-based quantile regression}, which directly estimates the conditional CDF and avoids quantile crossing \citep{salari2022real}. 

For simplicity, the carrier index $k$ is omitted
in the notation, and the dataset for training the per-carrier model is written as $\mathbb{D} = \{ (\mathbf{s}_m, d_m) \}_{m \in [M]}$.

\subsection{Contextual Feature Selection}

Carefully selecting contextual features is crucial for constructing
informed covariates for delivery deviation predictions, as various stages in the
order life cycle can contribute to deviations. 

Candidate features fall into three categories:
\begin{enumerate}

\item {\em Order-level}: This includes customer
  location, SKU weight and dimensions, SKU type, quantities, release
  hour, release day, and planned lead time (i.e., the difference
  between the promised delivery date and the order release
  time). These factors can influence deviations through logistical
  complexity, handling requirements, and time-sensitive elements
  affecting the order's processing and delivery schedule.

\item {\em Fulfillment center}: This includes coordinates, location
  type, on-hand capacity, and on-hand inventory. The geographical
  location and type of fulfillment center can affect the efficiency of
  order processing, while capacity and inventory levels determine the
  speed and reliability of fulfilling orders.

\item {\em Carrier-related}: This includes
  service level, promised transit time, shipping charges, and carrier
  zone. These factors impact the carrier ability to meet delivery
  deadlines, with service level and transit times directly affecting
  deviations, while shipping charges and zones reflect logistical
  challenges and potential bottlenecks.
\end{enumerate}

\subsection{Probabilistic Multi-Class Classification} \label{sec:mc-clf}

This discrete nature allows the contextual distribution
learning problem to be framed as a multi-class classification (MC-CLF)
problem. Under MC-CLF, each class represents a possible delivery deviation value in number of days,
and the output class probabilities correspond directly to the
probability mass function of the discrete distribution. For a new
sample with covariates vector $\mathbf{s}$, let the conditional
probability of class $c$ be $p_c(\mathbf{s}) = \mathbb{P} (d = \xi_c \mid
\mathbf{s})$, with the corresponding predicted probability from MC-CLF
denoted as $\hat{p}_c(\mathbf{s})$.

Multinomial logistic regression (MLR) and classification tree-based
models are considered to solve MC-CLF. Specifically, MLR is a classic
classifier that aims to estimate the class probabilities by maximizing
the log-likelihood function, the negative of the log loss
function. The log loss is defined as $-\sum_{m = 1}^ {M} \sum_{c = 1}^C \mathbf{1}\{d_m = \xi_c\} \log (\hat{p}_c(\mathbf{s}_m)).$
On the other hand, classification trees use a nonparametric approach
that recursively partitions the data into subsets based on
feature values in a top-down manner, resulting in a tree-like model of
decisions. 

Let $\mathcal{R}_1, \ldots, \mathcal{R}_R$ be the regions partitioned by a classification tree. Classification tree estimates
\begin{equation}
    \hat{p}_c(\mathbf{s}) = \frac{\sum_{r = 1}^R \mathbf{1} \{\mathbf{s} \in \mathcal{R}_r\} \cdot \mathbf{1} \{d = \xi_c\}}{\sum_{r = 1}^R \mathbf{1} \{\mathbf{s} \in \mathcal{R}_r\}},
\end{equation}

To enhance predictive power, ensemble learning methods such as bagging
and gradient boosting are employed, through the use of random forests
(RF) and gradient boosted trees (GBT). The details for their calculation of the class probabilities are available in Appendix G.

The performance of the MC-CLF models are evaluated by log loss and Brier score. The latter is
defined as the average squared error between predicted probabilities
and actual outcomes $\sum_{m = 1}^M \sum_{c = 1}^C (\mathbf{1}\{d_m = \xi_c\} \hat{p}_c(\mathbf{s}_m))^2$. For both metrics, lower values indicate better probabilistic prediction performance. 

\revision{
\subsubsection{Probability Calibration}

Raw class‐probability estimates \(\hat p_c(\mathbf s)\) from classifiers often exhibit mis‐calibration: a predicted probability \(\alpha\) does not correspond to an empirical frequency \(\alpha\) \citep{niculescu2005predicting}.  \emph{Good calibration} means that, over all samples with \(\hat p_c(\mathbf s)=\alpha\), the true fraction belonging to class \(c\) is approximately \(\alpha\).  Formally, let 
\(
  f(\alpha)
  = \mathbb{P}\bigl(d=\xi_c \mid \hat p_c(\mathbf s)=\alpha\bigr)
\)
be the calibration function. Perfect calibration means \(f(\alpha)=\alpha\) for all \(\alpha\in[0,1]\). 

Isotonic regression \citep{zadrozny2002transforming} learns a non‐decreasing mapping \(g:[0,1]\to[0,1]\) that best aligns raw scores to observed frequencies, preserving the rank order of \(\hat p_c\).  Given a calibration set \(\mathbb{D}_{cal} = \{(\hat p_c(\mathbf s_m),\,y_{c, m})\}_{m = 1}^{M_{cal}}\) with \(y_{c, m}=\mathbf{1}\{d_m=\xi_c\}\), isotonic regression solves
\[
  \min_{g}\;\sum_{m=1}^{M_{cal}}\bigl(y_{c, m} - g(\hat p_c(\mathbf s_m))\bigr)^2
  \quad\text{s.t.}\quad g(\alpha)\le g(\alpha')\;\text{whenever}\;\alpha\le\alpha'.
\]
The calibrated probabilities are \(\tilde p_c(\mathbf s)=g(\hat p_c(\mathbf s))\).  By enforcing monotonicity, isotonic calibration ensures that higher raw scores remain higher after mapping, and empirical frequencies within each score bin match the calibrated values.

In this paper, a $k$-fold cross‐validation procedure is employed: for each fold $i$, the classifier is trained on its training split, raw scores are obtained on the held‐out calibration split, and an isotonic regression model $g^{(i)}$ is fit to those scores. On unseen test instance with covariates $\mathbf{s}$, each fold’s classifier and corresponding $g^{(i)}$ yield calibrated probabilities $\tilde p_c^{(i)}(\mathbf{s})$, which are then averaged across all $k$ folds to produce the final estimate
\(
\tilde p_c(\mathbf{s}) = \frac{1}{k} \sum_{i=1}^k \tilde p_c^{(i)}(\mathbf{s})\,.
\)
}
\revision{Since per-class calibration does not guarantee the final probabilities will sum to one, a normalization step is applied. The vector of calibrated class probabilities $\bigl(\tilde p_c(\mathbf{s})\bigr)_{c=1}^C$ is divided by its sum to yield a valid multi-class probability distribution.}

\subsection{Tree-Based Quantile Regression}
\revision{
Unlike multi‑class classification methods that treat each deviation value as an unordered category, quantile regression (QR) exploits the inherent ordering of deviations by framing the task as ordinal regression and directly estimating the conditional cumulative distribution function (CDF) 
$F(u \mid \mathbf{s}) = \mathbb{P}(d \leq u \mid \mathbf{s})$.

Standard
QR models minimizes, for each quantile level $q$, the pinball loss 
\begin{equation}
\rho_q(d, \hat{d}) \;=\; q\,\max(d-\hat{d},0) \;+\; (1-q)\,\max(\hat{d}-d,0),
\end{equation}
where \(\hat d\) denotes the model’s predicted \(q\)‑th quantile.
 Collecting these quantile estimates across \(q\) yields a stepwise approximation of the CDF. However, estimating each quantile independently can violate the ordering
\(\hat d_{(q_1)}\le \hat d_{(q_2)}\) for \(q_1<q_2\), resulting in a non‑monotonic—and thus invalid—CDF.

Tree‐based QR methods---namely regression trees extended to CDF estimation, and their ensemble version, quantile random forests (QRF) \citep{meinshausen2006quantile}---avoid this “quantile crossing” problem by extracting quantiles from the empirical distribution of training sample in each leaf region.}

A regression tree partitions the
covariate space into disjoint regions to minimize mean squared error (MSE). Although regression trees are often used for point forecasts (e.g., predicting the mean response in each leaf), they can
be adapted for distributional forecasting by preserving all
observations within each partitioned region. For a new sample with covariates
$\mathbf{s}$, let $\mathcal{R}(\mathbf{s})$ be the leaf region into which $\mathbf{s}$ falls. The conditional CDF of the deviation can then be approximated by the empirical distribution of the training samples in that leaf:
\begin{equation}
    \hat{F}(u \mid \mathbf{s}) = \sum_{m = 1}^M \frac{\mathbf{1} \{\mathbf{s}_m \in \mathcal{R}(\mathbf{s})\} }{ |\{j: \mathbf{s}_j \in \mathcal{R}(\mathbf{s})\}|} \mathbf{1}\{d_m \leq u\}.
\end{equation} 
A subsequent kernel density estimation applied to \(\{d_m:\mathbf{s}_m\in\mathcal{R}(\mathbf{s})\}\) can further smooth this step‐function CDF \citep{salari2022real}.

QRF extends this approach by pooling the leaf‐level samples across an ensemble of $T$ trees. For each tree $t$, let $\mathcal{R}^{(t)}(\mathbf{s})$ be leaf region containing $\mathbf{s}.$ The QRF estimate of the conditional CDF is
\begin{equation}
    \hat{F}(u \mid \mathbf{s}) = \frac{1}{T} \sum_{t= 1}^T  \sum_{m = 1}^M \frac{\mathbf{1} \{\mathbf{s}_m \in \mathcal{R}^{(t)}(\mathbf{s})\} }{ |\{j: \mathbf{s}_j \in \mathcal{R}^{(t)}(\mathbf{s})\}|} \mathbf{1}\{d_m \leq u\}.
\end{equation}
\revision{Because this procedure pools ordered observations directly, it naturally preserves the ordinal support and never produces crossing quantiles.}

To measure the distributional forecasting performance, the paper uses two metrics---the Continuous Ranked Probability Score (CRPS) and pinball loss. CRPS evaluates the closeness of the estimated CDF $\hat{F}$ given covariates $\mathbf{s}$ to the observed ground truth value $d$
\begin{equation}
\text{CRPS}(\hat{F}, d, \mathbf{s}) = \int_{-\infty}^{\infty} (\hat{F}(u \mid \mathbf{s}) - \mathbf{1}\{u \ge d\})^2 du.
\end{equation}
Let $\xi$ and $\xi'$ be two independent random variables with distribution $\hat{F}(u \mid \mathbf{s})$. The CRPS can alternatively be expressed as
\begin{equation} \label{eq: CRPS}
\text{CRPS}(\hat{F}, d, \mathbf{s}) = \mathbb{E}(|\xi - d|) - \frac{1}{2} \mathbb{E}(|\xi - \xi'|).
\end{equation}
 This equation allows for an approximation of CRPS by sampling from $\hat{F}(u \mid \mathbf{s})$. Pinball loss measures the accuracy of quantile predictions by penalizing errors asymmetrically, depending on whether the predictions are overestimates or underestimates.

\section{Case Study}
\label{sec: experiments}

This section presents computational results that demonstrate the
potential of the proposed CSO framework when evaluated on real-world
industrial instances.

\subsection{Data Overview}

The industrial partner provides a month-long dataset in 2023, covering
more than 1 million online orders for home delivery and over 20,000
unique SKUs. Each order can be sourced from roughly 250 to 1,000 possible fulfillment centers across the
company's network and shipped via one of 10 to 30 distinct carriers.

The target variable in the machine learning problem, delivery deviation, is defined as the difference between the actual
delivered date and the promised delivery date. To enhance the quality of the analysis, order lines with extreme deviations were excluded---specifically, those with absolute values greater than 10 days. Such outliers often reflect exceptional circumstances, including returns to sender, lost shipments, or atypical business orders, which can introduce noise and bias into the model.

\subsection{Evaluation of Contextual Distributional Learning of Delivery Deviations}
\label{Evaluation of Contextual Distributional Learning of Delivery Deviations}
This section evaluates the contextual distributional learning models presented in Section \ref{sec: learning-app} for the case study.

\subsubsection{Machine Learning Model Settings}
For the evaluation of the learning models, the dataset was partitioned into a training set comprising orders from the first three weeks (approximately 1.3 million order lines) and a test set consisting of orders from the final week (approximately 0.4 million order lines). For multi-SKU orders, features were aggregated by retaining only those of the SKU with the highest shipping charge and quantity. The MC-CLF models were implemented using scikit-learn
\citep{scikit-learn}
and CatBoost \citep{CatBoost}, with probability calibration applied to the best‐performing uncalibrated models based on cross‐validated metrics. Quantile regression models, including
regression tree and QRF, were implemented
via the quantile-forest package \citep{Johnson2024}. All hyperparameters were tuned using three‑fold cross‑validation in Optuna \citep{akiba2019optuna}: MC‑CLF models were tuned to minimize log loss, whereas quantile regression models were tuned to minimize CRPS.  

\subsubsection{Learning Performance Analysis}
In Table \ref{tab:log_loss_brier_score}, two standard metrics---log loss and
Brier score---are reported to measure the accuracy of probabilistic
predictions for each tuned MC-CLF model on the test
set. Calibration yields a clear improvement: the calibrated classifier achieves the lowest log loss (0.965) and Brier score (0.031), down from 1.050 and 0.034 for the uncalibrated gradient‑boosted tree (GBT‑CLF), and substantially outperforming both the random‑forest classifier (RF‑CLF; 1.593, 0.054) and the simpler baselines (MLR and classification tree).

\begin{figure}[!ht]
\centering
    \begin{subfigure}[b]{0.32\textwidth}
        \centering
        \includegraphics[width=\textwidth]{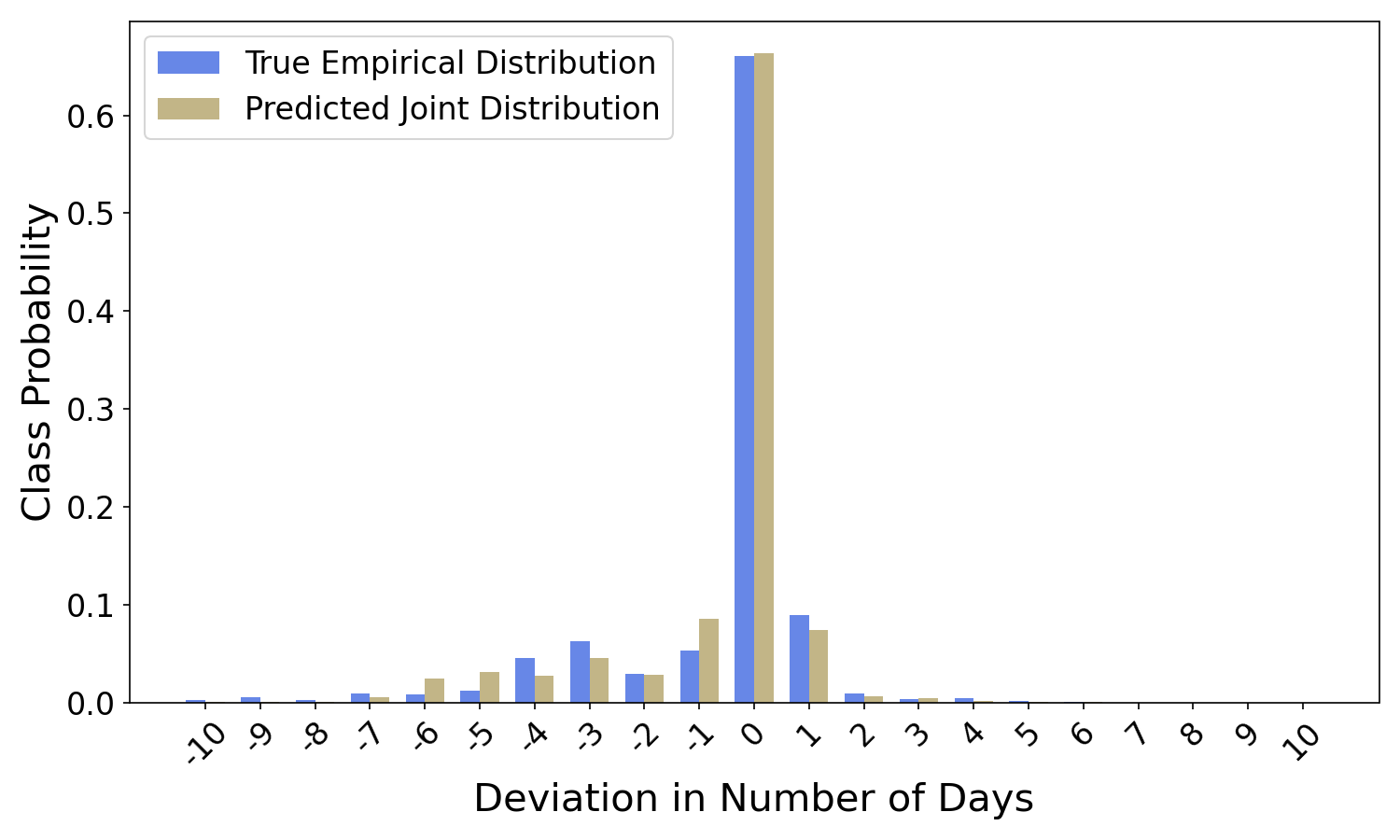}
        \caption{Calibrated-CLF}
        \label{fig:clf_dist}
    \end{subfigure}
    \begin{subfigure}[b]{0.32\textwidth}
        \centering
        \includegraphics[width=\textwidth]{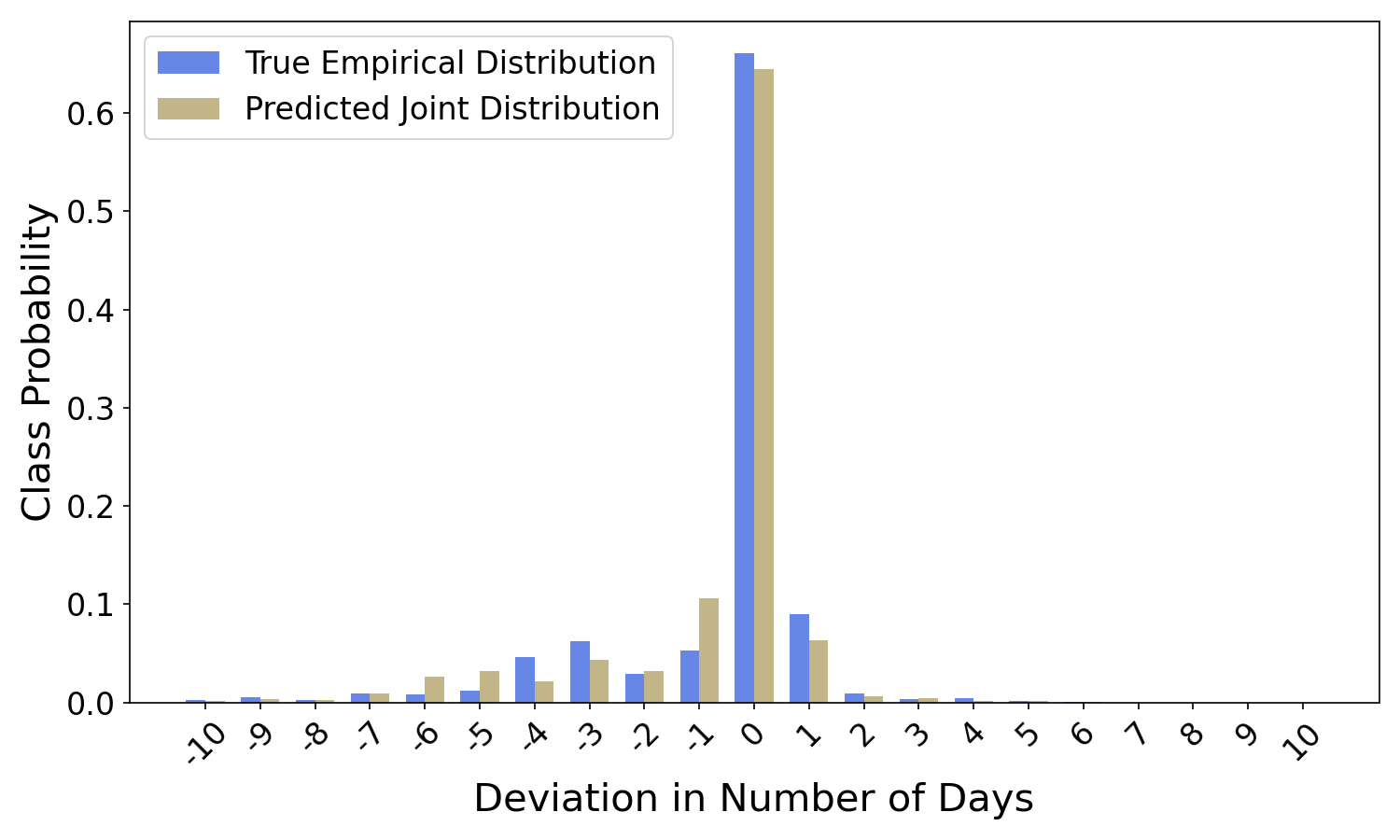}
        \caption{QRF}
        \label{fig:qrf_dist}
    \end{subfigure}
    \begin{subfigure}[b]{0.32\textwidth}
        \centering
        \includegraphics[width=\textwidth]{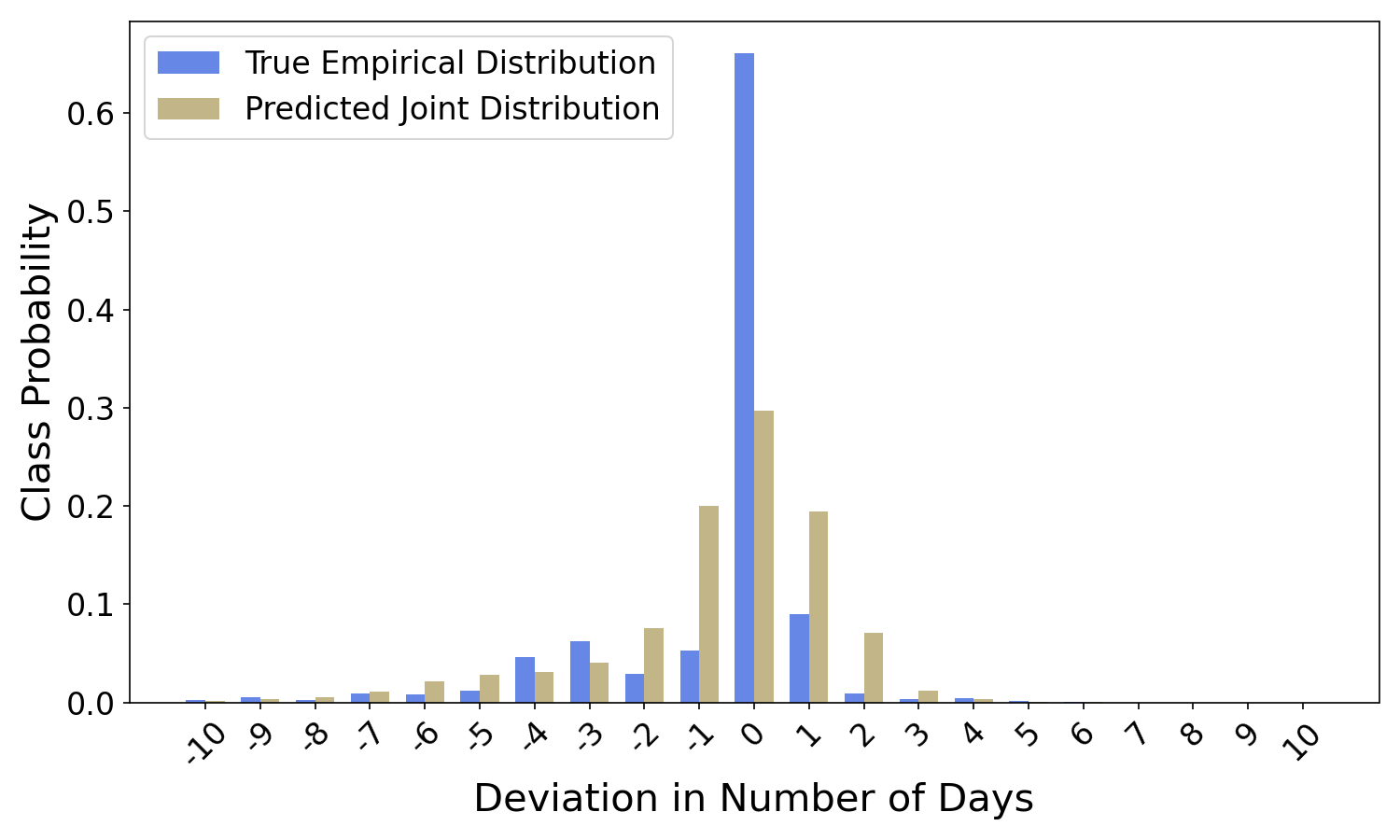}
        \caption{Regression Tree (discretized)}
        \label{fig:rt_dist}
    \end{subfigure}    
    \caption{\revision{Out‐of‐Sample Comparison of the Joint Predicted Distribution and True Empirical Distribution for an Example Carrier.}}
    \label{fig:dist_compare}
\end{figure}

\revision{
Figure \ref{fig:dist_compare} presents the joint predicted distributions of (a) the calibrated‑CLF, (b) the QRF, and (c) the discretized regression tree, each plotted against the true empirical distribution for orders fulfilled by a representative carrier. The calibrated‑CLF closely reproduces the empirical mass at zero deviation (on‑time deliveries) and tightly follows the tails corresponding to early and late deliveries. The QRF exhibits slightly greater discrepancies, particularly in the majority classes (-1 to 1), while the discretized regression tree shows the largest misalignment from the empirical distribution.
}

\begin{table}[!ht]
\centering
\caption{Out-of-Sample Average Log Loss and Brier Score Across All Carriers.}
\begin{tabular}{lccccc}
\toprule
Model & MLR & Classification Tree & RF-CLF & GBT-CLF & Calibrated-CLF \\
\midrule
Log loss & 1.950 & 2.242 & 1.593 & 1.050 & \textbf{0.965} \\
Brier score & 0.073 & 0.072 & 0.054 & 0.034 & \textbf{0.031} \\
\bottomrule
\end{tabular}
\label{tab:log_loss_brier_score}
\end{table}

Table \ref{tab:quantile_metrics} reports pinball losses at selected quantiles and CRPS values for the tuned quantile regression and calibrated-CLF models. CRPS was approximated by drawing 1,000 samples to compute Equation (\ref{eq: CRPS}). The calibrated-CLF model exhibit higher pinball losses than the two quantile regression models, reflecting its inherently discrete output, but nonetheless achieves the lowest CRPS overall, even though it was not directly tuned for this metric. 
\revision{This superior CRPS performance likely results from calibration producing sharper probability estimates, especially for the majority deviation classes.}

\begin{figure}[!ht]
\centering
    \begin{subfigure}[b]{0.32\textwidth}
        \centering
        \includegraphics[width=\textwidth]{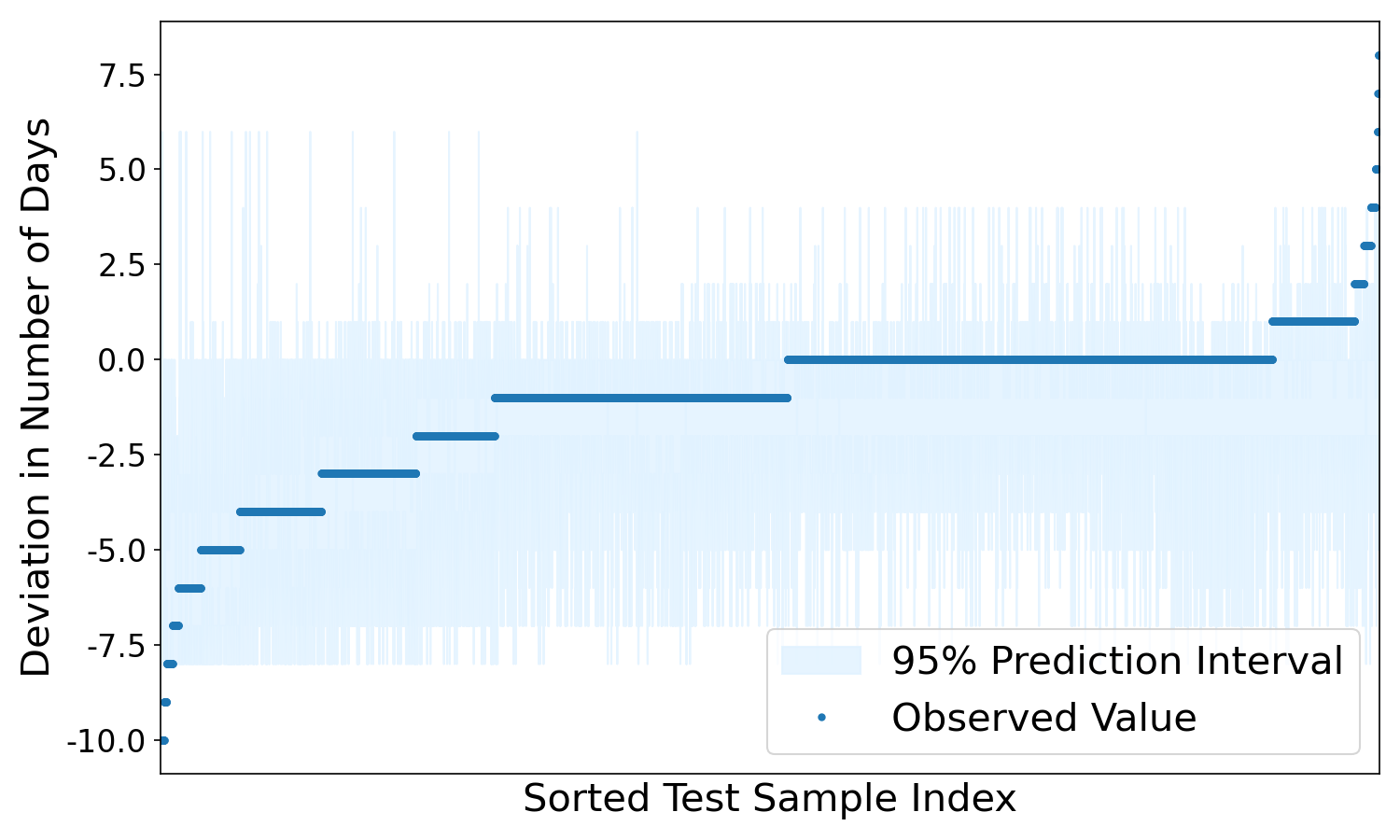}
        \caption{Calibrated-CLF Prediction Intervals}
        \label{fig:mc_clf_interval}
    \end{subfigure}
    \begin{subfigure}[b]{0.32\textwidth}
        \centering
        \includegraphics[width=\textwidth]{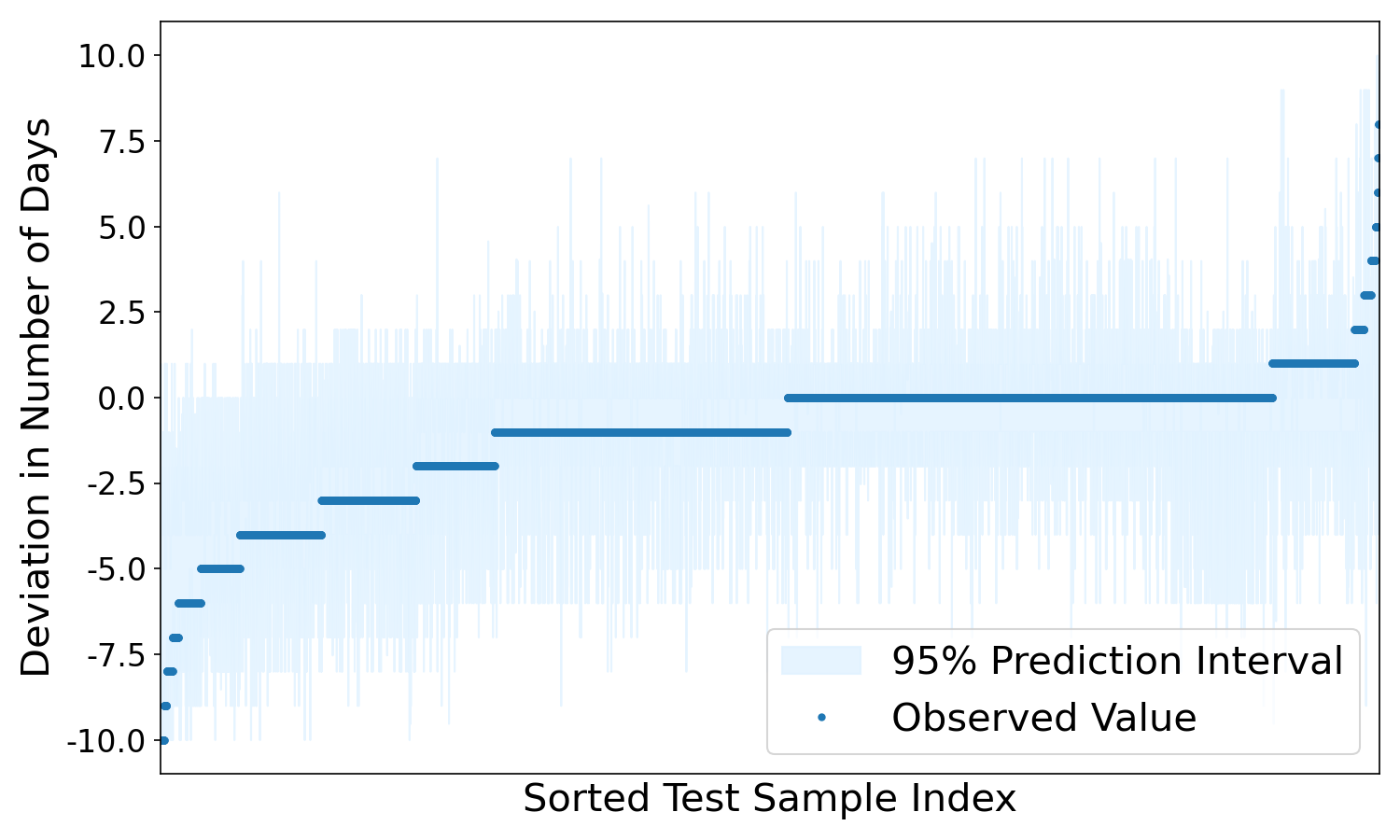}
        \caption{QRF Prediction Intervals}
        \label{fig:qrf_interval}
    \end{subfigure}
    \begin{subfigure}[b]{0.32\textwidth}
        \centering
        \includegraphics[width=\textwidth]{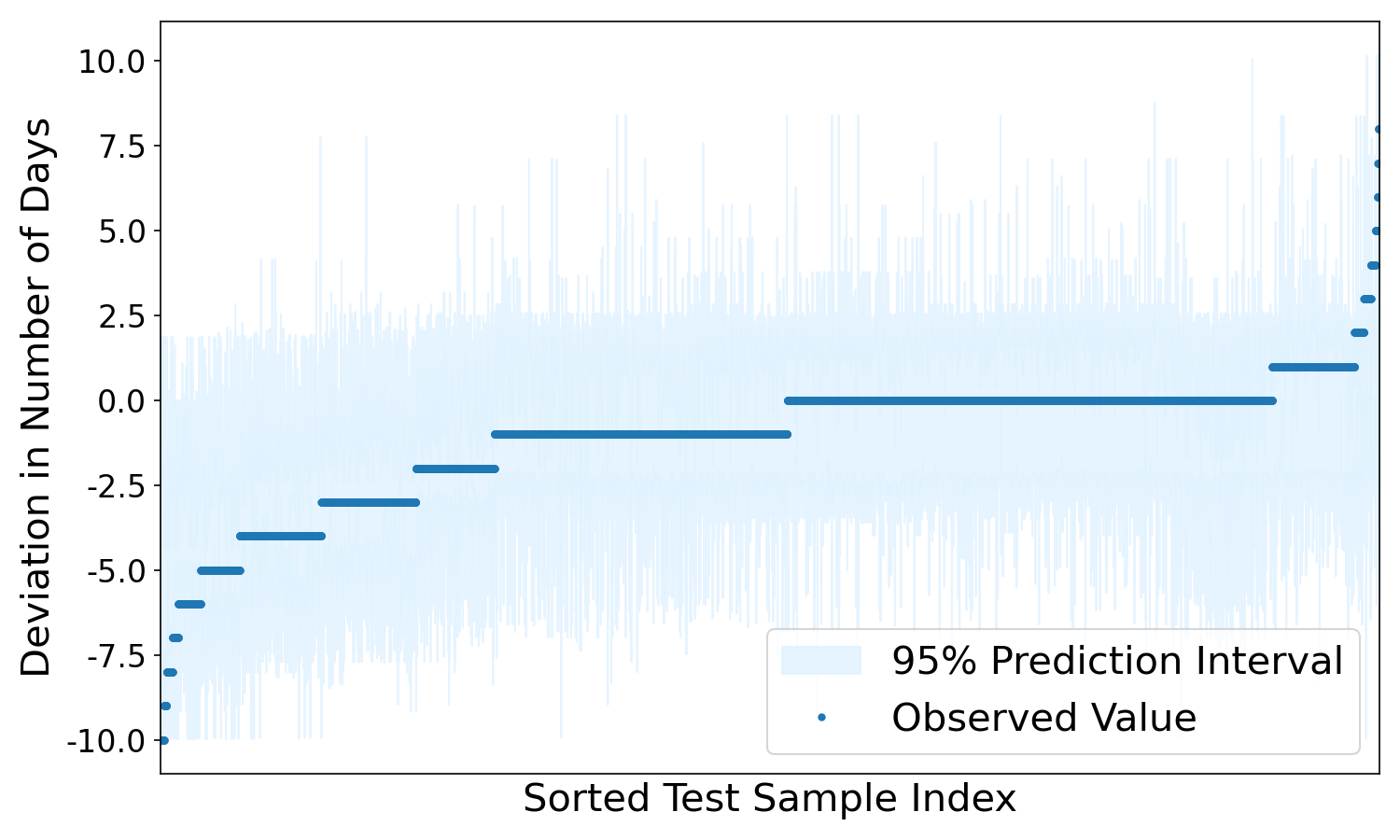}
        \caption{Regression Tree Prediction Intervals}
        \label{fig:rt_interval}
    \end{subfigure}
    
    \caption{Sorted Observed Deviation Values versus 95\% Prediction Intervals for an Example Carrier on Out‑of‑Sample Data.}
    \label{fig:interval_compare}
\end{figure}

Additionally, Figure \ref{fig:interval_compare} displays the 95\% prediction intervals created by the three contextual distribution learning methods for an example carrier. 
The visualization aligns with the pinball loss results, showing that QRF generally provides the narrowest interval widths. While the intervals generated by calibrated-CLF cover most of the observed values, their coverage is less effective at the two extremes compared to the two QR methods. The inferior performance of calibrated-CLF at extreme quantiles is likely due to the discretization of the output space, which limits precision, and the cumulative probability method, which can cause jumps when selecting quantiles. In contrast, regression tree and QRF provide finer-grained predictions, leading to more accurate and narrower intervals at the extremes.

In summary, MC-CLF models excel at providing distributional forecasts, whereas QR models are better suited to generating prediction intervals on this dataset. Since QRF benefits from ensembling multiple regression trees to reduce variance and improve accuracy and consistently outperform a single regression tree, only QRF is kept for QR in the subsequent CSO evaluations.

\begin{table}[!ht]
\caption{Out-of-Sample Average Pinball Loss and CRPS Across All Carriers.}
\centering
\begin{tabular}{lrrrrrrrrrr}
\toprule
 {} & \multicolumn{5}{c}{Pinball Loss} & CRPS \\
Model &        $q_{0.025}$ &   $q_{0.05}$ &   $q_{0.5}$ &     $q_{0.95}$ &     $q_{0.975}$ &  \\
\midrule
Calibrated-CLF          & 0.199 & 0.198 & 0.676 & 0.526 & 0.304 & \textbf{0.386} \\
QRF             & \textbf{0.033} & \textbf{0.057} & \textbf{0.263} & \textbf{0.143} & 0.096 & 0.402 \\
Regression Tree & 0.064 & 0.109 & 0.302 & 0.147 & \textbf{0.093} & 0.515 \\
\bottomrule
\end{tabular}
\label{tab:quantile_metrics}
\end{table}

\subsection{Baseline Models}

To compare with the proposed contextual optimization models, four baseline order-fulfillment algorithms are introduced below.

\begin{enumerate}
    \item {\em Pure Cost-Driven Greedy Heuristic (Greedy)}: Mimics the industrial partner’s fulfillment optimization engine by selecting, for each unit of an SKU, the carrier–location pair that minimizes the fulfillment cost. Formally, for each order it solves
      \(
      \min_{\mathbf{z} \in \mathcal{Z}} \{\sum_{i, k, \ell } (c_{\ell}^{fixed} + c_{ik \ell}^{ship})z_{ik
        \ell}\}.
      \)

    \revision{
    \item \emph{Empirical‑SAA}: Follows the SAA framework of Section~\ref{sec:c-saa}, except that each deviation sample is drawn from the overall empirical distribution $P^{E}$ of the training set, without using any contextual information.  For each \(n\in[N_1]\) and every \((k,\ell)\in\mathcal K\times\mathcal L\),
  \(
    ~\xi_{k\ell}^n \sim P^E,
    ~\bm\xi^n = \bigl(\xi_{k\ell}^n\bigr)_{k\in\mathcal K,\,\ell\in\mathcal L}.
  \)
    The empirical SAA problem is
    \(
    \min_{\bm z\in\mathcal Z}\;\frac{1}{N_1}\sum_{n=1}^{N_1} g\bigl(\bm z,\bm{\xi}^n\bigr).
    \)

      \item {\em C-Empirical-SAA}: A slightly smarter SAA baseline that maintains separate empirical distributions \(P^E_{k}\) for each carrier.  For each \(n\) and \((k, \ell)\),
      \(
        \xi_{k \ell}^n \sim P^E_{k},  ~\bm\xi^n = \bigl(\xi_{k\ell}^n\bigr)_{k\in\mathcal K,\,\ell\in\mathcal L},
      \)
      and it solves the same SAA problem as above.}

   \item \emph{Point Predict–Then–Optimize (PTO)}: Trains per-carrier regression models \(h_{k}(\cdot)\) to estimate the conditional mean deviations for each \((k,\ell)\in\mathcal K\times\mathcal L\), producing the point‐prediction vector \(\hat{\mathbf{d}}=\{h_{k}(\mathbf{s}_{k\ell})\}_{k \in \mathcal{K}, \ell \in \mathcal{L}}\).  It then solves the deterministic problem
  \(
    \min_{\bm z\in\mathcal Z}g\bigl(\bm z,\hat{\mathbf{d}}\bigr).
  \)
  
   This method is effective only if the objective function
    $g(\mathbf{z}, \tilde{\mathbf{d}})$ is linear in
    $\tilde{\mathbf{d}}$, in which case estimating the conditional
    distribution simplifies to estimating the conditional mean. This
    is because, by linearity of expectation,
    $\mathbb{E}_{\tilde{\mathbf{d}} \sim
      \mathbb{P}(\mathbf{d} \mid
      \mathbf{s})}[g(\mathbf{z}_o, \tilde{\mathbf{d}})] =
    g(\mathbf{z}, \mathbb{E}_{\tilde{\mathbf{d}} \sim
      \mathbb{P}(\mathbf{d} \mid
      \mathbf{s})}[\tilde{\mathbf{d}}])$. 
   However, \(g(\bm z,\tilde{\mathbf{d}})\) is nonlinear due to asymmetric penalty costs of the deviations, hence replacing the distribution by its mean may bias the solution.  For each carrier, the regression model with the lowest cross‐validated mean squared error is selected (see Appendix H for details).
\end{enumerate}

\subsection{Experimental Setting}
This section outlines the computational environment, simulation framework, and parameter configurations used in the experiments.

\subsubsection{Implementation Environment}
The optimization models in this study were implemented in Python 3.9 using Gurobi 11.0.0, with eight threads allocated to each order instance. A termination criterion of a 1\% MIP gap or a 5-minute time limit was applied. All experiments were conducted on a server
equipped with dual-socket Intel Xeon Gold 6226 CPUs, each featuring 24
cores running at 2.7 GHz.

\subsubsection{Simulation Framework}
One of the difficulties in this application is that the ground truth
is not available for options that were not selected. To overcome this
limitation, a simulation environment was built to generate synthetic
``ground truth'' when an algorithm selects an option. Specifically, delivery deviations
were drawn from a multinomial distribution, with probabilities
estimated by fitting a fine-tuned, calibrated multi-class
classification model based on orders from the last week of the dataset (testing period). Obviously, this data was not accessible to the algorithm and only used by the simulator to replicate real-world conditions. \revision{The training and calibration procedures of the simulator closely followed the methodology outlined in Section \ref{sec:mc-clf}. An evaluation of the simulator is provided in Appendix J.

For the following experiments, as mentioned in Section \ref{Evaluation of Contextual Distributional Learning of Delivery Deviations}, all ML models were trained on order data from the first three weeks. Evaluation was conducted on 5,000 multi-item orders, randomly sampled from the final week. To account for variability, solution metrics were computed across 50 random realizations. \textit{All reported metrics were normalized by setting the Greedy baseline to 1 for each random seed; true performance values were omitted to due to data confidentiality.}}

\subsubsection{Model Names and Parameters}

In the following, the C-SAA problems using samples from the
contextual distribution estimates provided by the MC-CLF and QRF
models are referred to as C-SAA-CLF and C-SAA-QRF,
respectively. Similarly, the C-RO problems are referred
to as C-RO-\{B,D\}-CLF and C-RO-\{B,D\}-QRF, respectively.

By default, a constant $\beta_k = 0.5$ consolidation discount factor for all $k \in \mathcal{K}$ was considered. The default penalty for early delivery $\gamma^-$ was $0.1$ and that for late delivery $\gamma^+$ was $10$. For SAA-based methods, the default was to have $Q = 20$
independent replications each with a sample size of $N_1 = 30$. The
obtained solutions were then evaluated with a default sample size of
$N_2 = 300$. The default sample size used for C-RO-D was $N = 200$. As
for C-RO-B, a default of $95\%$ interval was considered. Moreover, a
scaling factor $\eta \in [0, 1]$ (default of 0.5) determined the
uncertainty budget in C-RO-B, defined as \(
B = \eta\,|\mathcal{K}|\,|\mathcal{L}|.
\)

\revision{
\subsection{Overall Results}

Table \ref{tab:methods_comparison} reports the normalized average and worst‐case objective values across 50 random seeds and all order instances. Based on these results, the following observations are drawn.

\paragraph{Average Performance}
The C-SAA variants and C-RO-B-CLF method lead to more  reduction in the expected objective than other methods. C‑SAA‑CLF achieves the greatest reduction, with a normalized mean of 0.820 $\pm$ 0.003 (–18.0\% versus Greedy), followed by C‑SAA‑QRF at 0.844 $\pm$ 0.003 (–15.6\%) and C‑RO‑B‑CLF at 0.842 $\pm$ 0.003 (–15.8\%). By comparison, Empirical‑SAA (0.962 $\pm$ 0.002) and PTO (0.878 $\pm$ 0.003) offer smaller improvements, and C‑Empirical‑SAA (0.870 $\pm$ 0.003) performs only marginally better than PTO.

\paragraph{Worst-Case Performance}
The C-RO-D methods deliver the strongest guarantees under the worst simulated outcome: C‑RO‑D‑CLF attains 0.779 (–22.1\%), C‑RO‑D‑QRF 0.787 (–21.3\%), and C‑RO‑B‑CLF 0.811 (–18.9\%). These three methods uniformly dominate the baseline and other approaches in the worst‐case metric.

Overall, the C-SAA approaches excel in minimizing the expected objective, while the C‑RO‑D variants excel in worst‐case robustness. Notably,  C‑RO‑B‑CLF strikes a strong balance between average efficiency and tail‐risk protection.

\begin{table}[!ht]
\centering
\begin{threeparttable}
\caption{Average and Worst-Case Realized Objective Values Across Different Order Fulfillment Methods.}
\label{tab:methods_comparison}
\begin{tabular}{lcc} 
\toprule
Method & Average Objective (95\% CI) & Worst-Case Objective\\ \midrule
Greedy            & 1.000 $\pm$ 0.000 & 1.000 \\
PTO               & 0.878 $\pm$ 0.003 & 0.919 \\
Empirical-SAA     & 0.962 $\pm$ 0.002 & 1.049 \\
C-Empirical-SAA   & 0.870 $\pm$ 0.003 & 0.896 \\
C-SAA-CLF      & \textbf{0.820 $\pm$ 0.003} & 0.823 \\
C-SAA-QRF         & \textbf{0.844 $\pm$ 0.003} & 0.842 \\
C-RO-B-CLF     & \textbf{0.842 $\pm$ 0.003} & \textbf{0.811} \\
C-RO-B-QRF        & 0.861 $\pm$ 0.003 & 0.824 \\
C-RO-D-CLF     & 0.928 $\pm$ 0.003 & \textbf{0.779}\\
C-RO-D-QRF        & 0.914 $\pm$ 0.003 & \textbf{0.787} \\
\bottomrule
\end{tabular}
\footnotesize
\begin{tablenotes}
\item Note: Top 3 methods for each metric are highlighted in bold. 
\end{tablenotes}
\end{threeparttable}
\end{table}
}

\subsection{Effectiveness of C-SAA}

This section further compares the performance of the C-SAA methods with the two Empirical-SAA baselines and the PTO.

\subsubsection{Value of Contextual Information}

\begin{figure}[!ht]
\centering
    \begin{subfigure}[b]{0.4\textwidth}
        \centering
        \includegraphics[width=\textwidth]{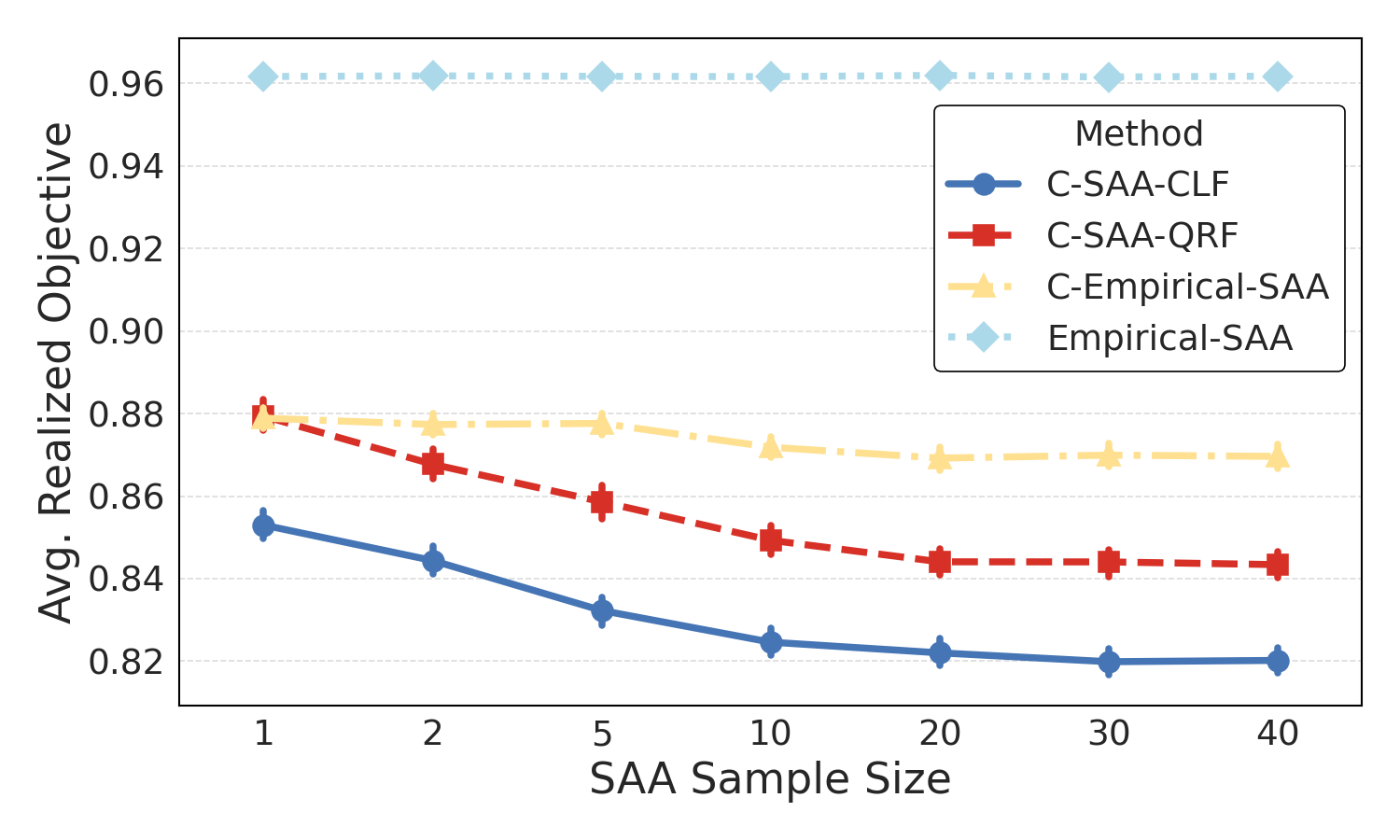}
        \caption{Performance in Average Realized Objective Value.}     
        \label{fig:saa_sample_size_objVal}
    \end{subfigure}
    \begin{subfigure}[b]{0.4\textwidth}
        \centering
        \includegraphics[width=\textwidth]{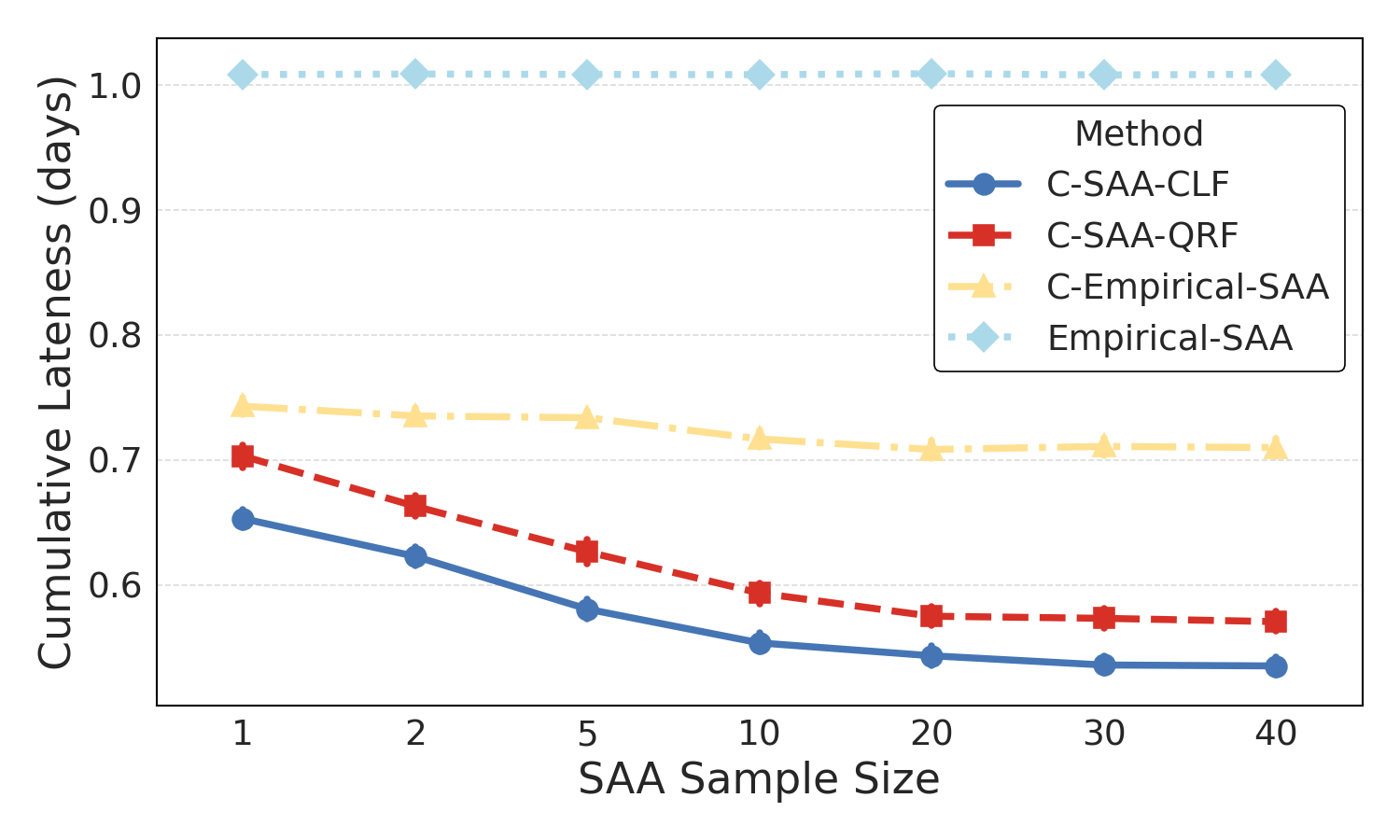}
        \caption{Performance in Cumulative Lateness.}
        \label{fig:saa_sample_size_cum_lateness}
    \end{subfigure}
    \caption{C-SAA vs. Empirical SAA Across Different Sample Sizes.}
    \label{fig:saa_sample_size}
\end{figure}

To facilitate a closer examination of the performance of the
SAA-based models, Figure \ref{fig:saa_sample_size} shows the impact of
varying sample sizes ($N_1$).  On the one hand, larger sample sizes should provide
a more accurate approximation of the expected value in Problem
(\ref{eq: CSOFP}), assuming the delivery deviation distributions were
well-estimated. On the other hand, it is important to maintain a
reasonably small sample size to ensure computational efficiency. For
both C-SAA methods, the effect of increasing sample sizes is
obvious. The average objectives are reduced from 0.85 to 0.82 for
C-SAA-CLF and 0.88 to 0.84 for C-SAA-QRF, indicating their
increasing alignment with the true risk-neutral objective
function. This improvement is particularly evident in their enhanced
effectiveness at avoiding cumulative lateness, defined to be the total
number of late delivery days across all units in the orders. The
cumulative lateness decreases by 14\% for C-SAA-CLF and by 12\% for
C-SAA-QRF.

In contrast, both Empirical-SAA and C-Empirical-SAA struggles to improve there performance in
both the average realized objective value and the cumulative lateness as the sample size grows. This is because relying solely on the empirical distribution fails to capture the variability in delivery deviation distributions, which may depend on contextual information specific to each order. A biased distributional forecasts
thus guide the model wrongly to create solutions with inferior
quality.  Based on these observations, both contextual distributional
forecasting methods proposed in the paper are successful in this
application, demonstrating a stronger performance compared to using
empirical distribution for delivery deviation forecasting.  In summary,
the value of contextual information is particularly important in order
fulfillment decision-making, where the uncertainty of delivery deviation
can vary drastically across different states of the environment.

\subsubsection{Value of Distributional Forecasts}

\begin{figure}[!ht]
\centering
    \begin{subfigure}[b]{0.45\textwidth}
        \centering
        \includegraphics[width=\textwidth]{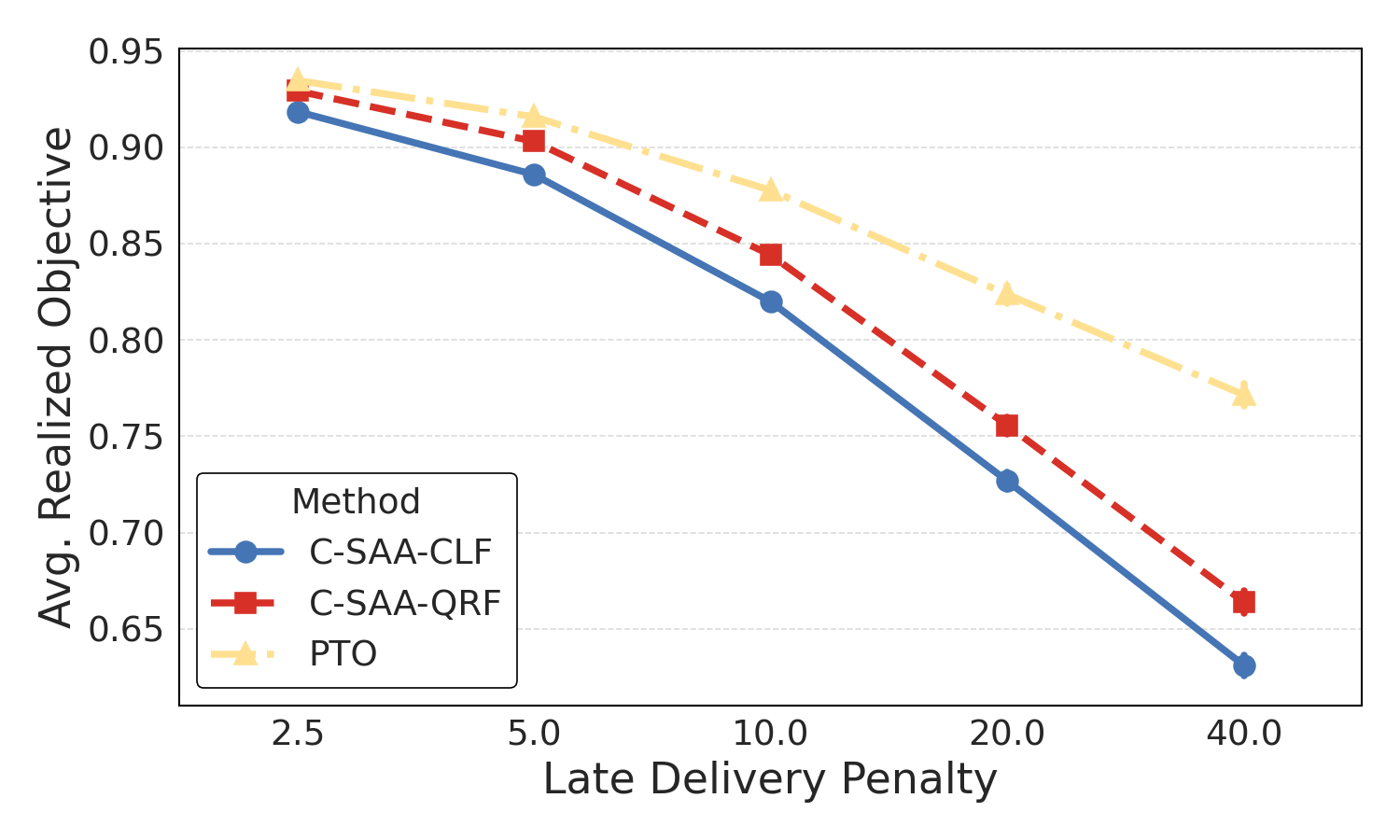}
        \caption{Performance in Average Realized Objective Value.}     
        \label{fig:csaa_vs_pointpred_objVal}
    \end{subfigure}
    \hfill
    \begin{subfigure}[b]{0.45\textwidth}
        \centering
        \includegraphics[width=\textwidth]{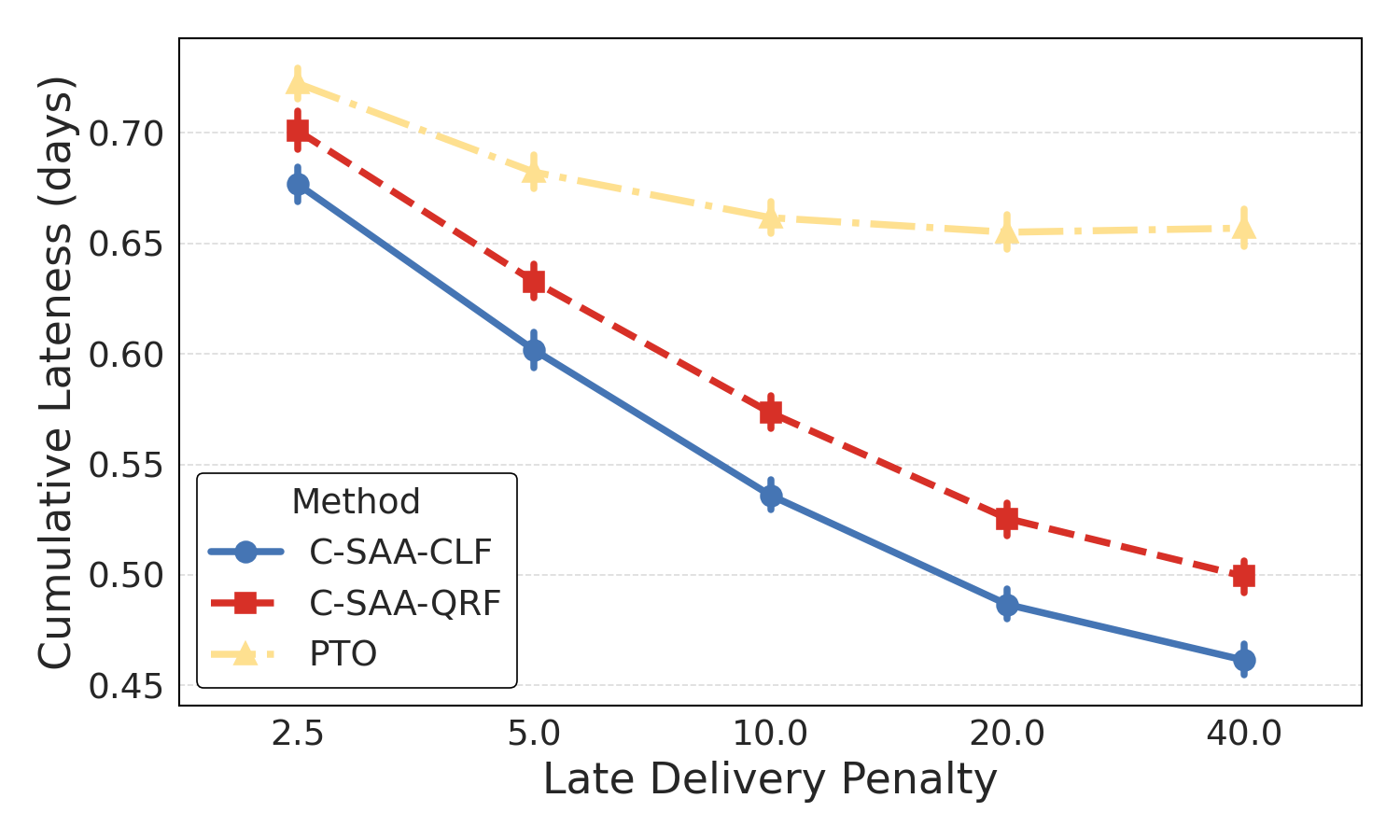}
        \caption{Performance in Cumulative Lateness.}
        \label{fig:csaa_vs_pointpred_cum_lateness}
    \end{subfigure}
    \caption{C-SAA vs. PTO Across Different Late Delivery Penalties.}
    \label{fig:csaa_vs_pointpred}
\end{figure}

This section examines the advantage of employing distributional
forecasting models over the simpler and more common point prediction
model within the CSO framework. As supported by theory, simply using the conditional
mean estimate of deviations to replace the expected value in CSOFP can
introduce bias when approximating the true objective value, thus degrading solution quality. Figure \ref{fig:csaa_vs_pointpred} confirms that, under changing $\gamma^+$, CSO methods with distributional estimator (C‑SAA‑CLF and C‑SAA‑QRF) consistently outperform PTO in both average realized objective value (Figure \ref{fig:csaa_vs_pointpred_objVal}) and cumulative lateness (Figure \ref{fig:csaa_vs_pointpred_cum_lateness}).
In fact, the performance differences become
more pronounced as late delivery penalty ($\gamma^+$) increases.

When $\gamma^+$ raises from 2.5 to 40, PTO achieves only a 7\% reduction in cumulative lateness, whereas C-SAA-CLF and
C-SAA-QRF realize reductions of approximately 22\%
and 20\%, respectively. These results demonstrate that distributional forecasting not only enables lower overall lateness but also leverages higher penalties more effectively, yielding a sharper decline in average objective and making C‑SAA particularly well suited to settings where delivery timeliness is critical.

\revision{
\subsection{Cost of Robustness}
To quantify the trade‑off between robustness and efficiency, C‑RO‑B‑CLF and C‑RO‑D‑CLF were compared with C‑SAA‑CLF as the late‑delivery penalty ($\gamma^+$) increases (Figure \ref{fig:cro_vs_csaa}). 

As depicted in Figure \ref{fig:cro_vs_csaa_objVal}, at $\gamma^+ =$ 2.5 , average objectives are approximately 0.92 for C‑SAA‑CLF, 0.95 for C‑RO‑B‑CLF, and 1.00 for C‑RO‑D‑CLF. As $\gamma^+$ rises to 40, these values decline in parallel—C‑SAA‑CLF to 0.63, C‑RO‑B‑CLF to 0.67, and C‑RO‑D‑CLF to 0.73. Throughout this range, C‑SAA‑CLF consistently achieves the lowest average objective, followed by C‑RO‑B‑CLF and then C‑RO‑D‑CLF.

Worst‑case objectives (Figure \ref{fig:cro_vs_csaa_worst_objVal}) underscore the robustness of the C‑RO methods: C‑RO‑D‑CLF consistently yields the strongest tail protection, dropping from 0.88 at $\gamma^+ =$ 2.5 to 0.72 at $\gamma^+ =$ 40. C‑RO‑B‑CLF falls from 0.885 to 0.775, while C‑SAA‑CLF declines more steeply from 0.925 to 0.73. Notably, at $\gamma^+ =$ 40, C‑SAA‑CLF nearly matches C‑RO‑D‑CLF, demonstrating that C-SAA can approach robust performance under extreme penalty scenarios.

This robustness manifests in cumulative lateness (Figure \ref{fig:cro_vs_csaa_cum_lateness}). Both C‑RO variants tightly control lateness with minimal sensitivity to $\gamma^+ $, whereas C‑SAA‑CLF’s lateness drops substantially from 0.67 at $\gamma^+ =$ 2.5 to 0.46 at $\gamma^+ =$ 40. These robustness gains incur higher expected fulfillment costs (Figure \ref{fig:cro_vs_csaa_Fulfillment_Costs}). C‑RO‑D‑CLF’s average cost increases from 1.07 to 1.20, and C‑RO‑B‑CLF’s from 1.01 to 1.05, and C‑SAA‑CLF from 0.96 to 1.05. Thus, while all methods incur higher costs under stricter penalties, C‑RO‑B‑CLF exhibits the smallest relative cost increase, C‑SAA‑CLF and C‑RO‑D‑CLF have more aggressive increases in costs due to more conservative decisions.

A final strategic consideration lies in selecting $\gamma^+$ to reflect firm priorities. Lower $\gamma^+$ emphasizes immediate cost savings—maximizing short‑term profits—whereas higher $\gamma^+$ prioritizes on‑time delivery, fostering customer satisfaction, loyalty, and long‑term value. The CSOFP framework encodes this trade‑off directly via $\gamma^+$; robust formulations (C‑RO‑B‑CLF and C‑RO‑D‑CLF) further reduce the need for precise $\gamma^+$ tuning by delivering stable, worst‑case performance across penalty levels, while C‑SAA‑CLF can leverage high $\gamma^+$ most aggressively when strict timeliness is paramount.

\begin{figure}[!ht]
\centering
    \begin{subfigure}[b]{0.45\textwidth}
        \centering
        \includegraphics[width=\textwidth]{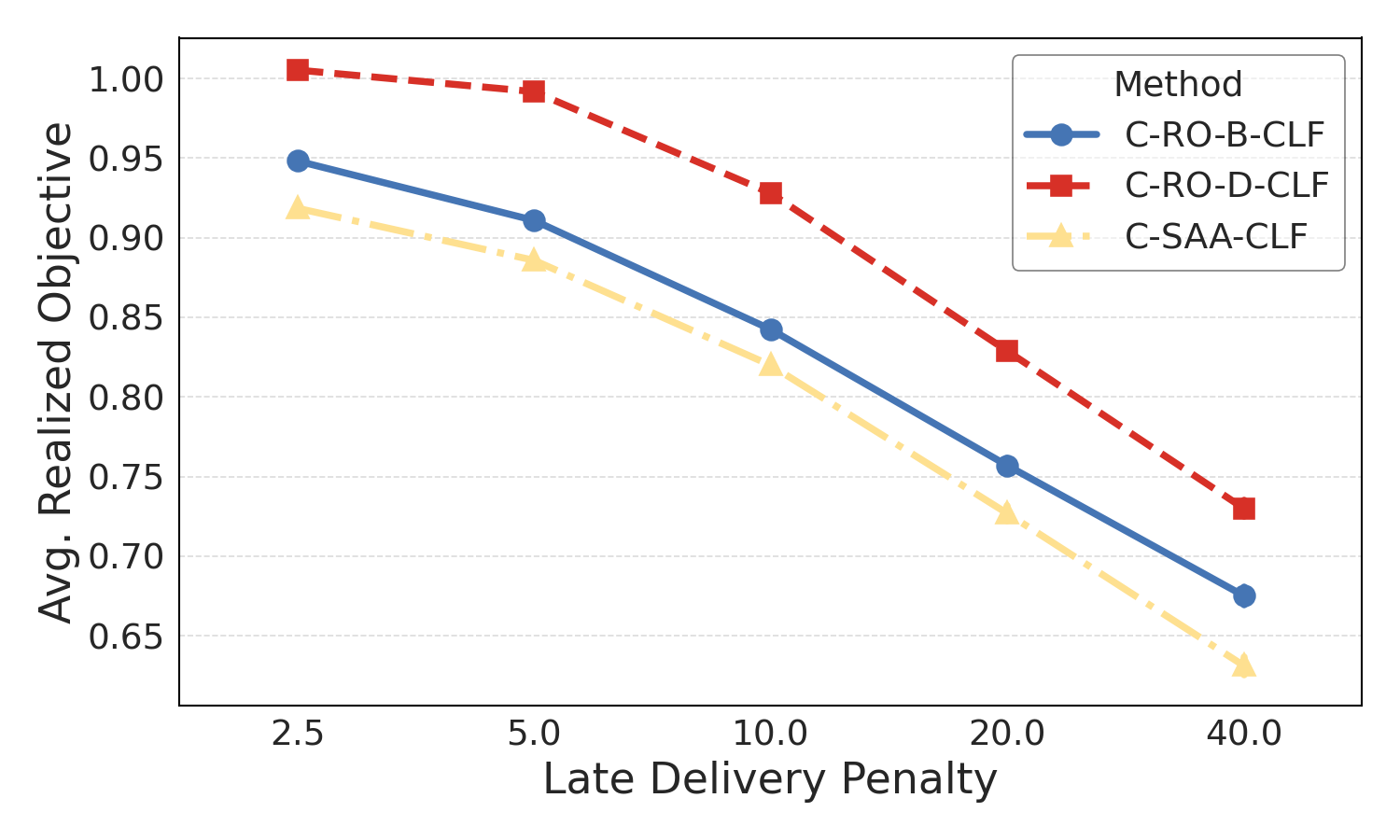}
        \caption{Performance in Average Realized Objective Value.}
    \label{fig:cro_vs_csaa_objVal}
    \end{subfigure}  
    \hfill
    \begin{subfigure}[b]{0.45\textwidth}
        \centering
        \includegraphics[width=\textwidth]{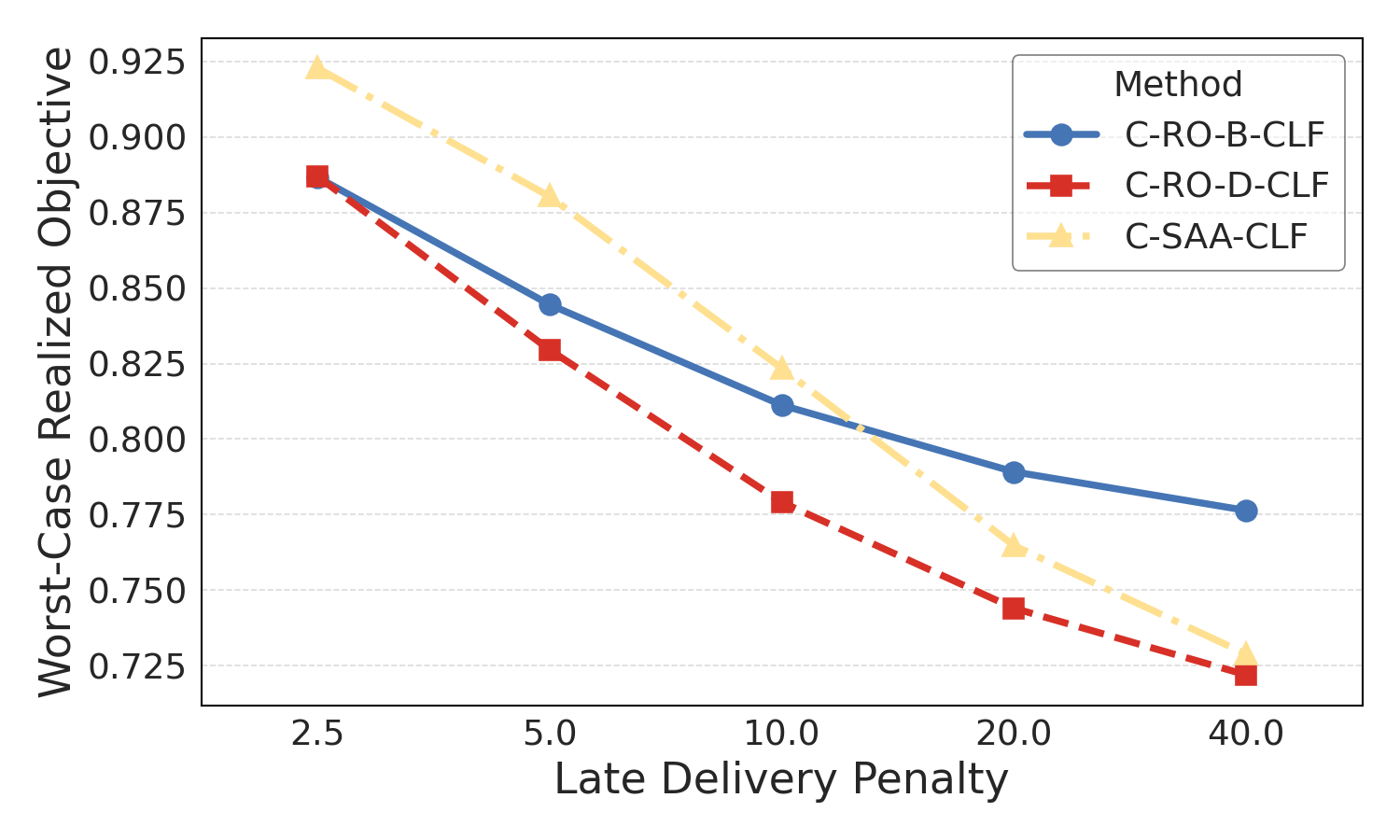}
        \caption{Performance in Worst-Case Realized Objective Value.}
    \label{fig:cro_vs_csaa_worst_objVal}
    \end{subfigure} 
    \begin{subfigure}[b]{0.45\textwidth}
        \centering
        \includegraphics[width=\textwidth]{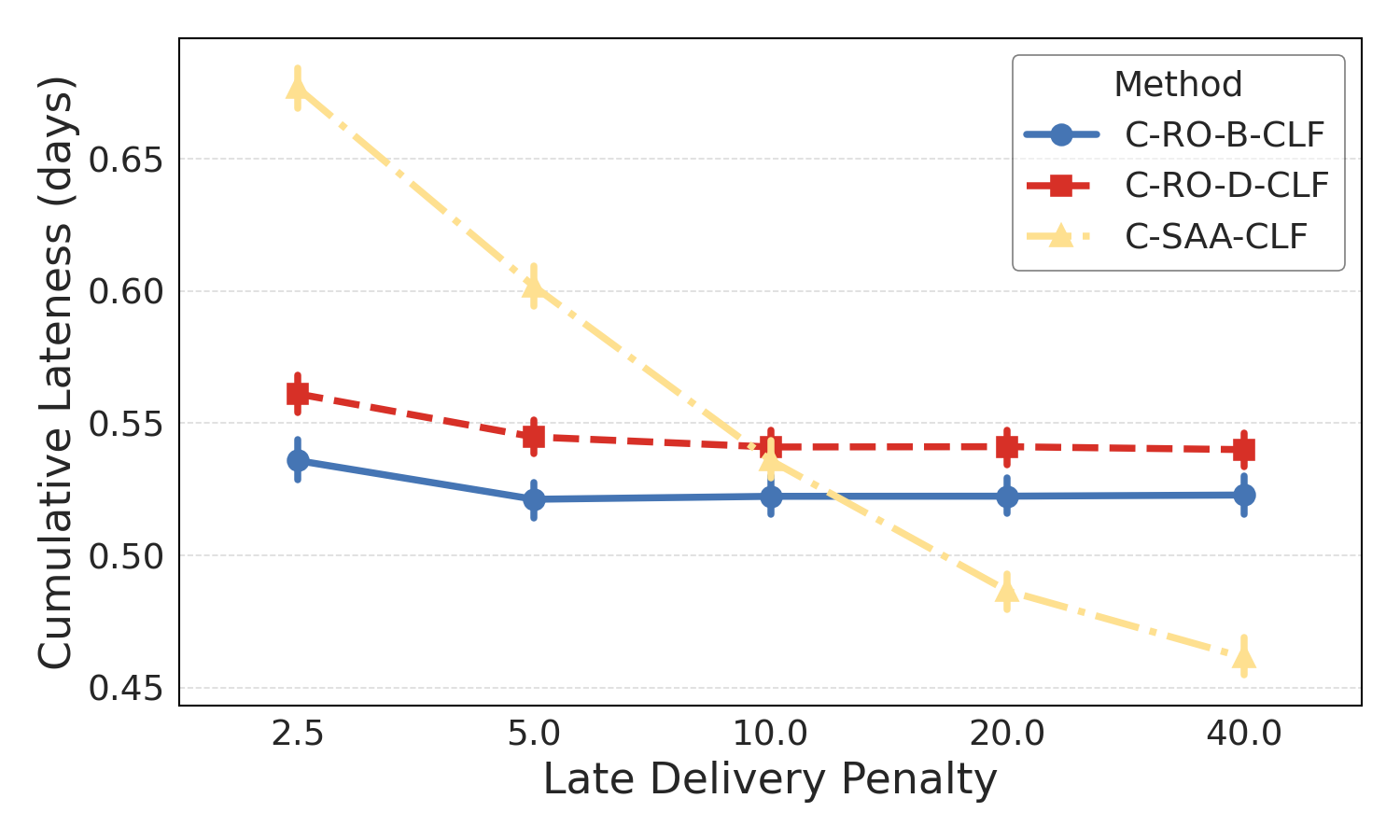}
        \caption{Performance in Cumulative Lateness.}
        \label{fig:cro_vs_csaa_cum_lateness}
    \end{subfigure}
    \hfill
    \begin{subfigure}[b]{0.45\textwidth}
        \centering
        \includegraphics[width=\textwidth]{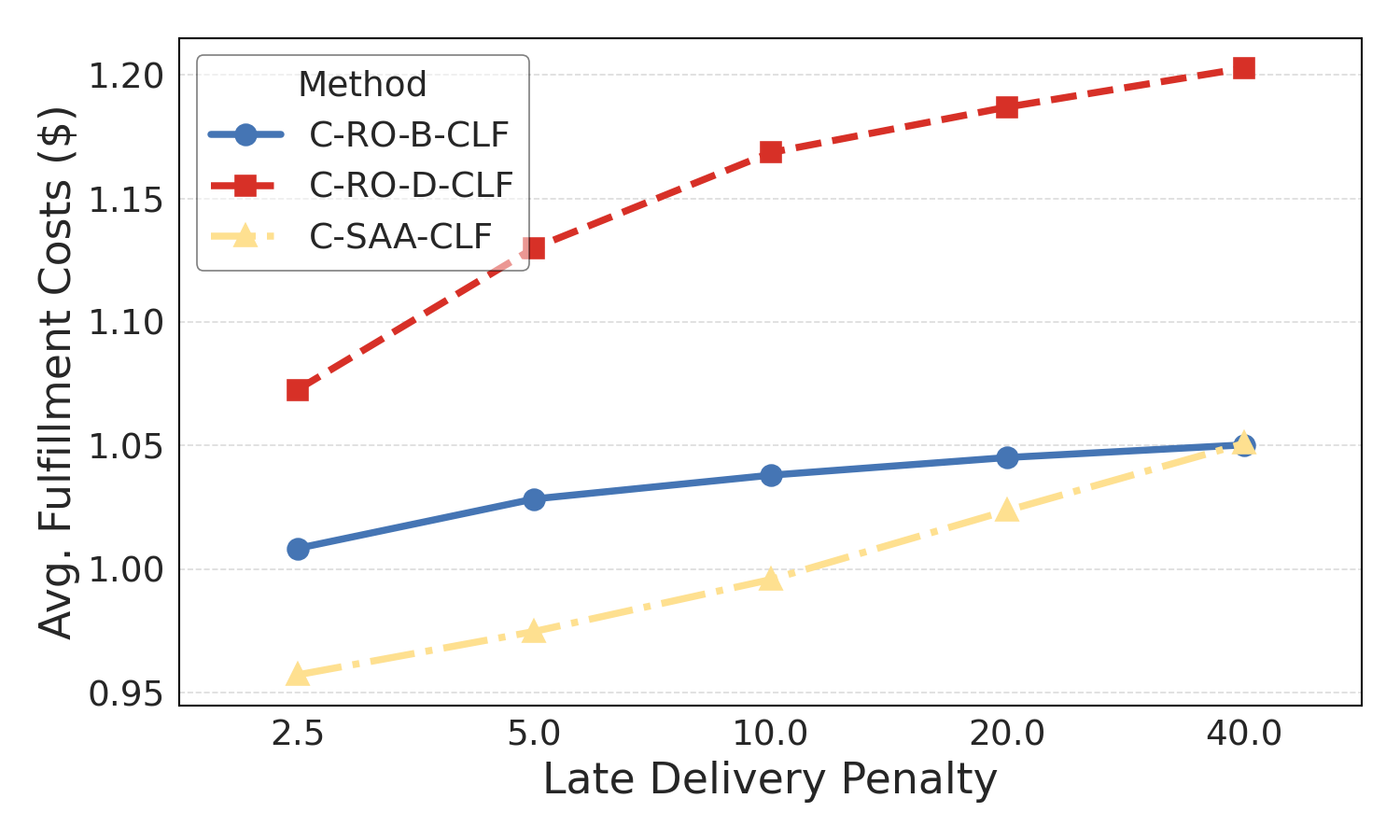}
        \caption{Performance in Fulfillment costs.}
        \label{fig:cro_vs_csaa_Fulfillment_Costs}
    \end{subfigure}  
        \caption{C-RO vs. C-SAA Across Different Late Delivery Penalties.}
\label{fig:cro_vs_csaa}
\end{figure}
}

\subsubsection{Flexibility of C-RO-B}
The C‑RO models effectively hedge against worst‐case delivery deviations and potential misspecifications in the prediction models. However, this robustness incurs more conservative decisions and higher fulfillment costs. The C‑RO‑B variants retain the distribution‐free guarantees while introducing two tuning parameters---prediction‐interval coverage and robust uncertainty budget---that directly control conservativeness. Lower coverage levels and smaller budgets yield less conservative, more cost‐efficient solutions (at the expense of risk protection), whereas higher values emphasize hedging against extreme deviations by sacrificing average performance.

\begin{figure}[!ht]
\centering
    \begin{subfigure}[b]{0.32\textwidth}
        \centering
        \includegraphics[width=\textwidth]{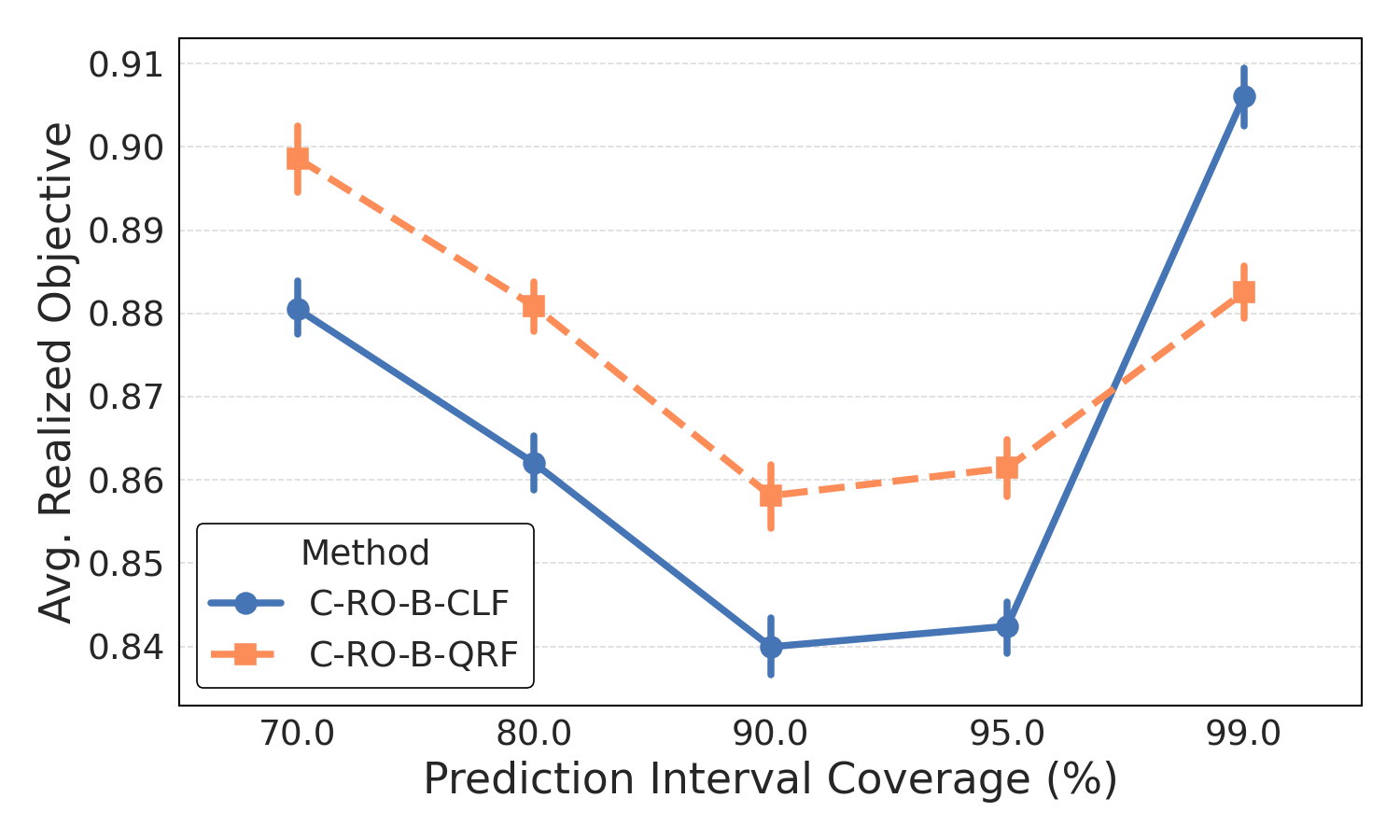}
        \caption{Performance in Avg. Obj.}
    \label{fig:sensitivity_coverage_objVal}
    \end{subfigure} 
    \begin{subfigure}[b]{0.32\textwidth}
        \centering
        \includegraphics[width=\textwidth]{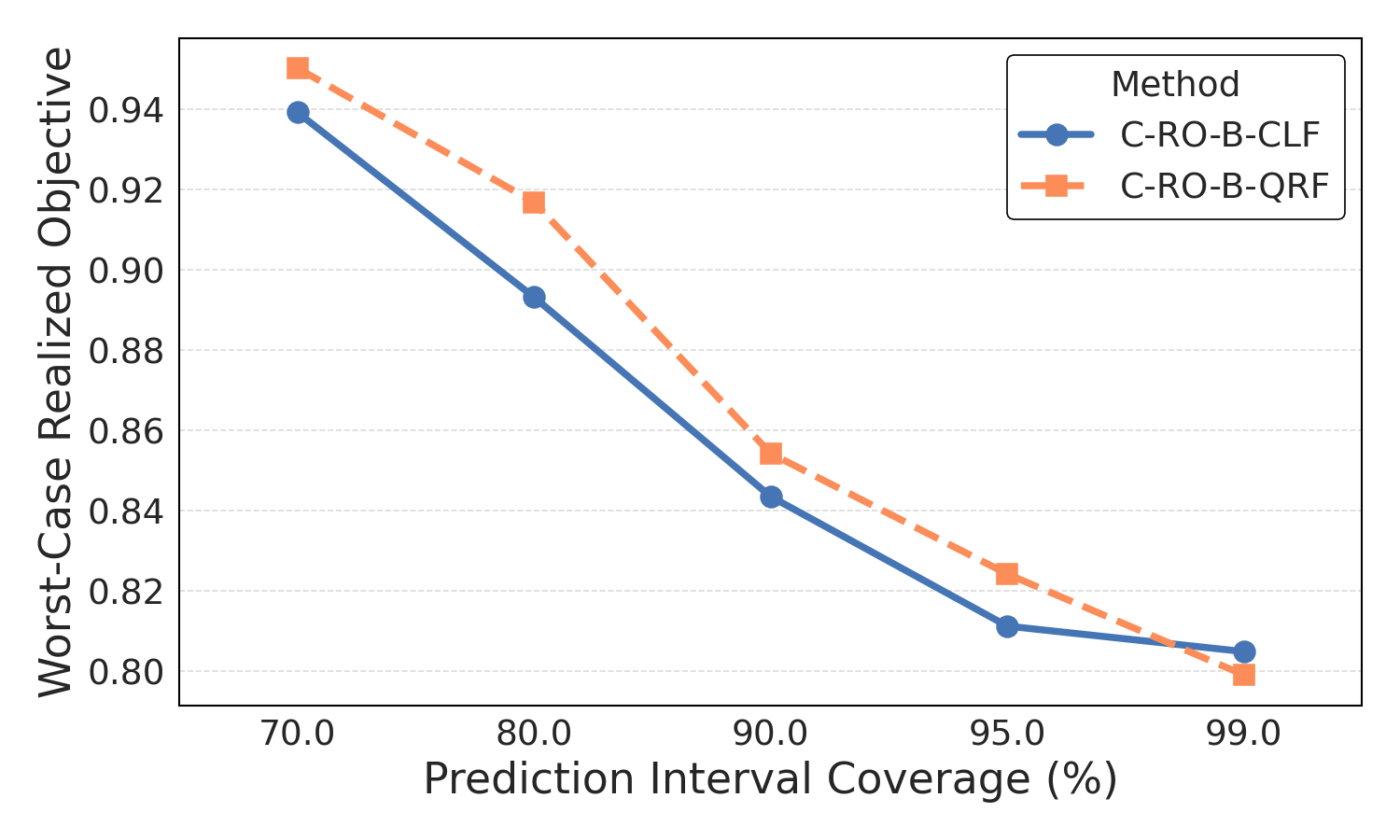}
        \caption{Performance in Worst-Case Obj.}
    \label{fig:sensitivity_coverage_worst_objVal}
    \end{subfigure} 
    \begin{subfigure}[b]{0.32\textwidth}
        \centering
        \includegraphics[width=\textwidth]{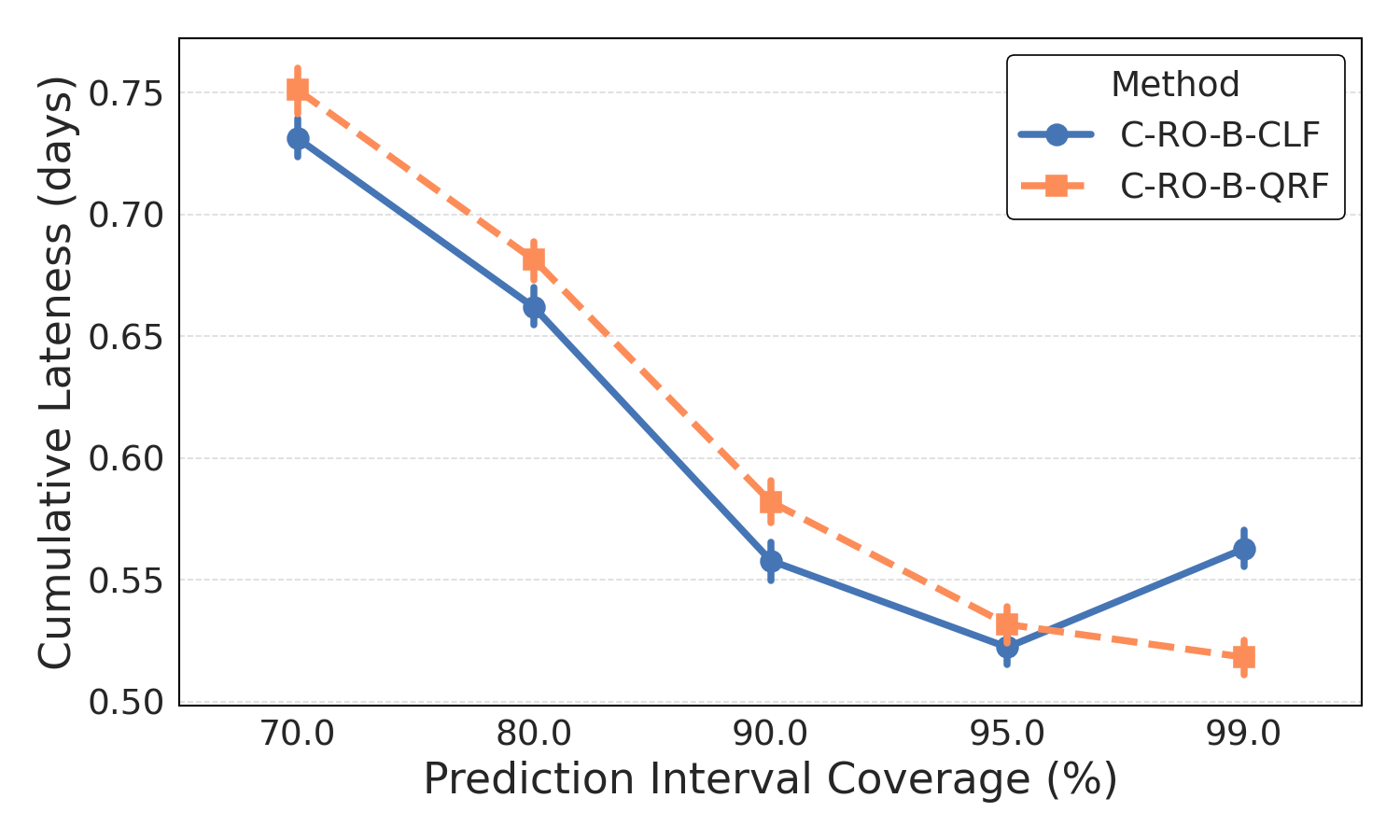}
        \caption{Performance in Cumulative Lateness.}
        \label{fig:sensitivity_coverage_cum_lateness}
    \end{subfigure}    
        \caption{Performance of C-RO-B Across Different Prediction Interval Coverages.}
\label{fig:sensitivity_coverage}
\end{figure}

\begin{figure}[!ht]
\centering
    \begin{subfigure}[b]{0.32\textwidth}
        \centering
        \includegraphics[width=\textwidth]{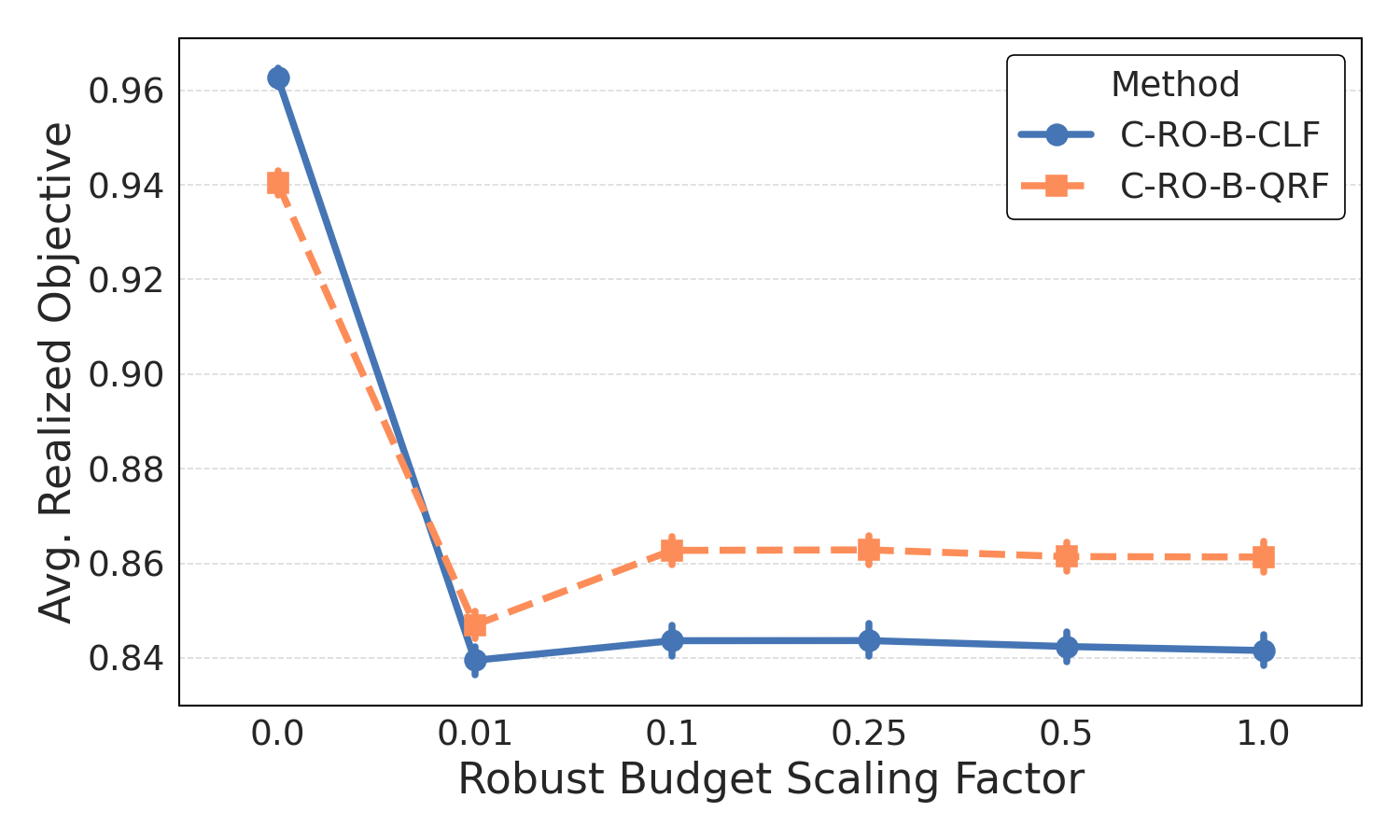}
        \caption{Performance in Avg. Obj.}
    \label{fig:sensitivity_robust_budget_objVal}
    \end{subfigure}
    \begin{subfigure}[b]{0.32\textwidth}
        \centering
        \includegraphics[width=\textwidth]{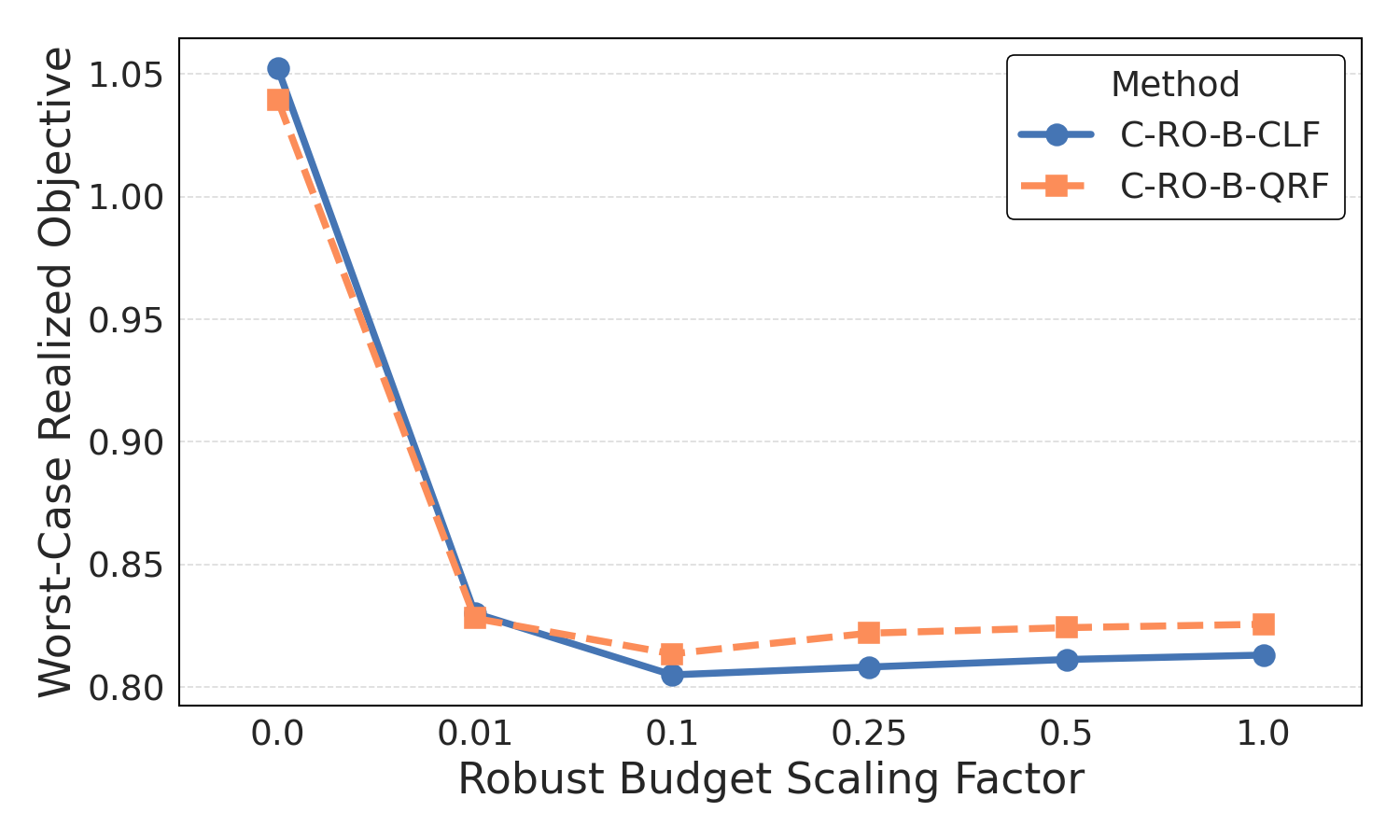}
        \caption{Performance in Worst-Case Obj.}
    \label{fig:sensitivity_robust_budget_worst_objVal}
    \end{subfigure}   
    \begin{subfigure}[b]{0.32\textwidth}
        \centering
        \includegraphics[width=\textwidth]{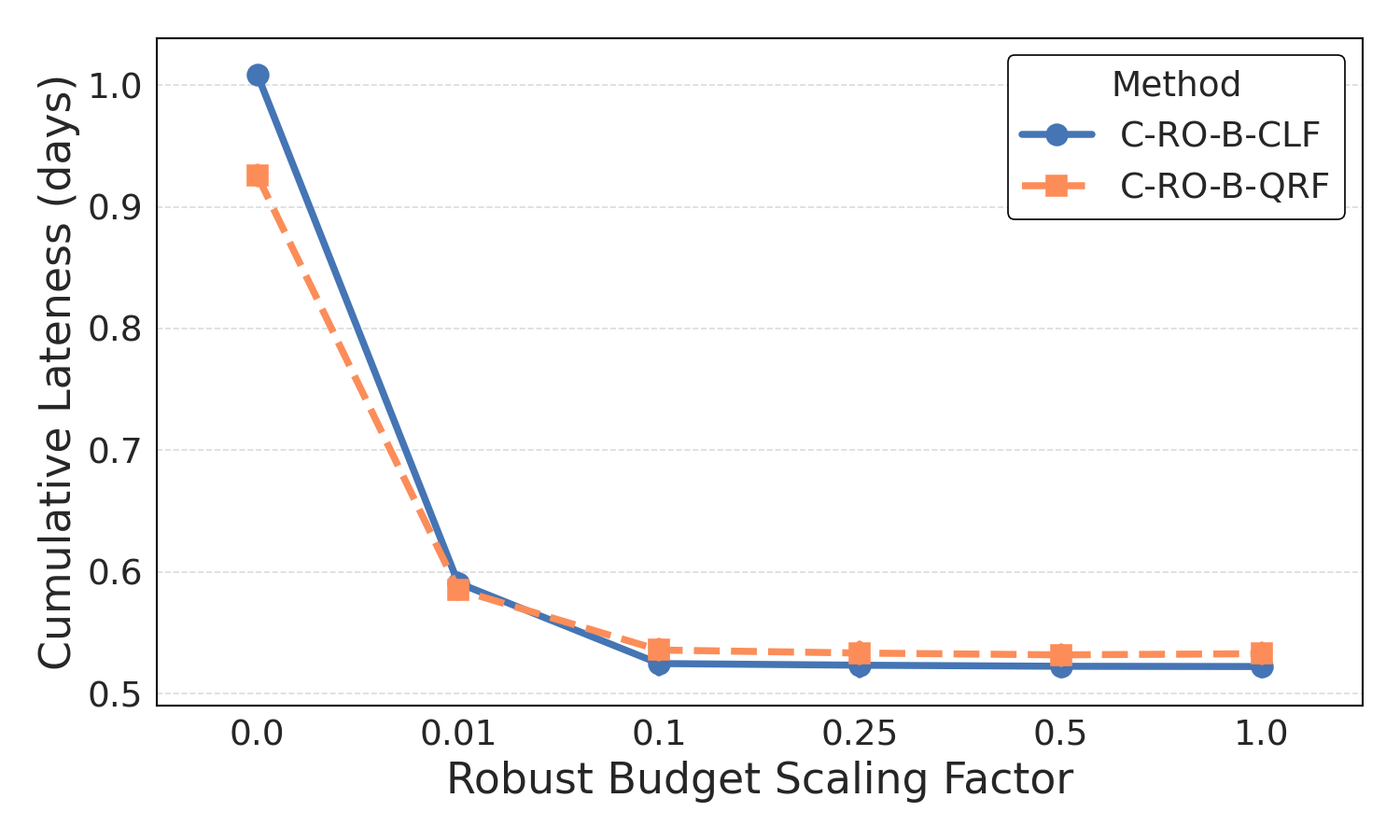}
        \caption{Performance in Cumulative Lateness.}
        \label{fig:sensitivity_robust_budget_cum_lateness}
    \end{subfigure}    
        \caption{Performance of C-RO-B Across Different Robust Uncertainty Budgets.}
\label{fig:sensitivity_robust_budget}
\end{figure}

Figures \ref{fig:sensitivity_coverage} and \ref{fig:sensitivity_robust_budget} summarize how these parameters affect solution quality under both C‑RO‑B‑CLF and C‑RO‑B‑QRF. As the prediction
interval coverage increases from 70\% to 90\%.  In Figure \ref{fig:sensitivity_coverage_objVal}, increasing coverage from 70\% to 90\% reduces the average objective to about 0.84 (C‑RO‑B‑CLF) and 0.86 (C‑RO‑B‑QRF), an approximately
4\% improvement relative to 70\% coverage. This is consistent with the concurrent decline in cumulative lateness shown in Figure \ref{fig:sensitivity_coverage_cum_lateness}. However, pushing coverage beyond 95\% reverses this trend: at 99\% coverage the average objective rises sharply due to the cost of extreme conservatism and, for C‑RO‑B‑CLF, degraded prediction accuracy (reflected in a spike in lateness). By contrast, C‑RO‑B‑QRF continues to reduce lateness at 99\%, illustrating its superior tail modeling. In terms of worst‐case performance (Figure \ref{fig:sensitivity_coverage_worst_objVal}), both methods improve with coverage, but C‑RO‑B‑CLF plateaus after 95\%, whereas C‑RO‑B‑QRF continues to improve at 99\%.

The effect of the robust budget, controlled by the scaling factor \(\eta\), appears in Figure~\ref{fig:sensitivity_robust_budget}. Increasing \(\eta\) from zero to 0.01 leads to a steep drop in the average objective for both methods (Figure~\ref{fig:sensitivity_robust_budget_objVal}). Both method reach their minimum values at \(\eta = 0.01\). At the same time, the worst-case objective follows a similar pattern (Figure~\ref{fig:sensitivity_robust_budget_worst_objVal}): it falls sharply up to \(\eta = 0.1\), then edges back up for both methods as the budget grows. This trend is also seen in cumulative lateness, which drops dramatically until \(\eta = 0.1\) and then levels off (Figure~\ref{fig:sensitivity_robust_budget_cum_lateness}). Together, these trends show that moderate budgets capture most of the gains in both average and worst case performance, while further increases deliver diminishing returns.

\subsection{Impact of Consolidation Factor}
This section provides a discussion on the impact of item
consolidation. CSOFP relies on the coefficients $\beta_k$ to control
the incentive to batch items. The industrial partner does not know the
exact percentage of discount to be offered for item consolidation,
which may depend on a variety of convoluted factors. Therefore, this
section performs sensitivity analysis on a range of possible discount
factors (from 0.1 to 0.7).

\begin{figure}[!ht]
\centering
    \begin{subfigure}[b]{0.45\textwidth}
        \centering
        \includegraphics[width=\textwidth]{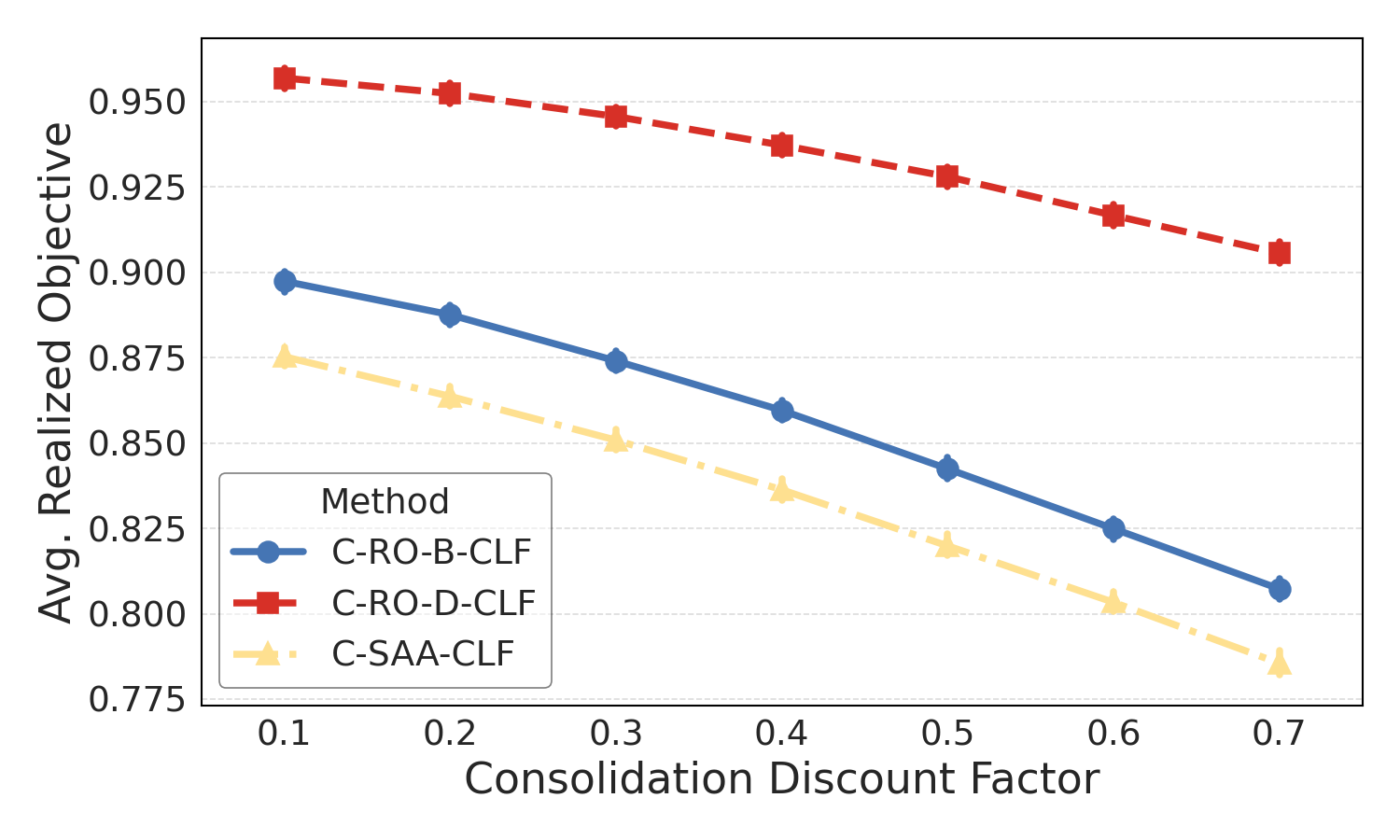} 
        \caption{Performance in Average Realized Objective Value.}
        \label{fig:sensitivity_discount_objVal}
    \end{subfigure}
    \hfill
    \begin{subfigure}[b]{0.45\textwidth}
        \centering
        \includegraphics[width=\textwidth]{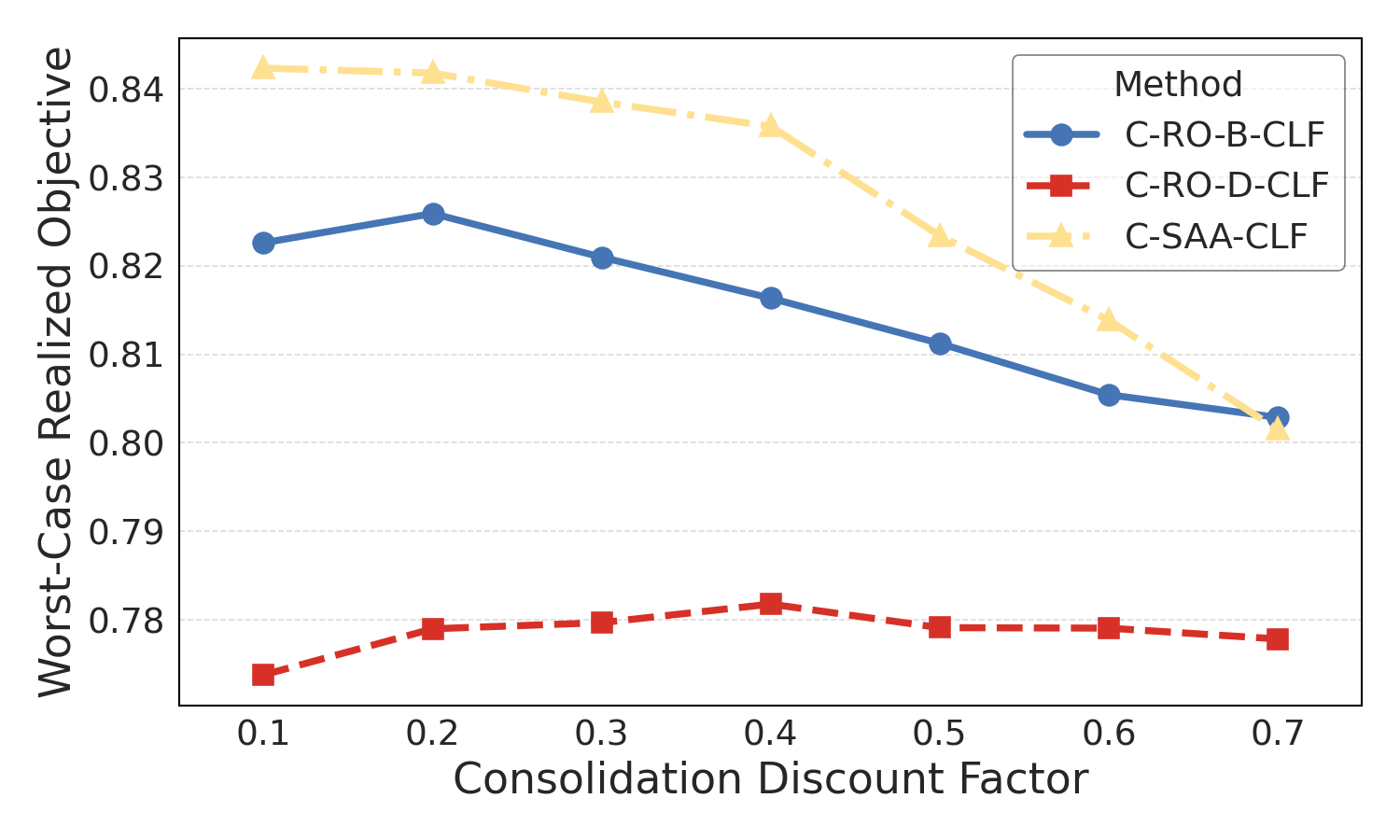} 
        \caption{Performance in Worst-Case Realized Objective Value.}
        \label{fig:sensitivity_discount_worstObj}
    \end{subfigure}
    \begin{subfigure}[b]{0.45\textwidth}
        \centering
        \includegraphics[width=\textwidth]{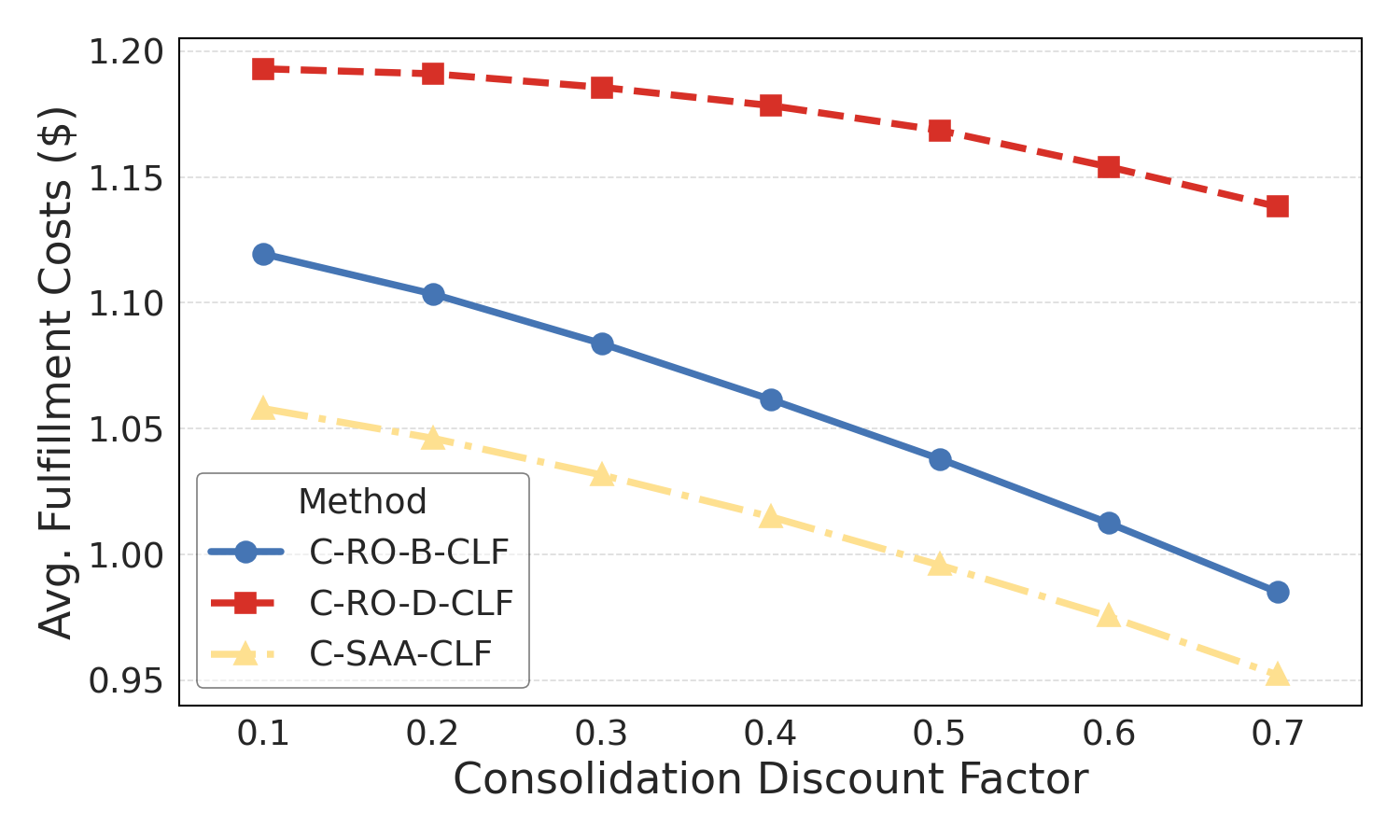} 
        \caption{Performance in Fulfillment Costs.}
        \label{fig:sensitivity_discount_fulfillment_cost}
    \end{subfigure}
    \hfill
    \begin{subfigure}[b]{0.45\textwidth}
        \centering
        \includegraphics[width=\textwidth]{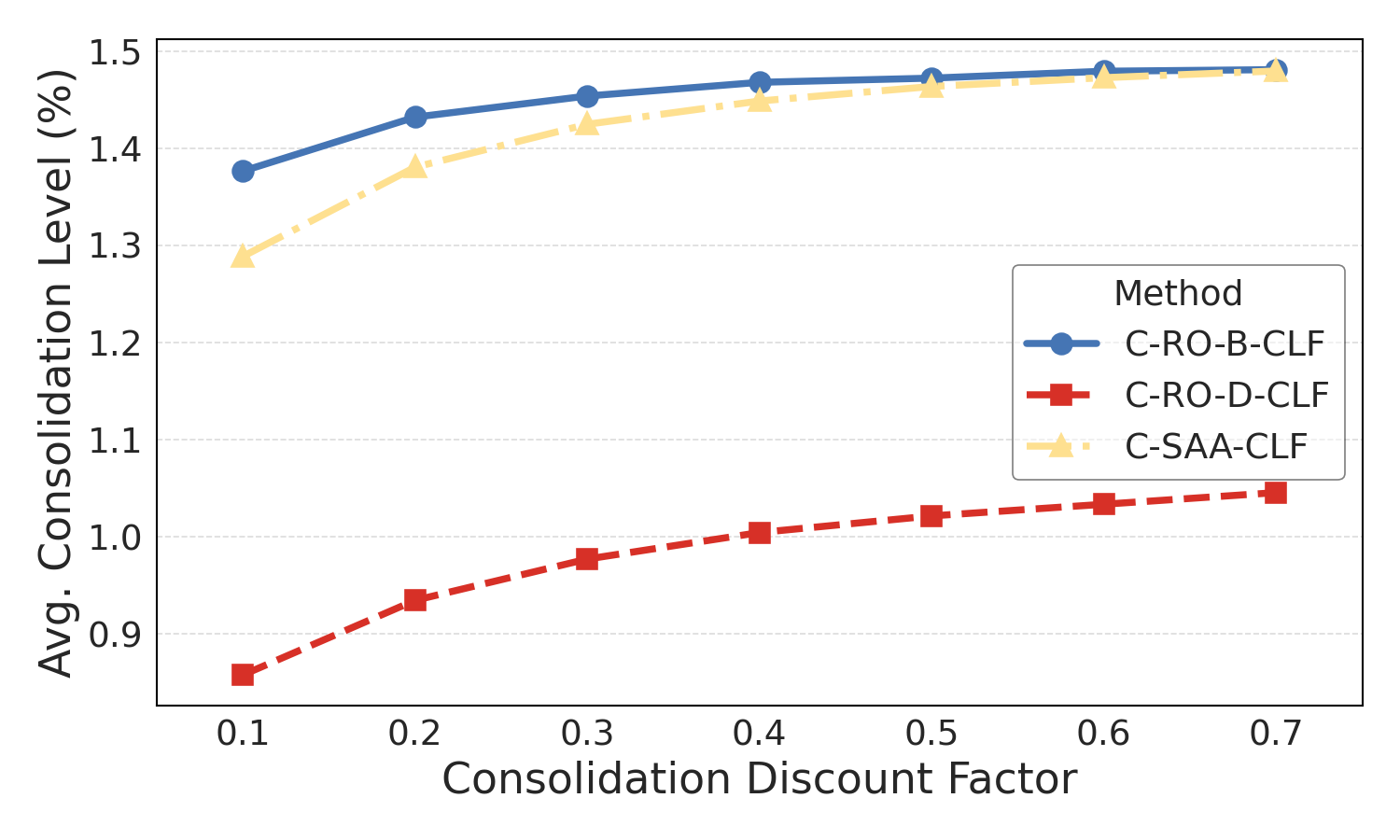} 
        \caption{Performance in Consolidation Level.}
        \label{fig:sensitivity_discount_consolidate_pct}
    \end{subfigure}
     \caption{Performance of CSO Methods Across Different Consolidation Discount Factors.}
     \label{fig:sensitivity_discount}
\end{figure}

Figure \ref{fig:sensitivity_discount} compares the performance of three leading methods---C-SAA-CLF, C-RO-B-CLF, and C-RO-D-CLF across average and worst-case realized objective, fulfillment cost, and consolidation level. Consolidation level is defined as the percentage of units batched (sourced and shipped by the same location–carrier pair) in each order.

As shown in Figure \ref{fig:sensitivity_discount_objVal}, deeper discounts lead each method to exploit batching incentives more aggressively.
C-SAA-CLF consistently attains the lowest average objective across
all discount factors, declining from approximately 0.875 at $\beta = $0.1 to about 0.785 at $\beta = $0.7. C‑RO‑B‑CLF follows closely, from 0.895 down to 0.805, while C‑RO‑D‑CLF remains the most conservative, decreasing from 0.955 to 0.905.  

The worst‑case objective exhibits a more differentiating pattern (Figure \ref{fig:sensitivity_discount_objVal}). Both C‑RO‑D‑CLF and C‑RO‑B‑CLF are less susceptible to change in the discount, whereas C‑SAA‑CLF starts at 0.842 and falls to 0.802 under high discounts.

\revision{
As expected, all three
methods exhibit lower average objective values due to lower fulfillment costs
as the discount factor increases (see Figure
\ref{fig:sensitivity_discount_fulfillment_cost}). This aligns with the
intuition that higher discounts encourage more cost-efficient
fulfillment strategies. The relative decrease in fulfillment costs across the methods can be attributed to their higher consolidation levels as the discount factor increases, as shown in Figure \ref{fig:sensitivity_discount_consolidate_pct}. Notably, C-SAA-CLF responds most sharply, with its consolidation rate rising the most rapidly as the discount grows, which is consistent with its pronounced improvement in worst‑case objective performance.
}

\revision{

\subsection{Perturbation Analysis under Varying Simulation Environments}
\revisiontwo{
A critical challenge in any simulation-based evaluation is that the simulator itself is an imperfect model of the real-world data-generating process. This is particularly true in this application, where the purely observational nature of the historical data, and the resulting lack of counterfactual outcomes, make it impossible to construct a perfect, clairvoyant simulator. Therefore, to rigorously assess the stability of the proposed solutions, this section performs a stress test designed to probe the framework's performance under a misspecified ground truth. This analysis is conducted by introducing systematic perturbations to the simulated delivery time deviation distributions, representing moderate to adversarial divergences from the test data. Specifically, two types of perturbations are considered: random uniform noise and a delay-specific bias applied to the simulation probability distributions. Together, these tests provide critical insight into the methods' stability under real-world distributional shifts and their robustness to the inevitable simulator-reality mismatch.
}

\paragraph{Random Uniform Noise Perturbation}
Let the original simulation distribution for each test instance be denoted by $P^{sim} =\{p_i\}_{i=1}^C$. For each class $i \in [C]$, a random noise term $\delta_i \sim \mathcal{U}(-\epsilon, \epsilon)$ is added to the probability $p_i$, resulting in:
\(
\tilde{p}_i = p_i + \delta_i.
\)
Negative values are clipped: $\tilde{p}_i = \max(\tilde{p}_i, 0)$, and the resulting vector is normalized to form a valid probability distribution: $\tilde{p}_i =  \tilde{p}_i / \sum_j \tilde{p}_j$.

Figure \ref{fig:uniform_perturb} reports the average and worst‐case realized objectives as the noise level \(\epsilon\) increases. In the average case (Figure \ref{fig:uniform_noise_realized_obj_normalized}), both C‑RO variants degrade more gradually than C‑SAA, with C‑RO‑B‑CLF becoming the top performer beyond \(\epsilon=10\). In the worst‐case scenario (Figure \ref{fig:uniform_noise_worst_case_obj_normalized}), while all methods experience performance degrade, the ranking remains unchanged: C‑RO‑D achieves the lowest realized objective, followed by C‑RO‑B, and then C‑SAA. These results indicate the superior robustness of the C‑RO methods under uniform noise.

\paragraph{Delay-Specific Bias Perturbation}
Define the set of late‐delay classes \(C_{\mathrm{late}} = \{\,i : \xi_i > 0\}\). A multiplicative bias \(\eta\) is applied to those classes:
\(\tilde{p}_i = p_i \cdot (1 + \eta).
\)
The probabilities for classes not in $C_{\mathrm{late}}$ remain unchanged. The resulting vector is clipped to the interval $[0, 1]:$ $\tilde{p}_i = \max(\min(\tilde{p}_i, 1), 0)$, followed by normalization.

Figure \ref{fig:bias_perturb} shows performance metrics for bias levels up to 250\%. In the average case (Figure \ref{fig:delay_bias_realized_obj_normalized}), all CSO methods demonstrate improved performance relative to the Greedy baseline, which deteriorates rapidly. This trend is evidenced by the overall decrease in normalized objective values as the bias parameter $\eta$ increases. Notably, C‑RO‑B overtakes C‑SAA once \(\eta\) exceeds 150\%, highlighting its robustness under adversarial delay distortions.  In the worst‐case objective (Figure \ref{fig:delay_bias_worst_case_obj_normalized}), C‑RO‑D consistently achieves the best performance, followed by C‑RO‑B‑CLF. Overall, the C‑RO methods exhibit lower variability in performance compared to C‑SAA, indicating enhanced stability under adversarial delay conditions.

\begin{figure}[!ht]
\centering
        \begin{subfigure}[b]{0.48\textwidth}
        \centering
        \includegraphics[width=\textwidth]{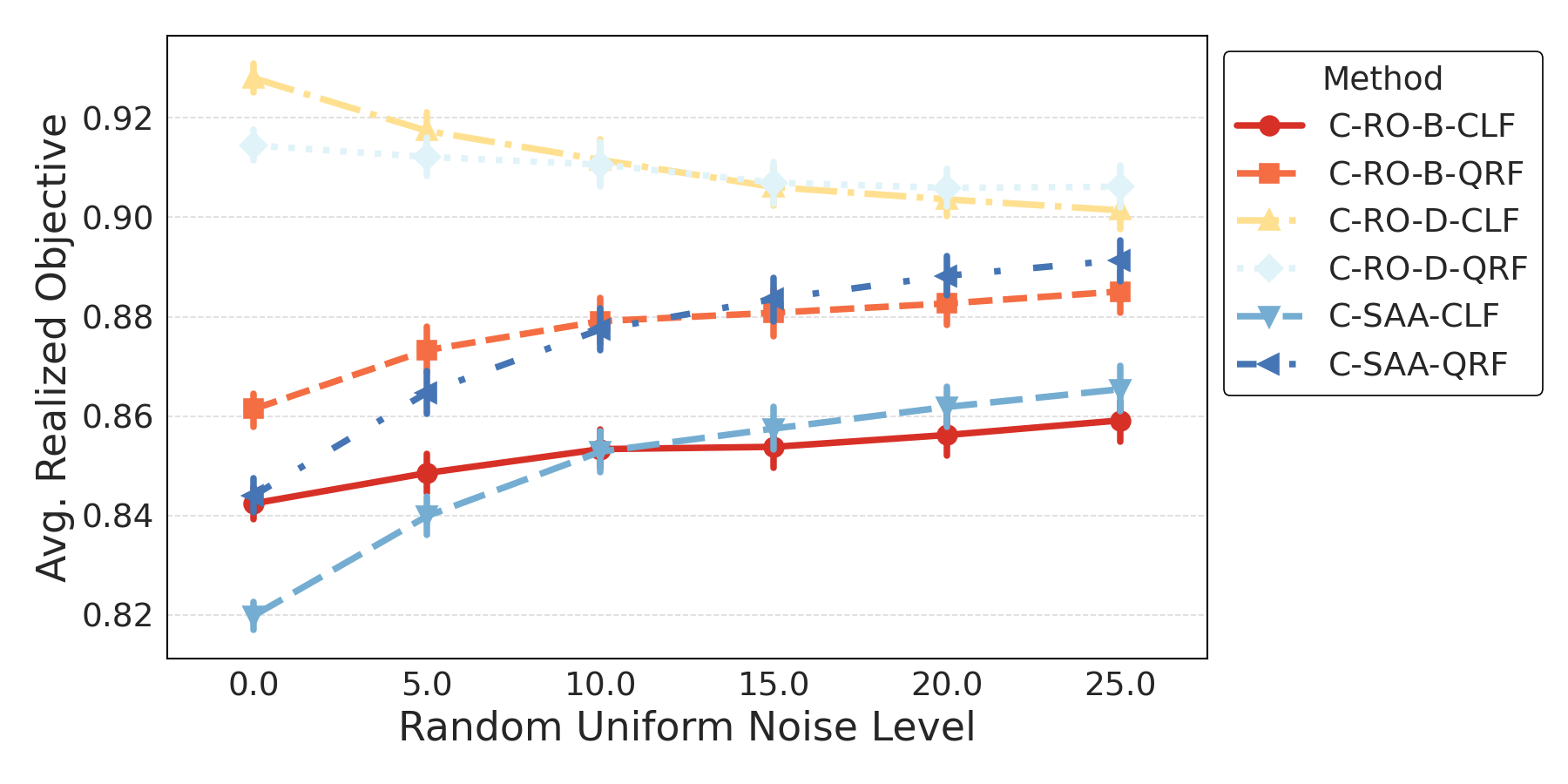}
        \caption{Change in Average Objective.}
    \label{fig:uniform_noise_realized_obj_normalized}
    \end{subfigure}   
    \hfill
    \begin{subfigure}[b]{0.48\textwidth}
        \centering
        \includegraphics[width=\textwidth]{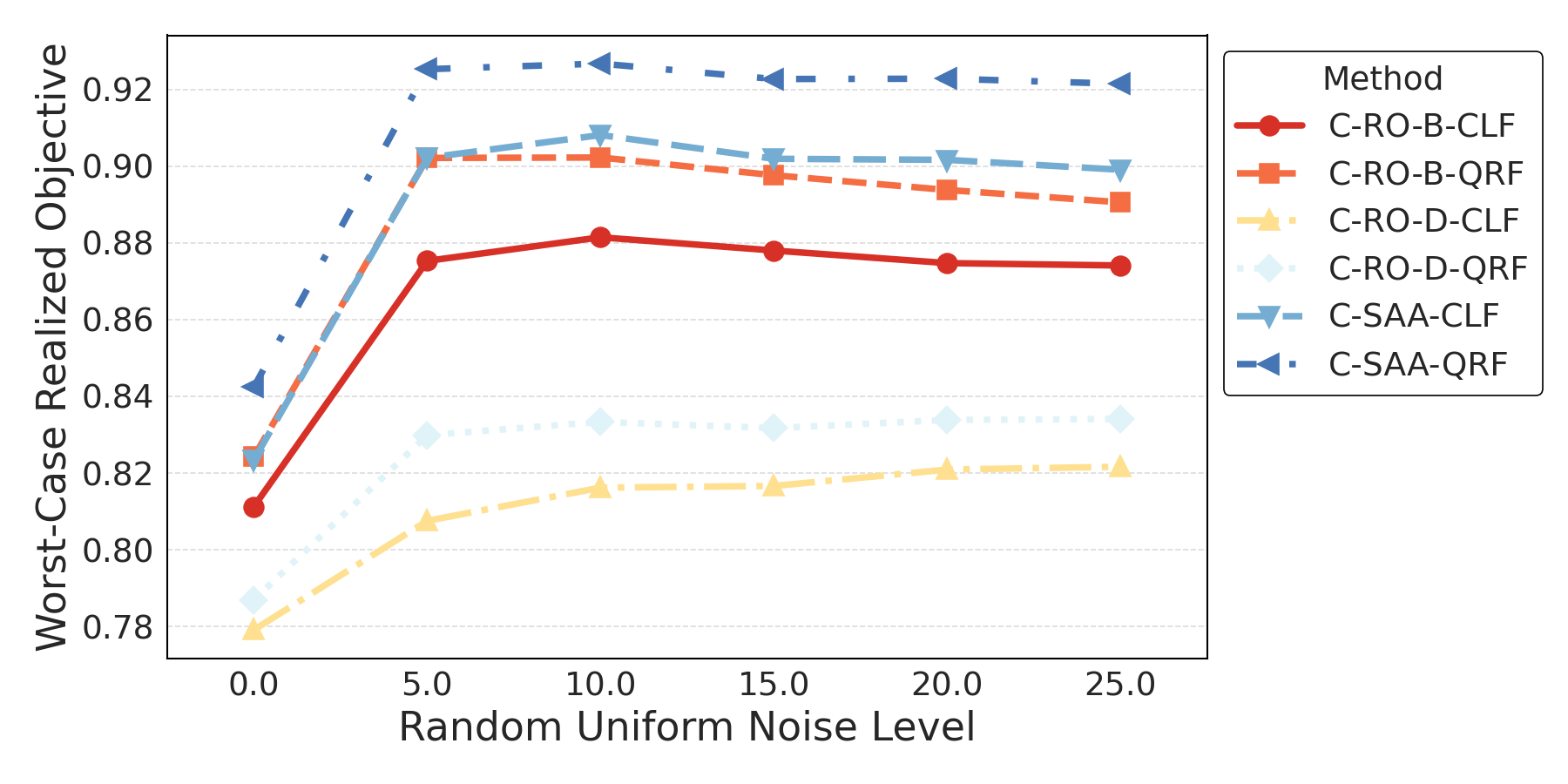}
        \caption{Change in Worst-Case Objective.}      \label{fig:uniform_noise_worst_case_obj_normalized}
    \end{subfigure}    
        \caption{Performance Across Different Random Uniform Noise Perturbations.}
\label{fig:uniform_perturb}
\end{figure}

\begin{figure}[!ht]
\centering
        \begin{subfigure}[b]{0.48\textwidth}
        \centering
        \includegraphics[width=\textwidth]{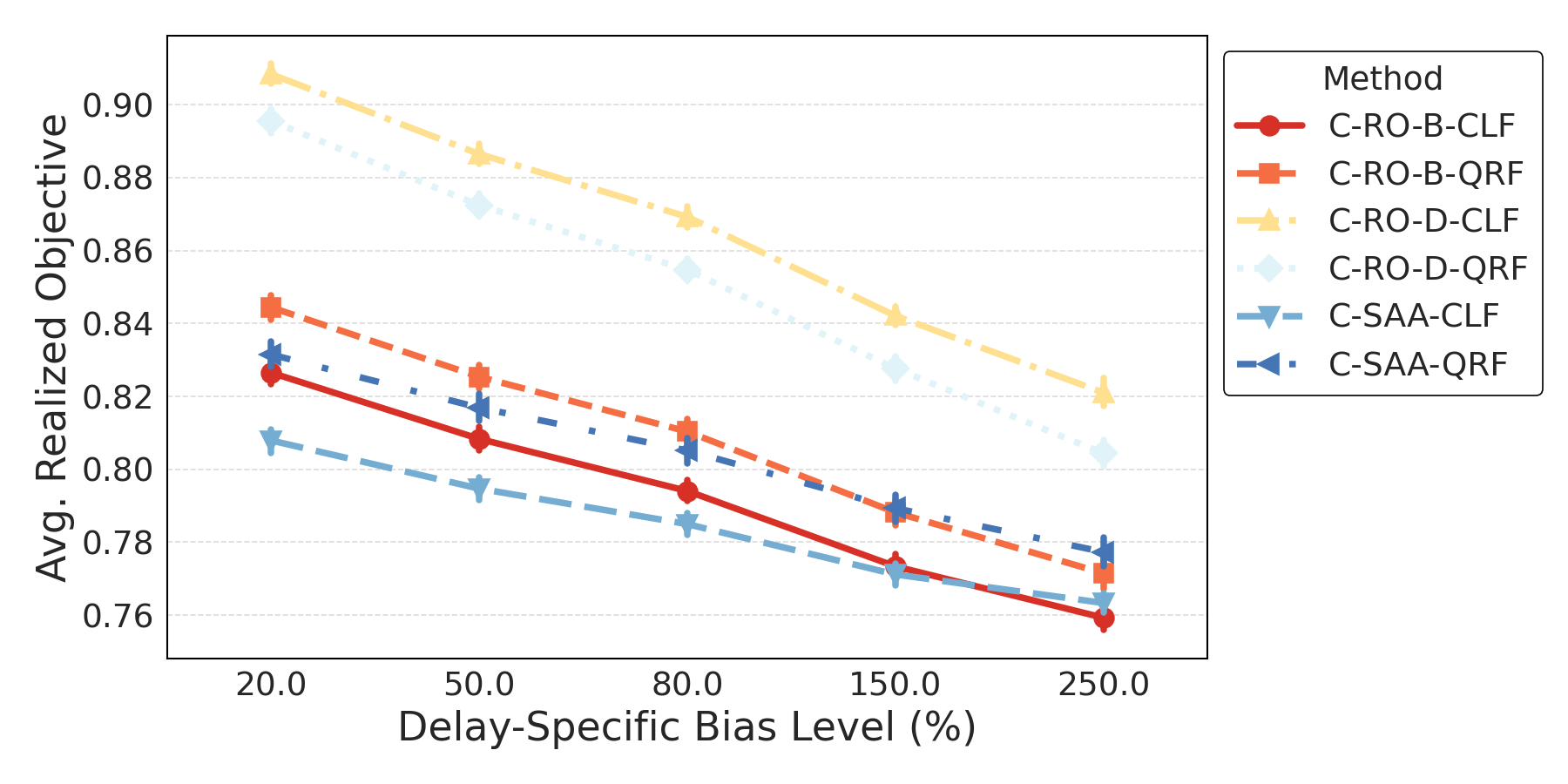}
        \caption{Change in Average Objective.}
    \label{fig:delay_bias_realized_obj_normalized}
    \end{subfigure}   
    \hfill
    \begin{subfigure}[b]{0.48\textwidth}
        \centering
        \includegraphics[width=\textwidth]{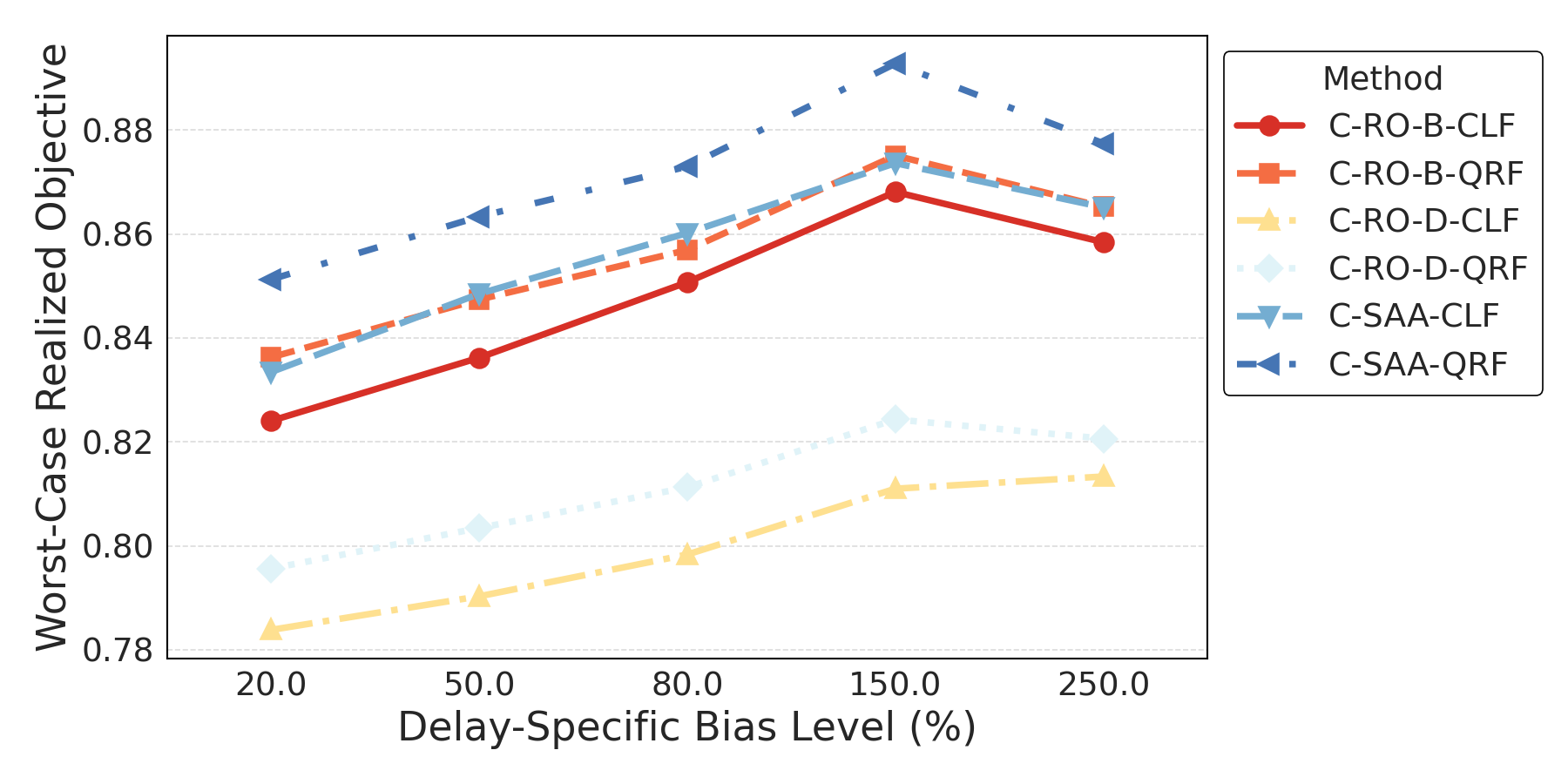}
        \caption{Change in Worst-Case Objective.}      \label{fig:delay_bias_worst_case_obj_normalized}
    \end{subfigure}    
        \caption{Performance Across Different Delay-Specific Bias Perturbations.}
\label{fig:bias_perturb}
\end{figure}

\subsection{Computational Scalability of the CSO Methods}

The computational scalability of the proposed framework was evaluated by measuring both the change in computation time (model construction + solve time) and solution quality as instance size grows. Each CSOFP instance is characterized by the number of items $|\mathcal I|$ and the number of location–carrier pairs (fulfillment options) $|\mathcal K|\times|\mathcal L|$. In practice, the eligibility constraint (Constraint~(\ref{eq:eligibility})) is used to prune infeasible location–carrier pairs in advance. Then, 
\(
  (\mathcal K',\mathcal L')
  \;=\;
  \bigl\{(k,\ell)\in\mathcal K\times\mathcal L \,\bigm|\, \exists\,i\in\mathcal I:\;e_{ik\ell}=1\bigr\}
\)
yields the final set of eligible fulfillment options and thus the effective problem size.

Figure~\ref{fig:solve_time_vs_eligible_options} shows computation time as the number of eligible fulfillment options increases. Only the MC‑CLF variants of C‑RO and C‑SAA are displayed, since their QRF counterparts exhibit nearly identical runtimes. All three methods achieve optimality on every instance within minutes, confirming their practical applicability. In particular, C‑RO‑B solves most instances in under one second (never exceeding three seconds), and C‑RO‑D solves most in under ten seconds.

These timing differences reflect each method’s formulation. Compared to the nominal problem, the robust counterpart in C‑RO‑B introduces only continuous variables whose count grows linearly with $|\mathcal K'|$ and $|\mathcal L'|$. By contrast, C‑RO‑D adds $N$ extra constraints from sampling, and C‑SAA requires solving $Q$ independent problems whose objective function complexity scales with the sample size $N_{1}$. Parallel computation of the $Q$ subproblems can further reduce C‑SAA’s runtime.

\revisiontwo{
For completeness, a comparison of computation times against the baseline models is provided (see Appendix K). The Greedy and PTO methods are computationally cheaper, solving most instances in fractions of a second. However, the proposed CSO methods, particularly the highly-efficient C-RO-B, remain well within practical time limits for real-world deployment (Figure~\ref{fig:solve_time_vs_eligible_options}). The substantial improvements in solution quality (Table~\ref{tab:methods_comparison}) justify this modest increase in computational overhead, representing a favorable trade-off between performance and speed.
}

Since single-item orders account for the majority of order volume, multi-item orders (what the methods are designed for) can be routed to a dedicated fulfillment optimization engine without disrupting standard operations. Under normal load, any of the three methods may be applied; under peak loads—such as holiday surges—C‑RO‑B is recommended to preserve fast solve times while maintaining high solution quality. Concurrent batches of orders can be processed in parallel when computational resources permit.

Figure~\ref{fig:obj_scalability} reports average and worst-case objective values for small ($\leq 400$ eligible options), medium ($\leq 1000$ eligible options), and large ($\leq 3000$ eligible options) instances. Average performance remains uniformly high across all size categories, whereas worst-case deviations grow with instance size. Collectively, these results demonstrate near-linear scalability in computation time, consistent high solution quality, and robustness across problem sizes, providing empirical evidence of the framework’s practical scalability.

\begin{figure}[!ht]
\centering
    \begin{subfigure}[b]{0.48\textwidth}
        \centering
        \includegraphics[width=\textwidth]{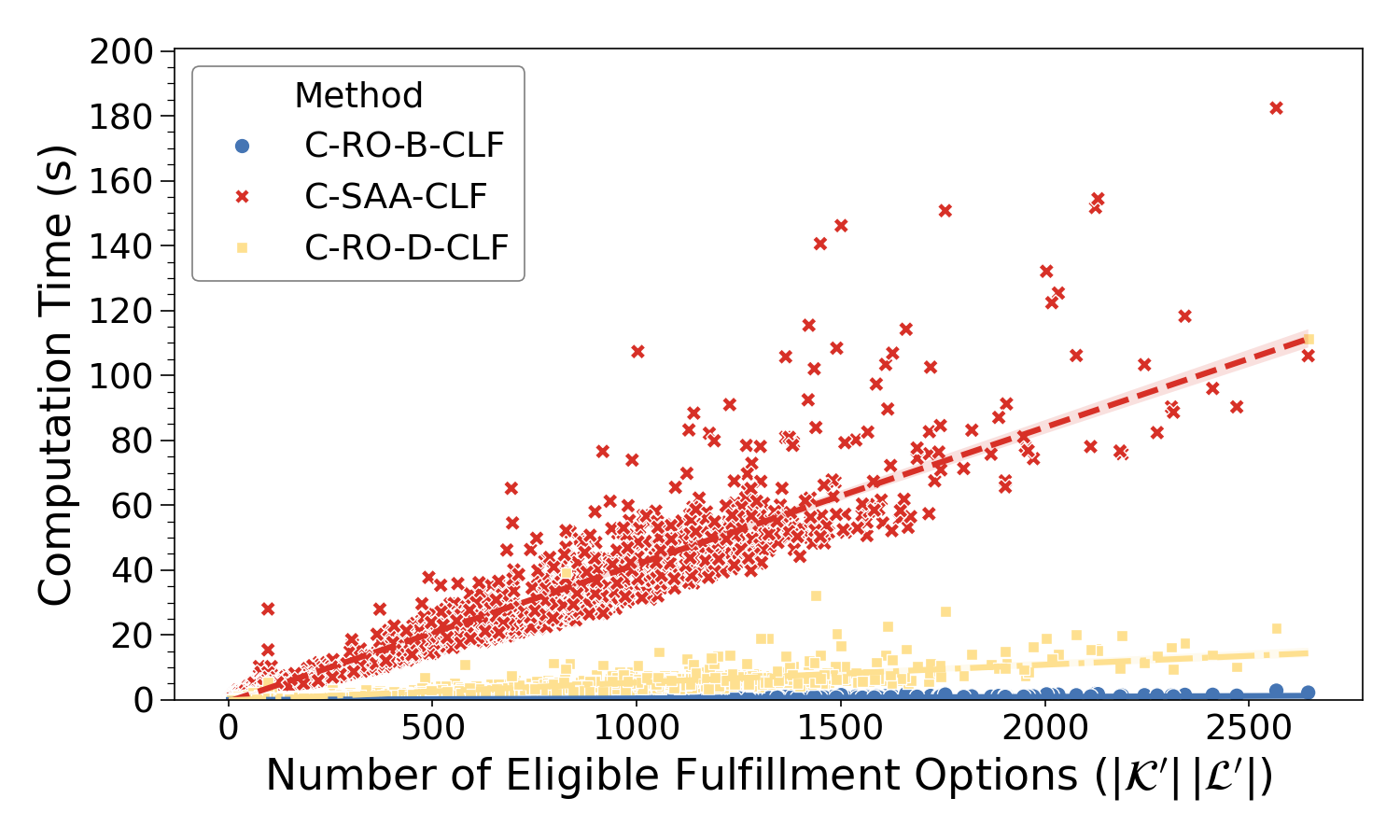}
        \caption{Change in the Computation Time vs. the Number of Eligible Fulfillment Options.}
        \label{fig:solve_time_vs_eligible_options}
    \end{subfigure}
    \hfill
    \begin{subfigure}[b]{0.48\textwidth}
        \centering
        \includegraphics[width=\textwidth]{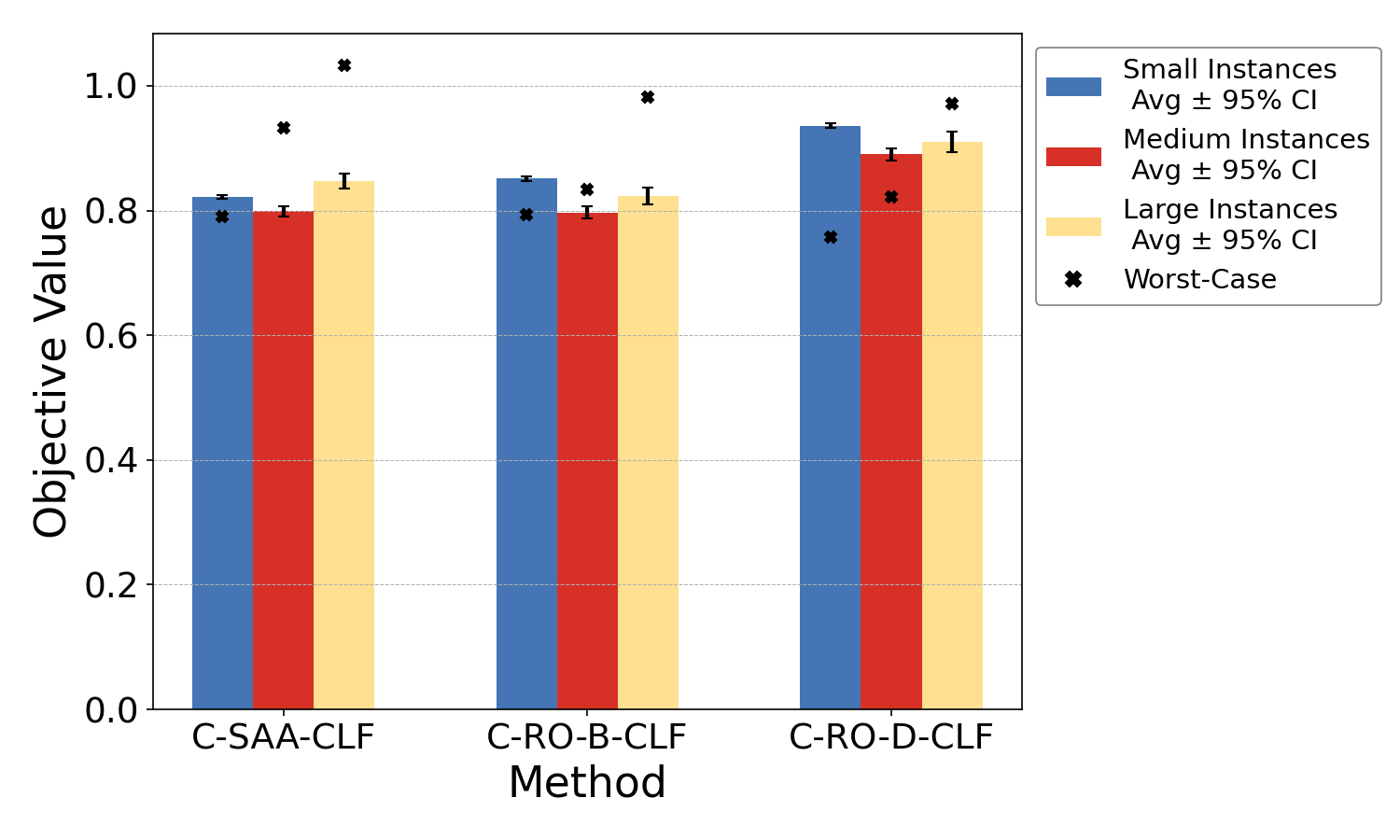}
        \caption{Objective Value by Instance Categories.}
        \label{fig:obj_scalability}
    \end{subfigure}
    \caption{Computational Scalability Analysis.}
    \label{fig:combined}
\end{figure}

\subsection{Managerial and Operational Implications}

The proposed fulfillment framework offers actionable guidance for practitioners seeking to align fulfillment strategies with demand variability, service-level requirements, and computational constraints. Depending on the operating context, different methods are preferable. When demand is stable and the goal is to minimize average fulfillment cost, C-SAA-CLF provides the most cost-effective solution. Under moderate runtime constraints, C-RO-B-CLF offers a balanced trade-off between efficiency and robustness. In contrast, during high-stakes periods such as holiday promotions or flash sales—where the cost of service-level violations is high—C-RO-D-CLF and C-RO-B-QRF (with 99\% coverage) are better suited to mitigate tail risks and ensure reliability. This scenario-driven approach allows e-commerce platforms to dynamically adapt fulfillment strategies to shifting conditions while leveraging a unified optimization infrastructure. As retailers increasingly face variable demand patterns and rising customer expectations, the ability to flexibly match fulfillment strategy to context becomes a key operational advantage.

Empirical evaluations indicate that the proposed methods can reduce combined fulfillment cost and delivery-timeliness penalties by an average of up to 18\% per order under the default settings. To illustrate the potential impact, consider a retailer  operating at a scale comparable to the industrial partner, with per-order costs ranging from \$5 to \$30 and monthly volumes of 0.1 to 1 million multi-item orders. Under these scenarios, the corresponding annual savings could range from several million to tens of millions of dollars. Moreover, this framework can be readily extended to single-item orders, which could yield further savings with minimal incremental investment.
}

\section{Conclusion} \label{sec:conclusion}
\revision{
This paper introduces a generic contextual stochastic optimization (CSO)
framework for data‑driven decision making when only observational data, hence missing counterfactual costs, are available. The framework is deliberately modular.  
First, any algorithm capable of returning a full predictive distribution can serve as the \textit{contextual distribution oracle}. In the empirical study, calibrated multi‑class classifiers and tree‑based quantile‑regression forests were selected to capture the discrete, ordered nature of delivery‑time deviations.  
Second, the oracle feeds directly into two tractable optimization engines: (i) a \emph{Contextual Sample‑Average Approximation} (C‑SAA) for risk‑neutral objectives and (ii) a \emph{Contextual Robust Optimization} (C‑RO) formulation for risk‑averse objectives.  
Because the learning and optimization layers communicate only through the oracle, the template can be readily translated to other operations‐management settings.

When applied to an omnichannel, multi‑courier order‑fulfillment problem, the CSO framework reduced the combined cost–service objective by up 18\% on average and up to 22\% in the simulated worst case on an industry‑scale data set, implying annual savings from several to tens of millions of dollars. This framework allows practitioners to tune robustness---raising it during promotional peaks and relaxing it in steadier periods---while still exploiting item‑consolidation discounts and heterogeneous carrier performance in the dynamic e-commerce landscape.

While the paper delivers an end‑to‑end prototype with demonstrable impact, several avenues remain open:

\begin{itemize}
    \item Performance guarantees under partial feedback: Establish theoretical performance guarantees, such as worst-case regret and post-decision surprise, for CSO methods under the partial‑feedback setting. Linking causal‑inference techniques with SAA and RO theory can be a promising direction, thereby quantifying the value of additional data when counterfactual costs remain unobserved.

    \item Hybrid distributional forecasters: Develop distributional models that fuse the calibration strength of classification with the sharp tails of quantile regression, yielding tighter and more informative uncertainty sets.

    \item End‑to‑end predict‑then‑prescribe learning for observational data: 
    Incorporate downstream optimization loss directly into training, extending recent ``integrated predict‑then‑optimize” approaches  \citep{qi2021integrated, elmachtoub2022smart} to settings where data are purely observational.

    \item Multi‑period, inventory‑aware CSO:
    Extend the framework to jointly optimize fulfillment, shipping, and stocking over rolling horizons, capturing inventory replenishment spillovers and elevate the model from a myopic to a long-term strategic planning tool.

    \item \revisiontwo{Generalization to broader problem classes: Extending the framework to broader settings like network flow optimization is a key direction for future work. The sampling-based C-SAA and C-RO-D approaches can be directly generalized. However, the C-RO-B formulation's tractability relies on the specific problem structure, and adapting it to more complex constraints would require new reformulations.}
\end{itemize}
}

\clearpage
\begin{APPENDIX}{}

\section{MILP Model Notations} \label{sec: notations}
Table~\ref{tab: opt-nomenclature} presents the notations used in the MILP optimization model.
\begin{table}[!ht] 
    \centering
    \caption{Optimization Model Nomenclature}
 \label{tab: opt-nomenclature}    
 \begin{tabular}{lp{10cm}}
    \hline
    \textbf{Sets} & \\
    $\mathcal{I}$ & Set of SKUs, indexed by $i$.\\
    $\mathcal{K}$ & Set of carriers, indexed by $k$. \\
    $\mathcal{L}$ & Set of locations, indexed by $\ell$. \\ 
    \hline
    \textbf{Parameters} & \\
    $\operatorname{inv}_{i\ell}$ & Inventory level of SKU $i$ at location $\ell$ when the order is placed.\\
    $\operatorname{cap}_{\ell}$ & Available capacity at location $\ell$ when the order is placed. \\
    $e_{ik\ell}$ & 1 if (carrier $k$, location $\ell$) is an eligible candidate location-carrier pair for SKU $i$ in the order and 0 otherwise. \\
    $c_{ik \ell}^{ship}$ & Per-unit shipping cost when using location $\ell$ and carrier $k$ to ship SKU $i$. \\
    $c_{\ell}^{fixed}$ & Per-unit fixed sourcing (non-parcel) cost when using location $l$. \\
    $\beta_{k }$ & Discount factor for shipping multiple units using carrier $k$. \\
    $q_{i}$ & Quantity of SKU $i$ in the order. \\ 
    $d_{k \ell}$ & Delivery time deviation when using location $\ell$ and carrier $k$. \\
    $\gamma^+$ & Constant that converts late delivery penalty into a per-unit cost. \\
    $\gamma^-$ & Constant that converts early delivery penalty into a per-unit cost. \\
    $N$ & A large number = $\sum_{i \in \mathcal{I}}q_{i} + 2$. \\
    $M_{k \ell}$ & A large number = $\sum_{i\in \mathcal{I}} c_{ik \ell}^{ship} q_{i}$.\\
    \hline
    \textbf{Decision Variables} & \\
    $z_{ik \ell}$ & Integer decision variable which measures the quantity of SKU $i$ in the order sourced from location $\ell$ and shipped by carrier $k$. \\
    $y_{k \ell}$ & Binary decision variable which equals 1 if more than one unit is sourced from location $\ell$ and carrier $k$ within the order, and equals 0 otherwise. \\
    $w_{k \ell }$ & Continuous decision variable which measures the cost to be discounted in the order using location-carrier pair $(k, \ell)$. \\
    \hline
    \end{tabular} 
\end{table}

\section{Data Preprocessing}

To optimize the dataset for machine learning and optimization models,
several preprocessing steps were implemented to clean the
data. Duplicate entries and rows with missing values were
removed. Geographic coordinates were generated from the zip codes of
the fulfillment centers and customer locations, with any invalid zip
codes discarded in the process. The geographic distance between
customer and fulfillment center locations was approximated using the
Haversine formula. Additionally, the dataset was filtered to include
only carriers with sufficient order history for training prediction
models, resulting in a subset of carriers.

\section{Descriptive Analytics and Insights}

The company receives tens of thousands of home delivery online orders
daily, with an average ranging from 30,000 to 70,000 orders. During
peak hours, 500 to 2,000 orders may be received within a 15-minute
time window. Direct information regarding which carriers are available
to ship specific SKUs from particular fulfillment centers is not
available in the dataset. However, heuristic rules can be derived
through analyses of historical orders. The customer’s location and the
dimensions/weight of the SKU(s) in an order are critical factors in
determining carrier suitability. Some carriers show a more scattered
distribution of customer locations, indicating limited coverage in
certain states of the U.S.. Additionally, coverage is dependent on the
location of the fulfillment center from which the order is shipped, as
some carriers can only serve a limited range around that center. In
contrast, other carriers have served locations spanning nearly the
entire US. Furthermore, the upper bounds of SKU length, height, and
weight in the historical data reflect the package size and weight
limits of different carriers.

\begin{figure}[!ht]
    \includegraphics[width=0.33\textwidth]{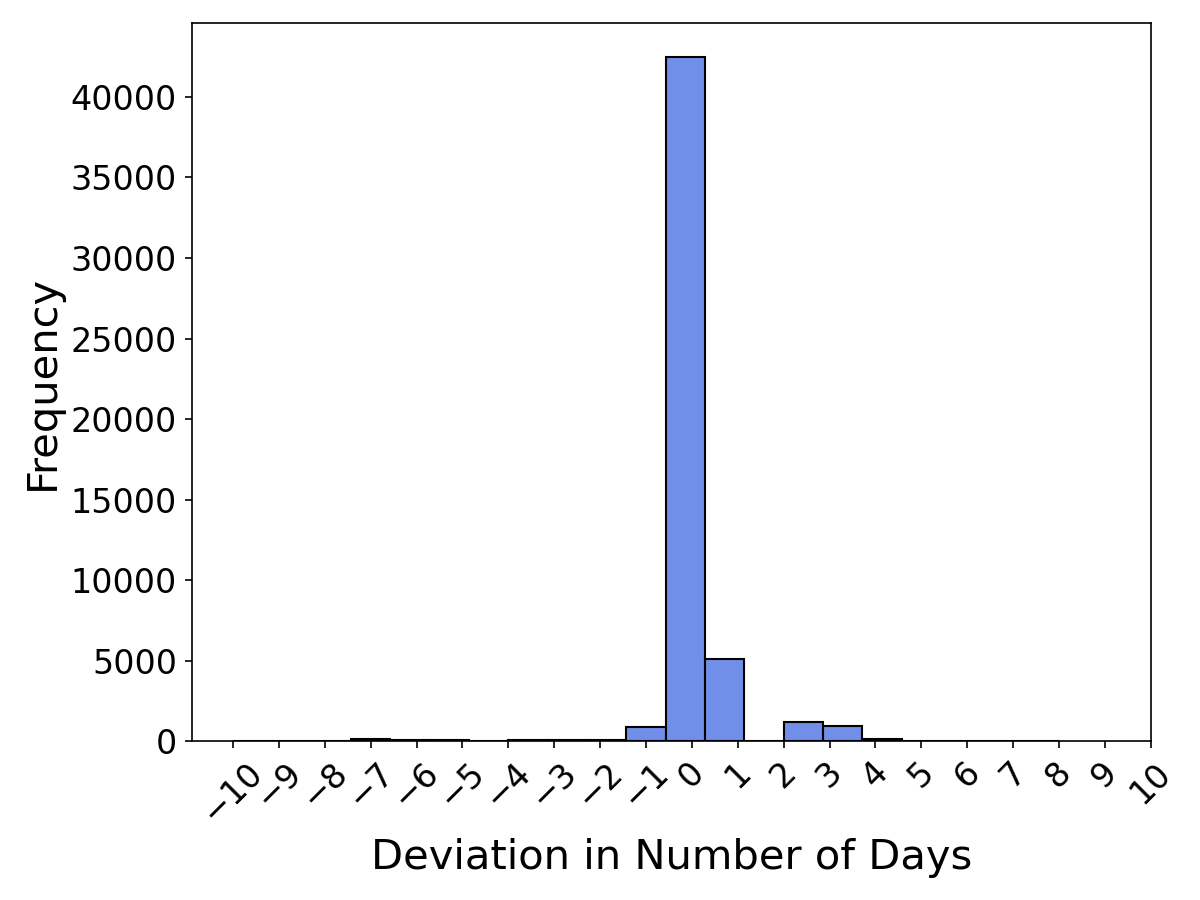} 
    \includegraphics[width=0.33\textwidth]{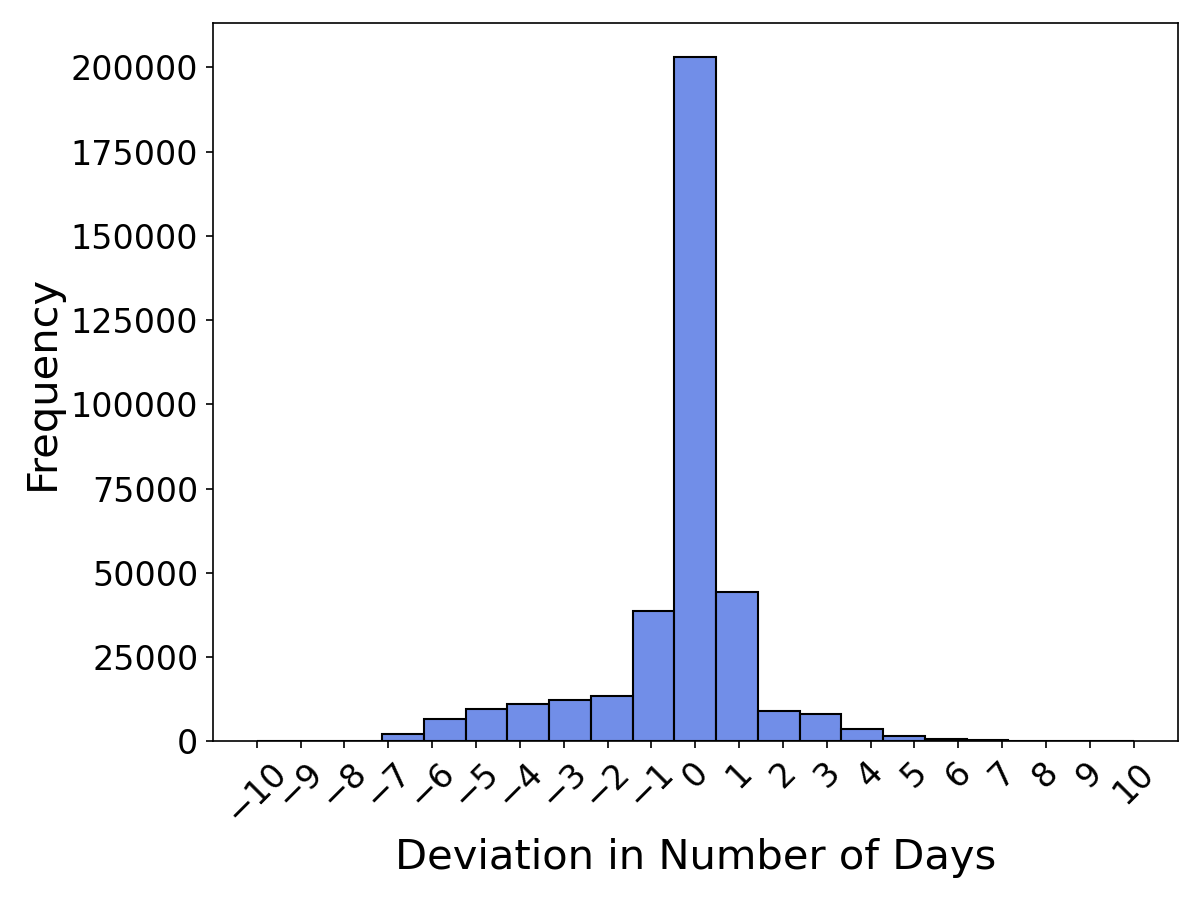}
    \includegraphics[width=0.33\textwidth]{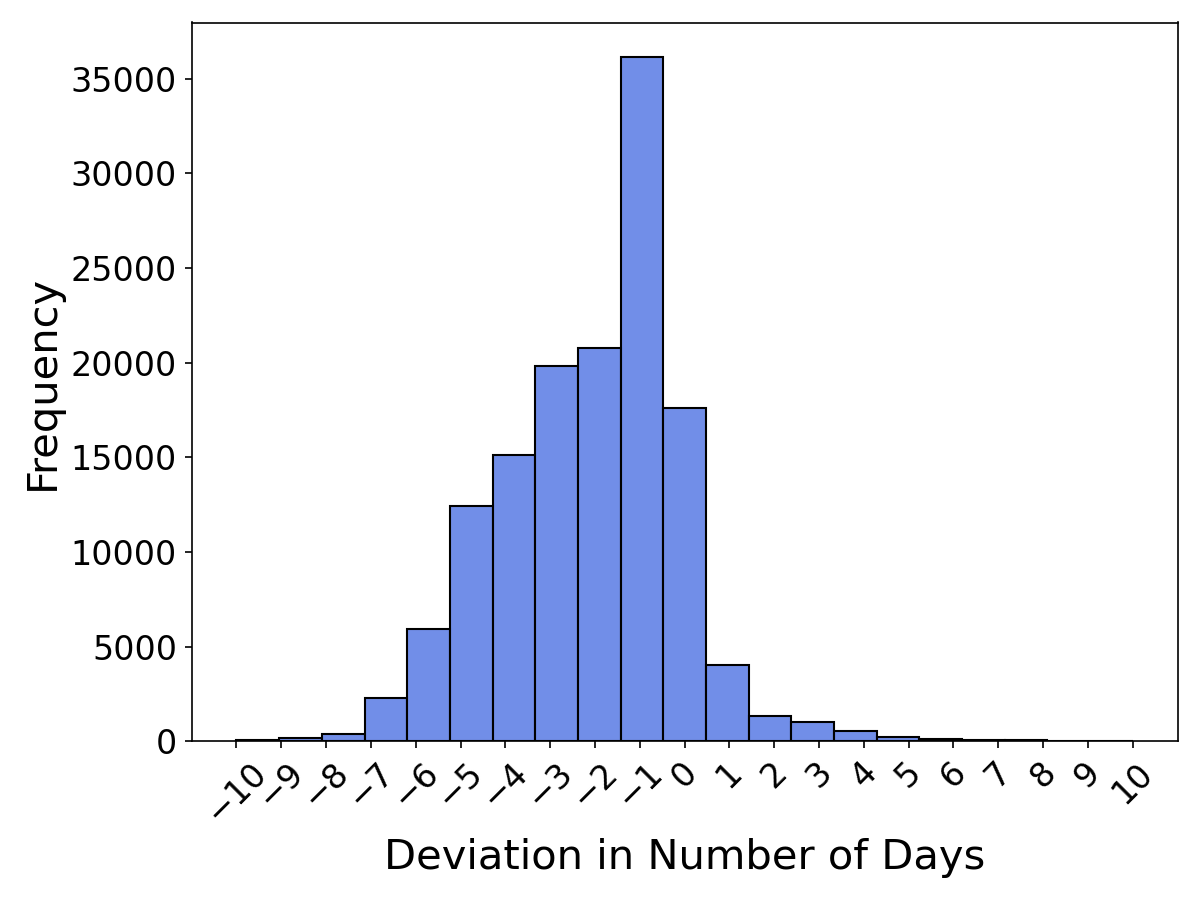}
    \caption{Delivery Deviation Distribution of Orders Shipped by Three Example Carriers.}
    \label{fig:delay_dist}
\end{figure}

The distribution of deviations varies significantly across carriers, as
illustrated in Figure \ref{fig:delay_dist}. For the first two
carriers, the majority of orders are delivered on time. The first
carrier shows nearly all on-time deliveries, with the remaining orders
arriving late. The second carrier has approximately half of the orders
delivered on time, with slightly more early deliveries than late
ones. In contrast, orders shipped by the last carrier predominantly arrive one
day early, and most deliveries arrive ahead of schedule. Given these
pronounced differences, separate prediction models were developed for
each carrier.

\section{Feature Importance}
To identify the most
important features for predicting deviations, a feature importance
analysis was conducted using a tree-based regression model. Two
methods were examined: Mean Decrease in Impurity (MDI) and Shapley
values (SHAP). MDI measures the importance of a feature by calculating
the total reduction in impurity (e.g., Gini impurity or entropy) it
brings when used in decision trees \citep{breiman2001random}. SHAP
values, on the other hand, provide a unified measure of feature
importance based on cooperative game theory, attributing to each
feature its contribution to the prediction by considering all possible
feature combinations \citep{lundberg2020local2global}.

Figure \ref{fig:feat_imp} highlights the top 10 most important
features evaluated, as evaluated using MDI and SHAP, for an example
carrier under a tree-based regression model. For this carrier, the
planned lead time (promised delivery time) emerges as the most significant
factor. Additionally, order information such as the release time and
SKU type, fulfillment center characteristics like on-hand capacity,
and carrier-specific attributes such as promised time-in-transit and
carrier zone also play important roles.

\begin{figure}[!ht]
\centering
    \begin{subfigure}{0.49\textwidth}
        \centering
        \includegraphics[width=\textwidth]{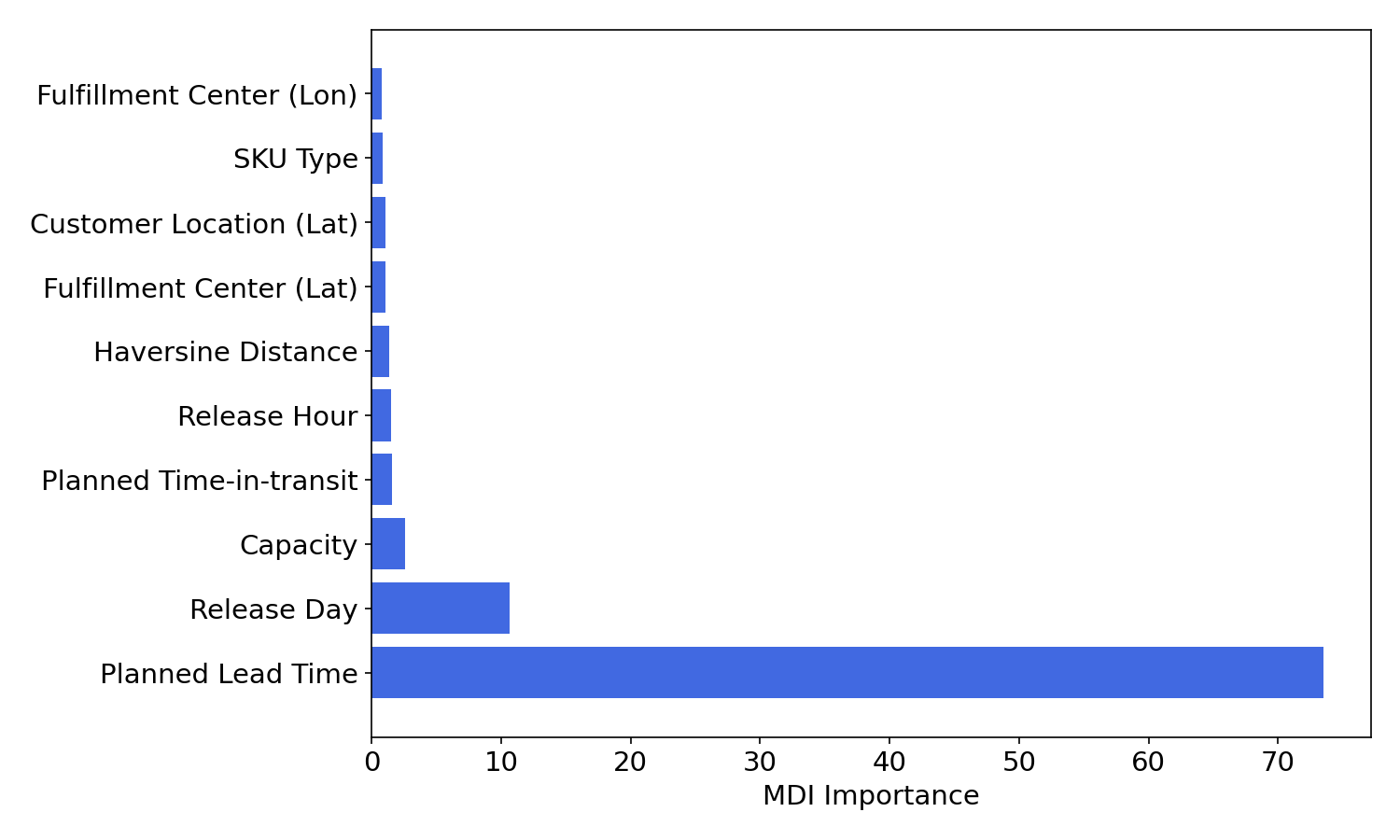}
        \caption{MDI Feature Importance.}
    \end{subfigure}
    \begin{subfigure}{0.49\textwidth}
        \centering
        \includegraphics[width=\textwidth]{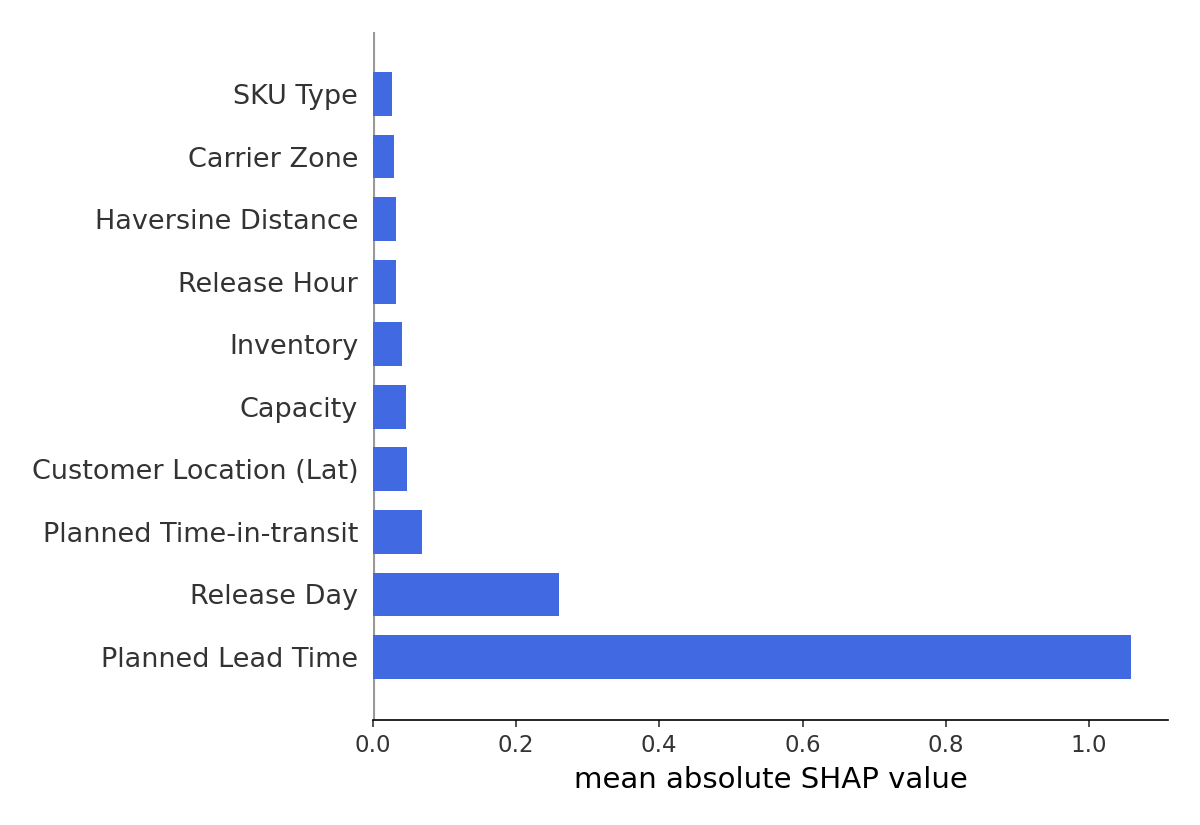}
        \caption{SHAP Value Feature Importance.}
    \end{subfigure}
    \caption{Top 10 Mean Regression Features for Delivery Deviation Measured by Two Feature Importance Metrics For an Example Carrier.}
    \label{fig:feat_imp}
\end{figure}

\section{Linearize the Objective Function of the Order Fulfillment Problem} \label{sec: linearize}

The nonlinear objective function $g(\mathbf{z}, \mathbf{d})$ can be linearlized using auxiliary variables. Let $\mathbf{y} = (y_{k \ell })_{k \in \mathcal{K}, \ell \in \mathcal{L}} \in \{0, 1\}^{(|\mathcal{K}||\mathcal{L}|)}$ be indicator variables, where $y_{k \ell} = 1$ indicates more than one unit is sourced from location $\ell$ and shipped by carrier $k$ within the order. Additionally, let $\mathbf{w}= (w_{k \ell})_{k \in \mathcal{K}, \ell \in \mathcal{L}} \in \mathbb{R}^{(|\mathcal{K}||\mathcal{L}|)+}$ be continuous variables, where $w_{k \ell}$ represents the shipping costs to be discounted in the order using carrier-location pair ($k, \ell$).

The linear objective function and the additional constraints for linearization are defined as follows:
\begin{equation}
   \begin{aligned}
        g(\mathbf{z}, \mathbf{d}) = & (\sum_{i \in \mathcal{I}} \sum_{k \in \mathcal{K}}  \sum_{\ell \in \mathcal{L}} (c_{\ell}^{fixed} + c_{ik \ell}^{ship} +  \gamma^+ d_{k \ell}^+ + \gamma^- d_{k \ell}^-) z_{ik \ell} - \sum_{k \in \mathcal{K}}  \sum_{\ell \in \mathcal{L}}  \beta_{k }  w_{k \ell}) \\
       \end{aligned}
\end{equation}
\begin{equation}
   \begin{aligned}
       \text{s.t.}  \quad 
       & \sum_{i \in \mathcal{I}} z_{ik \ell}  \geq 2 - N  (1- y_{k \ell})\\
       & \sum_{i \in \mathcal{I}} z_{ik \ell}  \leq 1 + N  y_{k \ell}\\
       & w_{k \ell} \leq M_{k \ell}  y_{k \ell} \\
       & w_{k \ell} \leq \sum_{i \in \mathcal{I}} c_{ik \ell}^{ship} z_{ik \ell} + (1 - y_{k \ell}) M_{k \ell} \\
       & w_{k \ell} \geq \sum_{i \in \mathcal{I}} c_{ik \ell}^{ship} z_{ik \ell} - (1 - y_{k \ell})  M_{k \ell} \\
       & y_{k \ell} \in \{0, 1\} \\
       & w_{k \ell} \geq 0 \\
       \end{aligned}
       \quad , \forall {k \in \mathcal{K}}, {\ell \in \mathcal{L}}, \\
\end{equation}
where $N$ and $M_{k \ell}$ are large constants for the big-M constraints, with their strictest values being  $\sum_{i \in \mathcal{I}}q_{i} + 2$ and $\sum_{i\in \mathcal{I}} c_{ik \ell}^{ship} q_{i}$, respectively.

\section{Applying the C-RO-B to CSOFP}

To successfully apply the budgeted interval uncertainty set, the uncertainty vector in the objective function $g(\mathbf{z}, \tilde{\mathbf{d}})$ should be element-wise positive. However, in the case of CSOFP, this requirement was not met due to the asymmetric penalty costs associated with deviations. Instead, the uncertainty in the problem can be represented by the vector of total costs. Let this uncertain cost vector be $\tilde{\mathbf{c}}= (\tilde{c}_{ik \ell})_{i \in \mathcal{I}, k \in \mathcal{K}, \ell \in \mathcal{L}}$, where each $\tilde{c}_{ik \ell} = (c_{\ell}^{fixed} + c_{ik \ell}^{ship} + \gamma^+ \tilde{d}_{k \ell}^+ + \gamma^- \tilde{d}_{k \ell}^-) \geq 0 $ is the sum of fulfillment costs and deviation penalties for the order.

Substituting $\tilde{\mathbf{c}}$ leads to:
\begin{equation} \label{eq: crofp2}
   \begin{aligned}
        \min_{\mathbf{z} \in \mathcal{Z}}  \max_{\tilde{\mathbf{c}} \in \mathcal{U}_b(\mathbf{s})} \quad & (\sum_{i \in \mathcal{I}} \sum_{k \in \mathcal{K}}  \sum_{\ell \in \mathcal{L}} \tilde{c}_{ik \ell} z_{ik \ell} -  \sum_{k \in \mathcal{K}} \sum_{\ell \in \mathcal{L}}  \beta_{k } w_{k \ell}).
   \end{aligned}
\end{equation}

Given that the prediction interval for each $\tilde{d}_{k \ell}$ is $[\underline{d}_{k \ell}, \underline{d}_{k \ell} + \hat{d}_{k \ell}]$, where $\underline{d}_{k \ell}$ and ($\underline{d}_{k \ell} + \hat{d}_{k \ell}$) are the estimated lower and upper quantiles of the deviation for the carrier-location pair $(k, \ell)$, the prediction interval for $\tilde{\mathbf{c}}$ can be expressed as: 
\begin{equation}
[\underline{\xi}_{ik\ell},\underline{\xi}_{ik\ell} + \hat{\xi}_{ik\ell}]= c_{\ell}^{fixed} + c_{ik \ell}^{ship} + 
   \begin{aligned}
       \begin{cases}
           [\gamma^- \underline{d}_{k \ell}^-, \gamma^- (\underline{d}_{k\ell} + \hat{d}_{k\ell})^-] & \text{if } \underline{d}_{k\ell} + \hat{d}_{k\ell} \leq 0 \\
           [0, \max \{ \gamma^- \underline{d}_{k \ell}^-, \gamma^+ (\underline{d}_{k \ell}+ \hat{d}_{k\ell}) \}] & \text{if } \underline{d}_{k \ell} \leq 0, \underline{d}_{k\ell} + \hat{d}_{k\ell} \geq 0 \\
           [\gamma^+ \underline{d}_{k\ell}, \gamma^+ (\underline{d}_{k\ell} + \hat{d}_{k\ell})] & \text{if } \underline{d}_{k\ell} \geq 0 \\
       \end{cases}
   \end{aligned}
\end{equation}

Let $B \in [0, (|\mathcal{K}||\mathcal{L}|)]$ be the uncertainty budget. The contextual budgeted interval uncertainty set
$\mathcal{U}_b(\mathbf{s})$ is expressed as follows:
\begin{equation}
\mathcal{U}_b(\mathbf{s}, B) = \{\bm{\xi} \in \mathbb{R}^{(|\mathcal{I}||\mathcal{K}||\mathcal{L}|)}: \xi_{ik\ell} =  \underline{\xi}_{ik\ell} + \hat{\xi}_{ik\ell} \delta_{k \ell}, \sum_{k \in \mathcal{K}, \ell \in \mathcal{L}} \delta_{k \ell} \leq B, \bm{\delta} \in [0, 1]^{(|\mathcal{K}||\mathcal{L}|)}\}.
\end{equation}

Then, with fixed $\mathbf{z}$, the inner maximization problem of Problem (\ref{eq: crofp2}) becomes:
   \begin{align}
       \max \quad &  \sum_{i \in \mathcal{I}} \sum_{k \in \mathcal{K}}  \sum_{\ell \in \mathcal{L}} \hat{\xi}_{ik \ell} z_{ik \ell} \delta_{k\ell} \\
       \text{s.t.} \quad &  \sum_{k \in \mathcal{K}}  \sum_{\ell \in \mathcal{L}}  \delta_{k\ell} \leq B \label{eq:imp:c1}\\
       & \delta_{k\ell} \in [0, 1], \quad \forall k \in \mathcal{K}, \ell \in \mathcal{L}. \label{eq:imp:c2}
   \end{align}

The robust counterpart of Problem (\ref{eq: crofp2}) is formulated as follows:
\begin{equation} \label{eq: rcp}
\begin{aligned} 
   \min_{\mathbf{z} \in \mathcal{Z}}  \quad & ( \sum_{i \in \mathcal{I}} \sum_{k \in \mathcal{K}}  \sum_{\ell \in \mathcal{L}} \underline{\xi}_{ik \ell} z_{ik \ell} - \sum_{k \in \mathcal{K}}  \sum_{\ell \in \mathcal{L}}  \beta_{k }  w_{k \ell}  + \pi B + \sum_{k \in \mathcal{K}}  \sum_{\ell \in \mathcal{L}} \lambda_{k \ell})\\
   \text{s.t.} \quad & \pi + \lambda_{k \ell} \geq \hat{\xi}_{ik \ell} z_{ik \ell},  \quad \forall i \in \mathcal{I}, k \in \mathcal{K}, \ell \in \mathcal{L} \\
   & \pi \geq 0\\
   & \lambda_{k \ell} \geq 0,  \quad \forall k \in \mathcal{K}, \ell \in \mathcal{L},
\end{aligned}
\end{equation}
where $\pi$ and $(\lambda_{k \ell})_{k \in \mathcal{K}, \ell \in \mathcal{L}}$ are the dual variables of Constraints (\ref{eq:imp:c1}) and (\ref{eq:imp:c2}), respectively.

\section{Probabilistic Multi-Class Classification Using Ensemble Learning} \label{sec: ensemble}

RF leverages multiple
classification trees trained on different random subsets of the data,
making predictions by averaging the results
\citep{breiman2001random}. Let $\mathcal{R}_1^{(t)}, \ldots,
\mathcal{R}_R^{(t)}$ correspond to the partitioned regions of the
$t$-th tree. Given covariates $\mathbf{s}$, the predicted probability
for class $c$ in an RF classifier is:
\begin{equation}
    \hat{p}_c(\mathbf{s}) = \frac{1}{T} \sum_{t=1}^T \frac{\sum_{r = 1}^R \mathbf{1} \{\mathbf{s} \in \mathcal{R}_r^{(t)}\} \cdot \mathbf{1} \{d = \xi_c\}}{\sum_{r = 1}^R \mathbf{1} \{\mathbf{s} \in \mathcal{R}_r^{(t)}\}}.
\end{equation}

In contrast, GBT build classification trees sequentially, where each new tree is trained to optimize a loss function based on the residual errors of the previous trees \citep{friedman2001greedy}. Let $R$ be the number of boosting iterations. 
Given covariates $\mathbf{s}$, for the $r$-th tree where $r \in [R]$, let
$h_r(\mathbf{s})$ be the prediction and $\gamma_{cr}$ be the weight for class $c$.  
The predicted probability for class $c$ in a GBT classifier is obtained using the softmax function of the log-odds:
\begin{equation}
    \hat{p}_c(\mathbf{s}) = \frac{e^{G_c(\mathbf{s})}}{\sum_{c'=1}^{C} e^{G_{c'}(\mathbf{s})}},
\end{equation}
with the log-odds defined as:
\begin{equation}
    G_c(\mathbf{s}) = \sum_{r=1}^R \gamma_{cr} h_r(\mathbf{s}).
\end{equation}

\section{Point Prediction ML Models} \label{sec:regression}
The regression models employed to generate point predictions are presented below. 
\begin{enumerate}
    \item Linear Regression (LinReg): LinReg solves the Ordinary Least Squares (OLS) problem.
    \item Ridge Regression (Ridge): Ridge adds an $L_2$ norm regularization term to the coefficients obtained from OLS.
    \item Least Absolute Shrinkage and Selection Operator Regression (Lasso): Lasso adds an $L_1$ norm regularization term to the coefficients obtained from OLS.
    \item Random Forest (RF-Reg): RF-Reg operates similarly to the Random Forest classifier described in Section \ref{sec: ensemble}, but it builds a collection of regression trees instead.
    \item Gradient Boosted Trees (CatBoost-Reg and XGBoost-Reg): Two versions of Gradient Boosted Trees (GBT) are implemented for regression. These GBT models are similar to the GBT classifiers described in Appendix \ref{sec: ensemble}, except they use regression trees. CatBoost-Reg is implemented using the CatBoost package, which builds symmetric trees and has native support for categorical features \citep{CatBoost}. XGBoost-Reg is implemented using the XGBoost package \citep{Chen:2016:XST:2939672.2939785}.
    \end{enumerate}

\begin{table}[htb]
    \centering
    \caption{Average Out‐of‐Sample MSE of Regression Models.}
    \begin{tabular}{lrrrrrr}
    \toprule
     Model & LinReg &    Ridge &   Lasso &  RF-Reg &  CatBoost-Reg &  XGBoost-Reg \\
    \midrule
     MSE & 2.324 & 2.324 & 3.149 &  1.791 &  1.668 & \textbf{1.666} \\
    \bottomrule
    \end{tabular}
    \label{tab:mse}
\end{table}

\section{Hyperparameters of the ML Models}
Tables~\ref{tab:regression_models}, \ref{tab:classification_models}, and \ref{tab:quantile_models} detail the hyperparameters and corresponding search spaces employed during the model selection and tuning processes for the ML algorithms considered in this study.

\begin{table}[!ht]
\centering
\caption{Hyperparameters for Point Prediction Regression Models}
\label{tab:regression_models}
\begin{tabular}{lll}
\toprule
\textbf{Model} & \textbf{Hyperparameter} & \textbf{Range} \\
\midrule
Ridge & $\alpha$ ($L_2$ regularization) & $1e-4 \leq \alpha \leq 10.0$ \\
\midrule
Lasso & $\alpha$ ($L_1$ regularization) & $1e-4 \leq \alpha \leq 10.0$ \\
\midrule
\multirow{4}{*}{RF-Reg} & n\_estimators & $100 \leq \text{n\_estimators} \leq 1000$ (step: 100) \\
& max\_depth & $3 \leq \text{max\_depth} \leq 10$ \\
& min\_samples\_split & $2 \leq \text{min\_samples\_split} \leq 10$ \\
& min\_samples\_leaf & $1 \leq \text{min\_samples\_leaf} \leq 20$ \\
\midrule
\multirow{8}{*}{CatBoost-Reg} & iterations & $50 \leq \text{iterations} \leq 300$ \\
& depth & $4 \leq \text{depth} \leq 10$ \\
& learning\_rate & $0.01 \leq \text{learning\_rate} \leq 0.3$ \\
& l2\_leaf\_reg & $1 \leq \text{l2\_leaf\_reg} \leq 10$ \\
& bagging\_temperature & $0.0 \leq \text{bagging\_temperature} \leq 1.0$ \\
& random\_strength & $0.0 \leq \text{random\_strength} \leq 10.0$ \\
& border\_count & $1 \leq \text{border\_count} \leq 255$ \\
& colsample\_bylevel & $0.5 \leq \text{colsample\_bylevel} \leq 1.0$ \\
\midrule
\multirow{13}{*}{XGBoost-Reg} & n\_estimators & $100 \leq \text{n\_estimators} \leq 1000$ (step: 100) \\
& max\_depth & $3 \leq \text{max\_depth} \leq 10$ \\
& learning\_rate & $0.01 \leq \text{learning\_rate} \leq 0.3$ \\
& subsample & $0.5 \leq \text{subsample} \leq 1.0$ \\
& colsample\_bytree & $0.5 \leq \text{colsample\_bytree} \leq 1.0$ \\
& gamma & $0 \leq \text{gamma} \leq 5$ \\
& reg\_alpha & $0.0 \leq \text{reg\_alpha} \leq 100.0$ \\
& reg\_lambda & $0.0 \leq \text{reg\_lambda} \leq 100.0$ \\
& min\_child\_weight & $1 \leq \text{min\_child\_weight} \leq 10$ \\
& max\_delta\_step & $0 \leq \text{max\_delta\_step} \leq 10$ \\
& colsample\_bylevel & $0.5 \leq \text{colsample\_bylevel} \leq 1.0$ \\
& colsample\_bynode & $0.5 \leq \text{colsample\_bynode} \leq 1.0$ \\
\bottomrule
\end{tabular}
\end{table}

\begin{table}[!ht]
\centering
\caption{Hyperparameters for Classification Models}
\label{tab:classification_models}
\begin{tabular}{lll}
\toprule
\textbf{Model} & \textbf{Hyperparameter} & \textbf{Range} \\
\midrule
MLR & C (regularization) & $1e-4 \leq \text{C} \leq 1e2$ \\
\midrule
\multirow{2}{*}{Classification Tree} & max\_depth & $3 \leq \text{max\_depth} \leq 10$ \\
& min\_samples\_leaf & $10 \leq \text{min\_samples\_leaf} \leq 30$ (step: 5) \\
\midrule
\multirow{4}{*}{RF-MC-CLF} & n\_estimators & $600 \leq \text{n\_estimators} \leq 1000$ (step: 200) \\
& max\_depth & $3 \leq \text{max\_depth} \leq 10$ \\
& min\_samples\_split & $2 \leq \text{min\_samples\_split} \leq 10$ \\
& min\_samples\_leaf & $10 \leq \text{min\_samples\_leaf} \leq 30$ (step: 5) \\
\midrule
\multirow{8}{*}{CatBoost-MC-CLF} & iterations & $50 \leq \text{iterations} \leq 300$ \\
& depth & $4 \leq \text{depth} \leq 10$ \\
& learning\_rate & $0.01 \leq \text{learning\_rate} \leq 0.3$ \\
& l2\_leaf\_reg & $1 \leq \text{l2\_leaf\_reg} \leq 10$ \\
& bagging\_temperature & $0.0 \leq \text{bagging\_temperature} \leq 1.0$ \\
& random\_strength & $0.0 \leq \text{random\_strength} \leq 10.0$ \\
& border\_count & $1 \leq \text{border\_count} \leq 255$ \\
& colsample\_bylevel & $0.5 \leq \text{colsample\_bylevel} \leq 1.0$ \\
\bottomrule
\end{tabular}
\end{table}

\begin{table}[!ht]
\centering
\caption{Hyperparameters for Quantile Regression Models}
\label{tab:quantile_models}
\begin{tabular}{lll}
\toprule
\textbf{Model} & \textbf{Hyperparameter} & \textbf{Range} \\
\midrule
Regression Tree & max\_depth & $3 \leq \text{max\_depth} \leq 10$ \\
& min\_samples\_leaf & $10 \leq \text{min\_samples\_leaf} \leq 30$ (step: 10) \\
\midrule
\multirow{4}{*}{QRF} & n\_estimators & $600 \leq \text{n\_estimators} \leq 1000$ (step: 200) \\
& max\_depth & $3 \leq \text{max\_depth} \leq 10$ \\
& min\_samples\_split & $2 \leq \text{min\_samples\_split} \leq 10$ \\
& min\_samples\_leaf & $10 \leq \text{min\_samples\_leaf} \leq 30$ (step: 5) \\
\bottomrule
\end{tabular}
\end{table}

\revision{
\section{Evaluation of the Simulator} \label{sec:simulator_eval}

The simulator is compared against a baseline that uses multinomial logistic regression (MLR), a common approach in the literature \citep{salari2022real}.

Table \ref{tab:simulator_crps} reports the average CRPS across all carriers. The proposed calibrated MC-CLF simulator achieves a CRPS of 0.314, compared with 0.779 for the MLR baseline, indicating a substantial improvement in distributional accuracy at the individual‐instance level.

Figure \ref{fig:overall_simulator_dist} overlays each simulator’s joint predicted distribution on the true empirical distribution from the test set. The proposed calibrated MC-CLF simulator closely matches the empirical mass at zero deviation and reproduces both early‐ and late‐delivery tails, whereas the MLR baseline systematically overpredicts late deviations and underpredicts on‐time deliveries.

Figure \ref{fig:example_simulator_dist} presents three carrier‐specific examples. In each case, the proposed calibrated MC-CLF simulator aligns tightly with the carrier’s empirical distribution across all deviation classes. By contrast, the MLR baseline exhibits clear misalignments---particularly in the majority classes and tails---underscoring the proposed simulator’s superior ability to capture carrier‐level heterogeneity.

\begin{table}[!ht]
\caption{Average In-Sample CRPS Across All Carriers.}
\centering
\begin{tabular}{lcc}
\toprule
Model &  Calibrated MC-CLF Simulator & MLR Simulator (Baseline)  \\
\midrule
CRPS & 0.314 & 0.779\\
\bottomrule
\end{tabular}
\label{tab:simulator_crps}
\end{table}
}

\begin{figure}[!ht]
    \centering
    \begin{subfigure}[b]{0.48\textwidth}
        \centering
        \includegraphics[width=\textwidth]{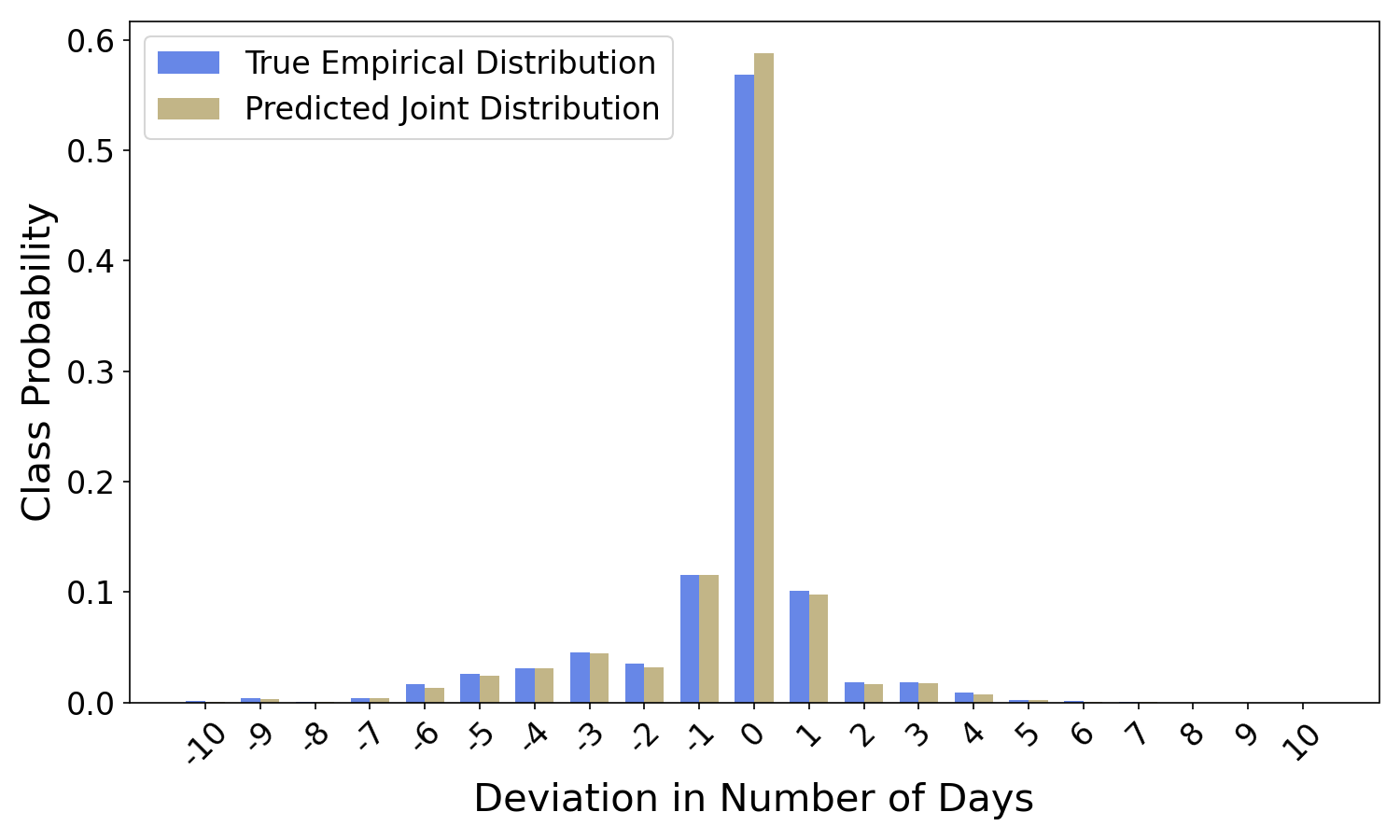}
        \caption{Calibrated MC-CLF Simulator}      \label{fig:overall_dist_comparison_simulator}
    \end{subfigure}    
        \hfill
        \begin{subfigure}[b]{0.48\textwidth}
        \centering
        \includegraphics[width=\textwidth]{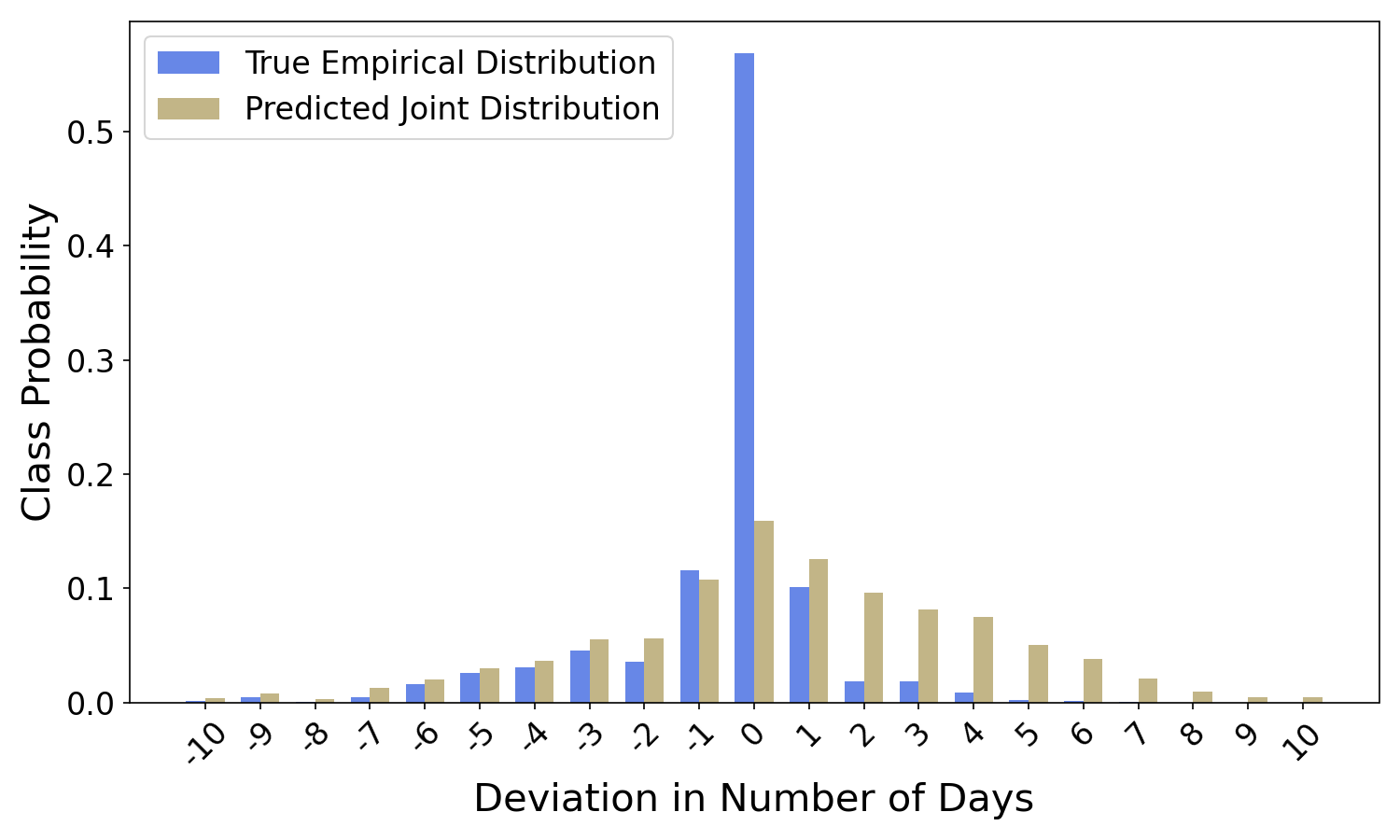}
        \caption{Baseline}
    \label{fig:overall_dist_comparison_baseline}
    \end{subfigure}   
        \caption{Overall Joint Predicted Distribution of the Simulators vs. the True Empirical Distribution}
\label{fig:overall_simulator_dist}
\end{figure}

\begin{figure}[!ht]
    \centering
    \begin{subfigure}[b]{0.48\textwidth}
        \centering
        \includegraphics[width=\textwidth]{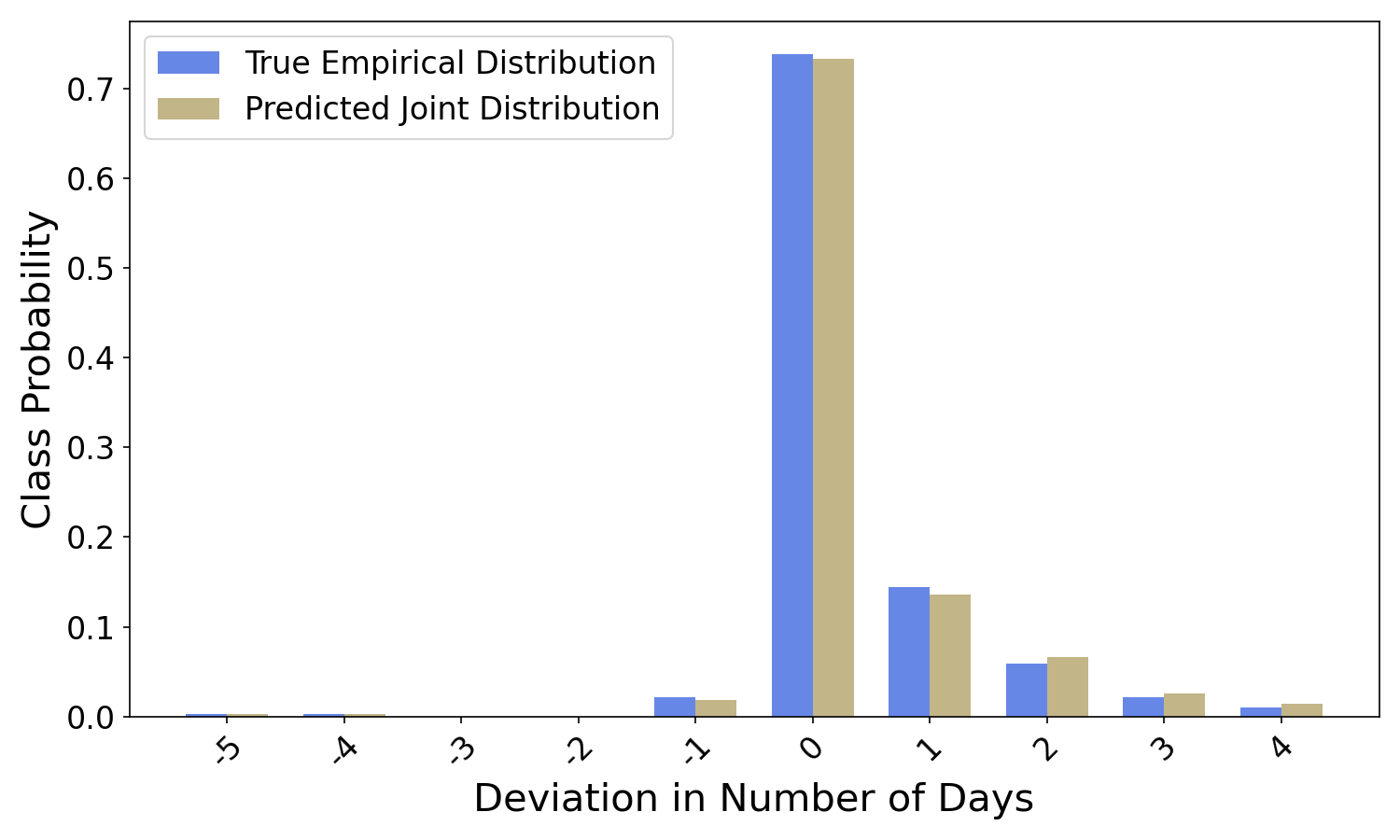}
        \caption{Example Carrier 1: Calibrated MC-CLF Simulator}      \label{fig:16_dist_comparison_simulator}
    \end{subfigure}    
        \hfill
        \begin{subfigure}[b]{0.48\textwidth}
        \centering
        \includegraphics[width=\textwidth]{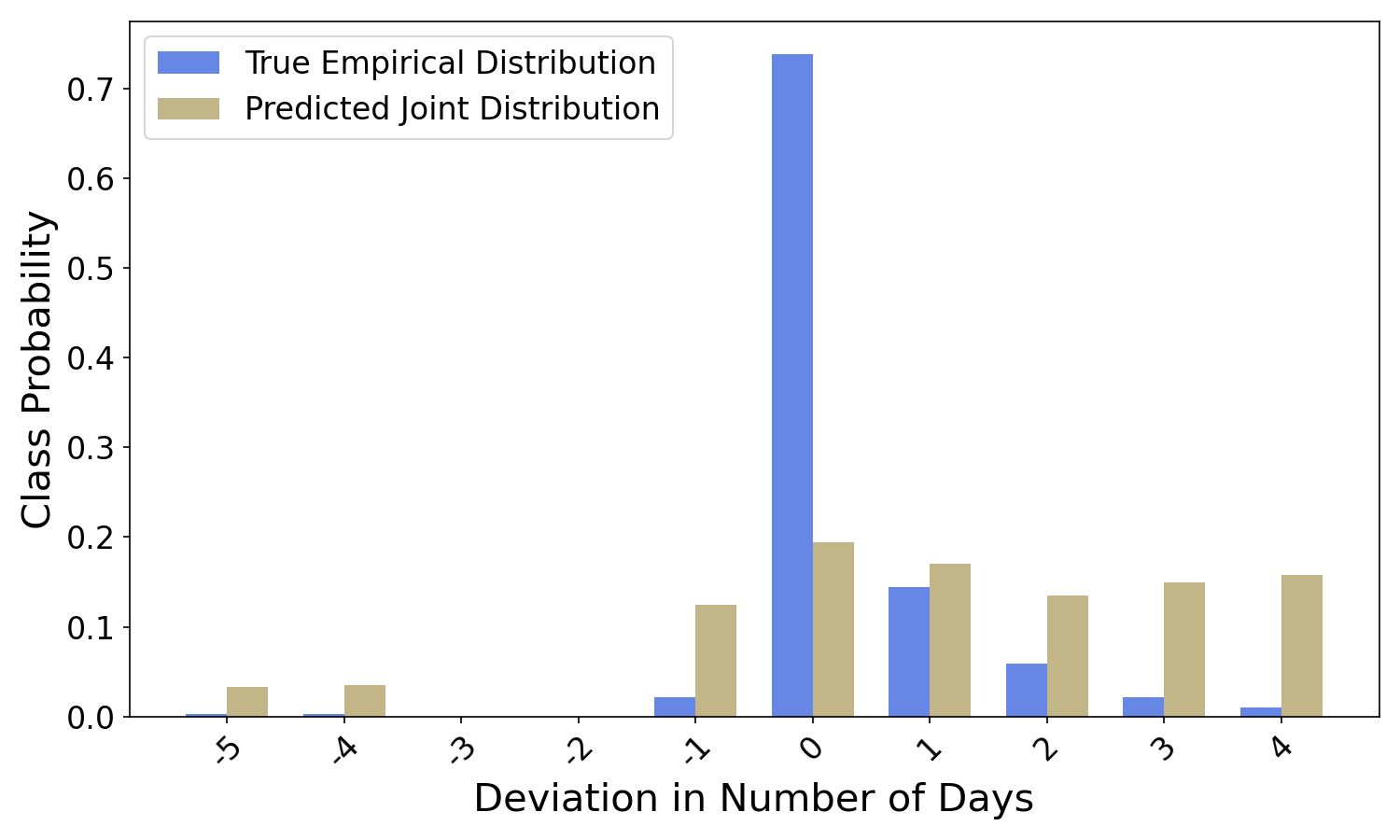}
        \caption{Example Carrier 1: Baseline}
    \label{fig:16_dist_comparison_baseline}
    \end{subfigure}  
        \begin{subfigure}[b]{0.48\textwidth}
        \centering
        \includegraphics[width=\textwidth]{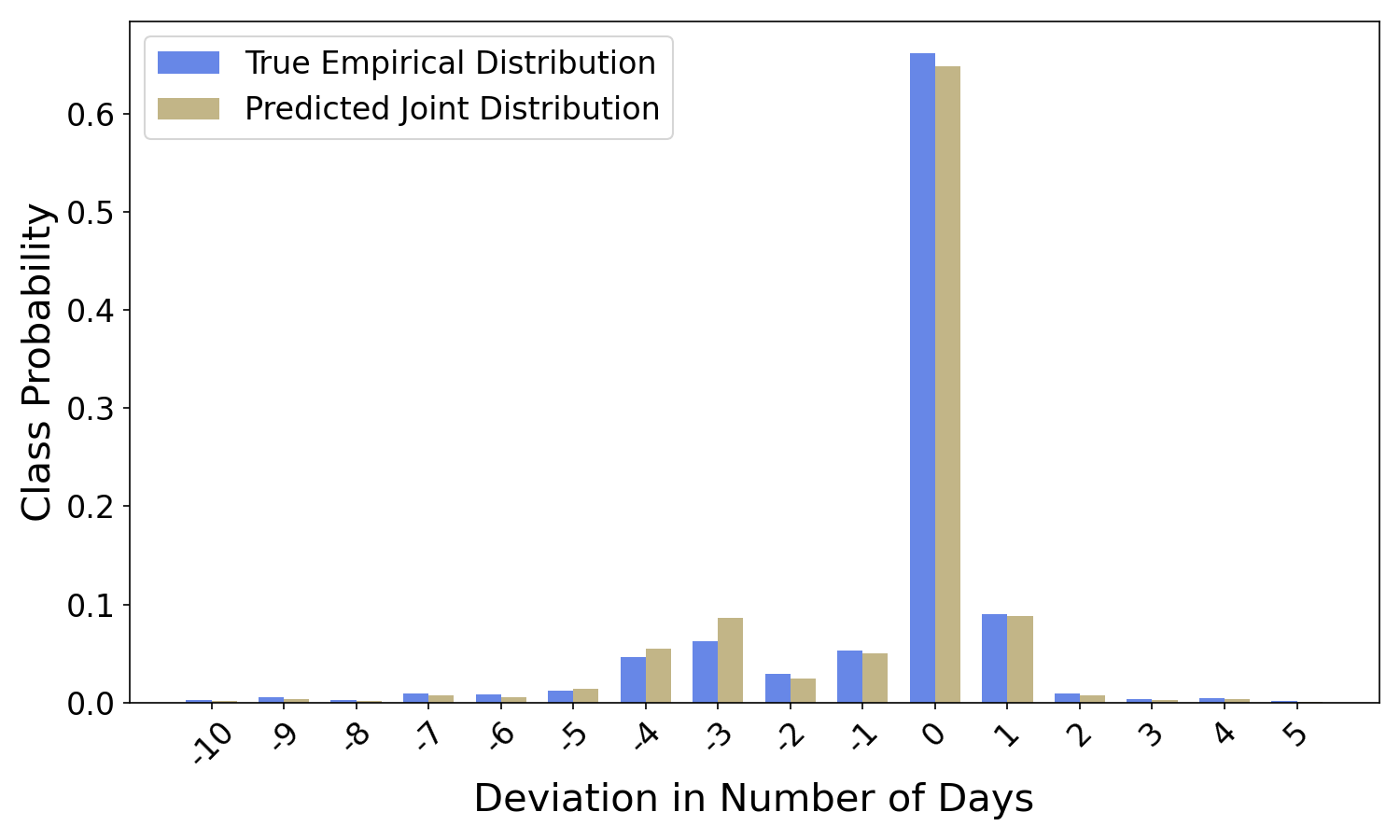}
        \caption{Example Carrier 2: Calibrated MC-CLF Simulator}      \label{fig:28_dist_comparison_simulator}
    \end{subfigure}    
        \hfill
    \begin{subfigure}[b]{0.48\textwidth}
        \centering
        \includegraphics[width=\textwidth]{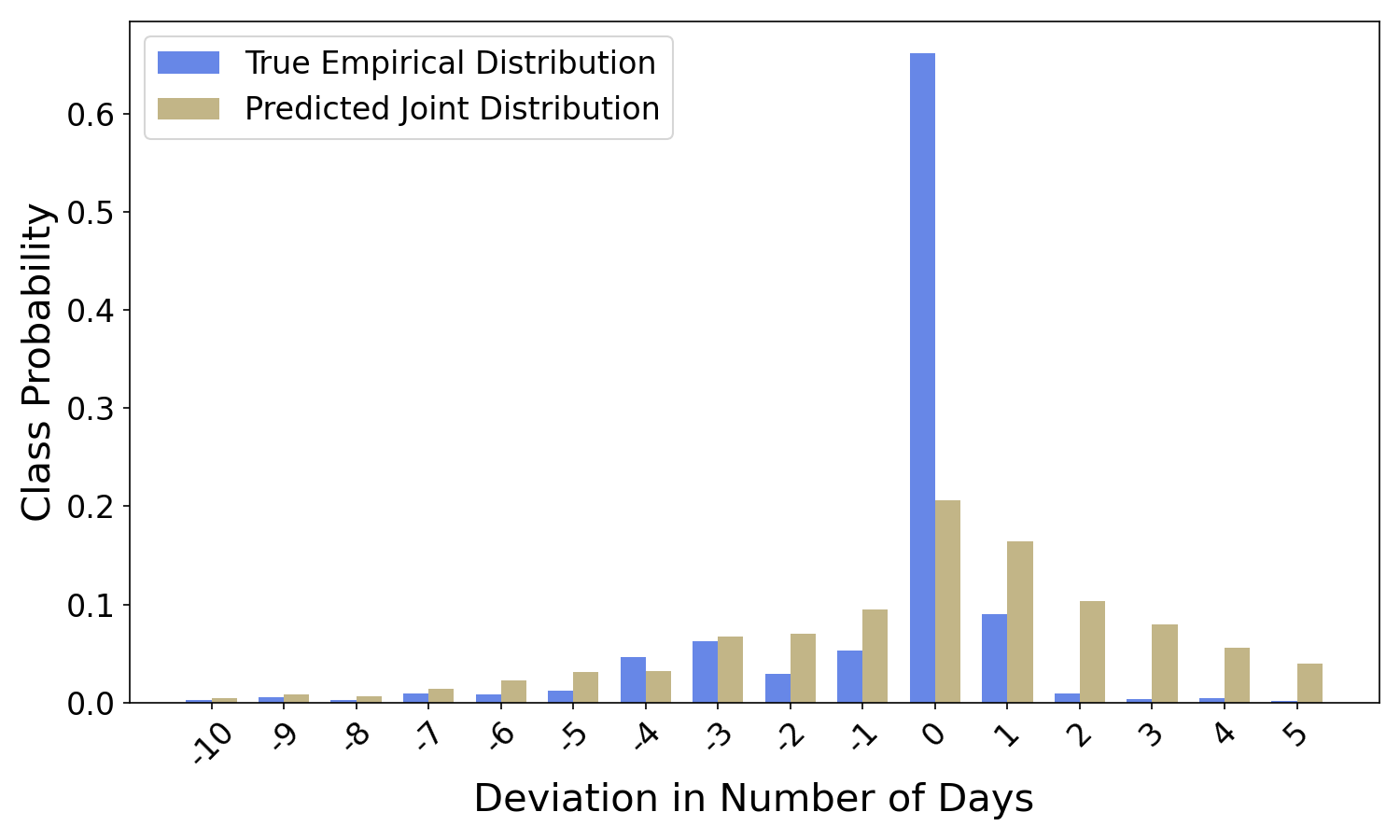}
        \caption{Example Carrier 2: Baseline}
    \label{fig:28_dist_comparison_baseline}
    \end{subfigure}   
    \begin{subfigure}[b]{0.48\textwidth}
        \centering
        \includegraphics[width=\textwidth]{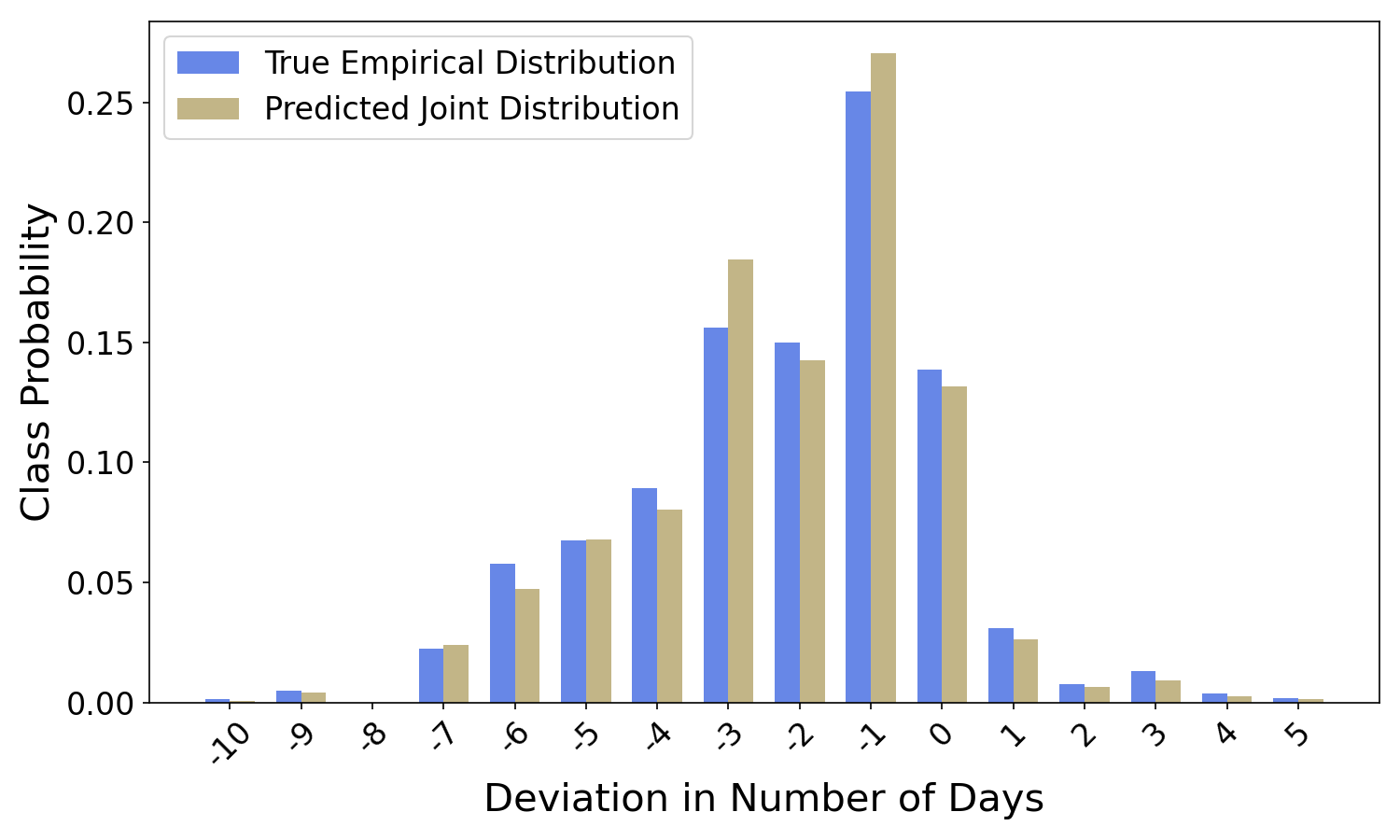}
        \caption{Example Carrier 3: Calibrated MC-CLF Simulator}      \label{fig:44_dist_comparison_simulator}
    \end{subfigure}    
        \hfill
        \begin{subfigure}[b]{0.48\textwidth}
        \centering
        \includegraphics[width=\textwidth]{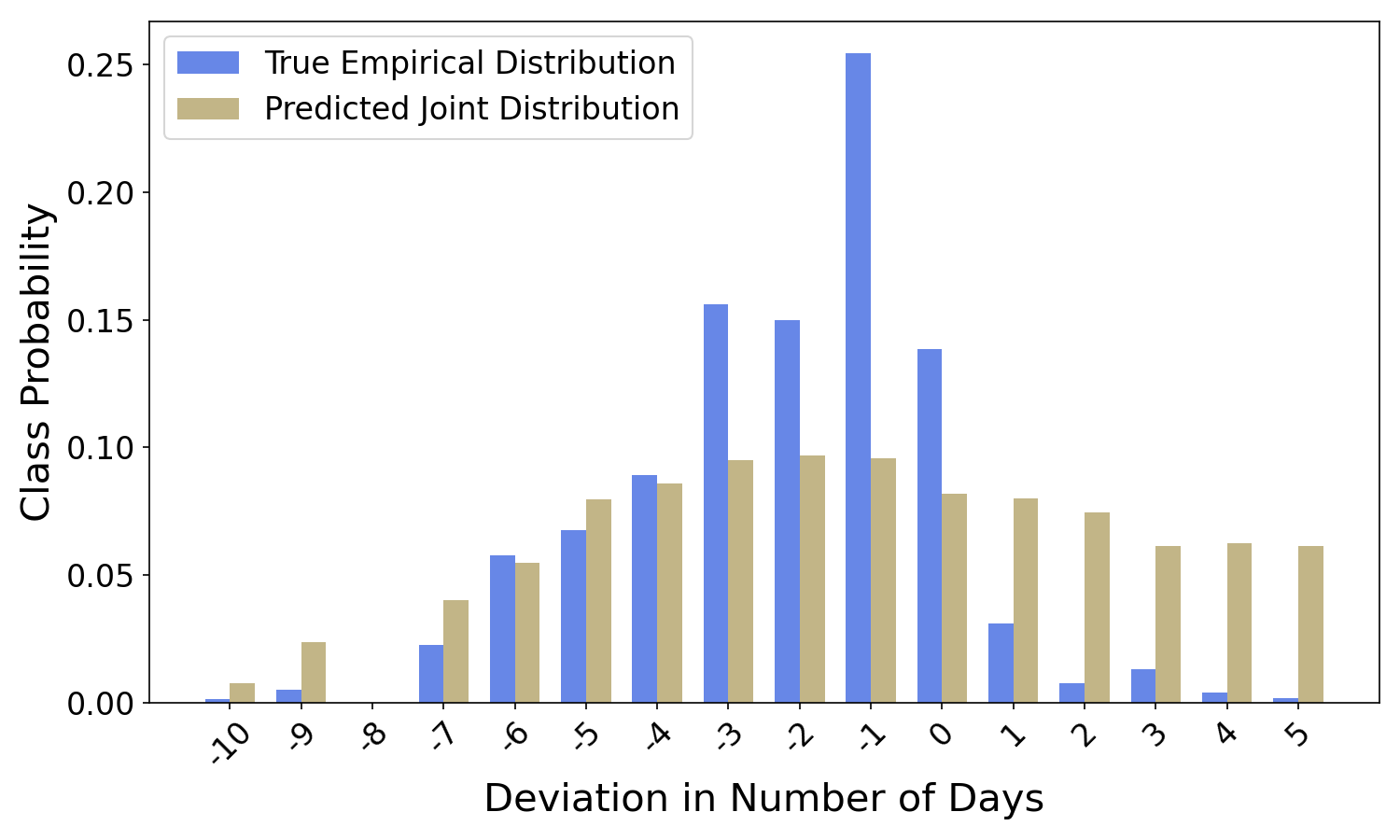}
        \caption{Example Carrier 3: Baseline}
    \label{fig:44_dist_comparison_baseline}
    \end{subfigure}   
        \caption{Selected Carrier Examples: Joint Predicted Distribution of the Simulators vs. the True Empirical Distribution}
\label{fig:example_simulator_dist}
\end{figure}

\revision{

\section{Additional Details on Computational Scalability}

Table \ref{tab:objective_comparison} provides a comprehensive comparison of objective values across different instance categories for the six proposed CSO methods. \revisiontwo{Table \ref{tab:time_statistics} summarizes the computation times for the baseline models and the proposed CSO-CLF methods.}

\begin{table}[ht]
\centering
\caption{Realized Average (\(\pm\) 95\% CI margin of error) and Worst‐Case Objective Values Across Instance Categories.}
\begin{adjustbox}{max width=\textwidth}
\begin{tabular}{lcc cc cc cc cc cc}
\toprule
Instance 
& \multicolumn{2}{c}{C-SAA-CLF} 
& \multicolumn{2}{c}{C-SAA-QRF} 
& \multicolumn{2}{c}{C-RO-B-CLF} 
& \multicolumn{2}{c}{C-RO-B-QRF} 
& \multicolumn{2}{c}{C-RO-D-CLF} 
& \multicolumn{2}{c}{C-RO-D-QRF} \\
\cmidrule(lr){2-3} \cmidrule(lr){4-5} \cmidrule(lr){6-7}
\cmidrule(lr){8-9} \cmidrule(lr){10-11} \cmidrule(lr){12-13}
& Avg & Worst 
& Avg & Worst 
& Avg & Worst 
& Avg & Worst 
& Avg & Worst 
& Avg & Worst 
\\
\midrule
Small  
& $0.821\pm0.003$ & 0.790 
& $0.848\pm0.004$ & 0.813 
& $0.851\pm0.004$ & 0.794 
& $0.871\pm0.003$ & 0.813 
& $0.936\pm0.004$ & 0.757 
& $0.920\pm0.004$ & 0.766 \\

Medium 
& $0.799\pm0.008$ & 0.932 
& $0.806\pm0.009$ & 0.925 
& $0.797\pm0.009$ & 0.834 
& $0.815\pm0.009$ & 0.821 
& $0.890\pm0.010$ & 0.822 
& $0.883\pm0.009$ & 0.825 \\

Large  
& $0.847\pm0.012$ & 1.032 
& $0.871\pm0.017$ & 1.055 
& $0.823\pm0.014$ & 0.982 
& $0.836\pm0.016$ & 0.981 
& $0.910\pm0.017$ & 0.972 
& $0.904\pm0.014$ & 0.973 \\

\midrule
Total  
& $0.820\pm0.003$ & 0.823 
& $0.844\pm0.003$ & 0.842 
& $0.842\pm0.003$ & 0.811 
& $0.861\pm0.003$ & 0.824 
& $0.928\pm0.003$ & 0.779 
& $0.914\pm0.003$ & 0.787 \\
\bottomrule
\end{tabular}
\end{adjustbox}
\label{tab:objective_comparison}
\end{table}

\begin{table}[ht]
\centering
\caption{\revisiontwo{Computation Time Statistics (in seconds) Across Instance Categories.}}
\begin{adjustbox}{max width=\textwidth}
\begin{tabular}{lccc ccc ccc ccc ccc ccc ccc}
\toprule
Instance 
& \multicolumn{3}{c}{C-SAA-CLF} 
& \multicolumn{3}{c}{C-RO-B-CLF} 
& \multicolumn{3}{c}{C-RO-D-CLF} 
& \multicolumn{3}{c}{C-Empirical-SAA} 
& \multicolumn{3}{c}{Empirical-SAA} 
& \multicolumn{3}{c}{PTO} 
& \multicolumn{3}{c}{Greedy}
\\
\cmidrule(lr){2-4} \cmidrule(lr){5-7} \cmidrule(lr){8-10} \cmidrule(lr){11-13} \cmidrule(lr){14-16} \cmidrule(lr){17-19} \cmidrule(lr){20-22}
& Avg & Std & Max 
& Avg & Std & Max 
& Avg & Std & Max 
& Avg & Std & Max 
& Avg & Std & Max 
& Avg & Std & Max 
& Avg & Std & Max 
\\
\midrule
Small  &  1.620 &  2.267 &  28.021 
       &  0.016 &  0.024 &   0.374 
       &  0.199 &  0.301 &   5.475 
       &  2.520 & 3.552  & 46.90
       &  2.510 & 3.547  & 46.86
       & 0.014 & 0.021 & 0.684
       & 0.001 & 0.001 & 0.008
       \\
Medium & 27.831 & 10.245 & 107.347 
       &  0.281 &  0.106 &   0.860 
       &  3.229 &  2.008 &  39.222 
       & 43.742 & 16.475 & 165.717
       & 43.864 & 16.678 & 179.577
       & 0.243 & 0.103 & 1.027
       &  0.013 & 0.007 & 0.060
       \\
Large  & 59.266 & 22.590 & 182.401 
       &  0.698 &  0.322 &   2.856 
       &  7.683 &  7.056 & 111.188 
       & 91.190 & 35.100 & 278.677
        & 92.225 & 35.421 & 295.212
        & 0.629 & 0.329 & 2.584
        & 0.052 & 0.034 & 0.317
       \\
\midrule
Total  &  8.974 & 17.433 & 182.401 
       &  0.096 &  0.203 &   2.856 
       &  1.101 &  2.791 & 111.188 
       & 13.959 & 27.052 & 278.677
        & 14.033 & 27.304 & 295.212
        & 0.086 & 0.187 & 2.584
        & 0.006 & 0.016 & 0.317
       \\
\bottomrule
\end{tabular}
\end{adjustbox}
\label{tab:time_statistics}
\end{table}

\subsection{Sample Size Selection: Balancing Solution Quality and Computational Efficiency}

Figures \ref{fig:saa_size_time} and \ref{fig:saa_size_obj} illustrate the trade-off between average realized objective value and computation time for the two C-SAA methods.

Similarly, Figures \ref{fig:cro_d_size_time} and \ref{fig:cro_d_size_obj} depict the trade-off between worst-case realized objective value and computation time for the two C-RO-D methods.

\begin{figure}[!ht]
\centering
    \begin{subfigure}[b]{0.32\textwidth}
        \centering
        \includegraphics[width=\textwidth]{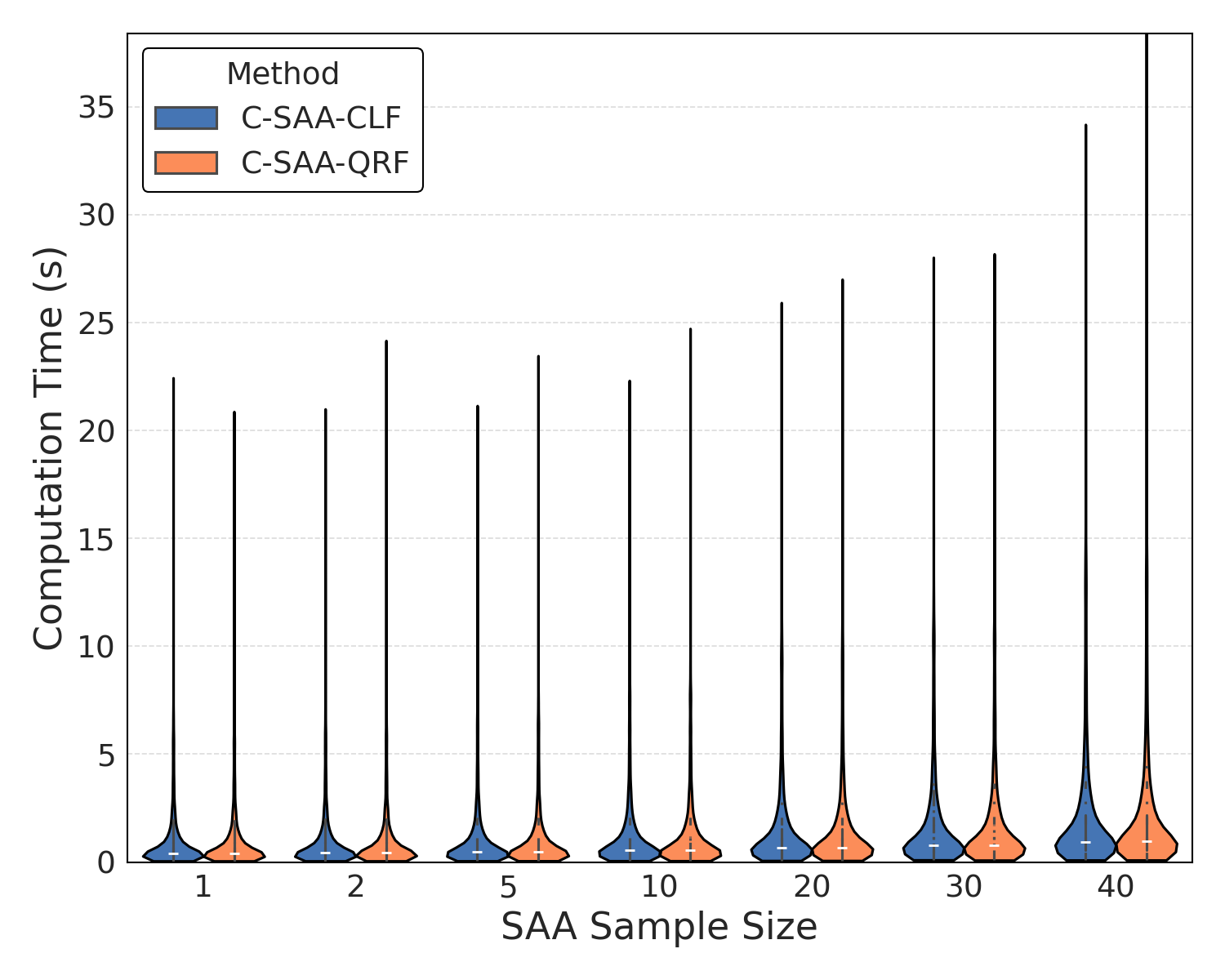} 
        \caption{Small Instances.}
        \label{fig:scens_violin_solve_time_small}
    \end{subfigure}
    \begin{subfigure}[b]{0.32\textwidth}
        \centering
        \includegraphics[width=\textwidth]{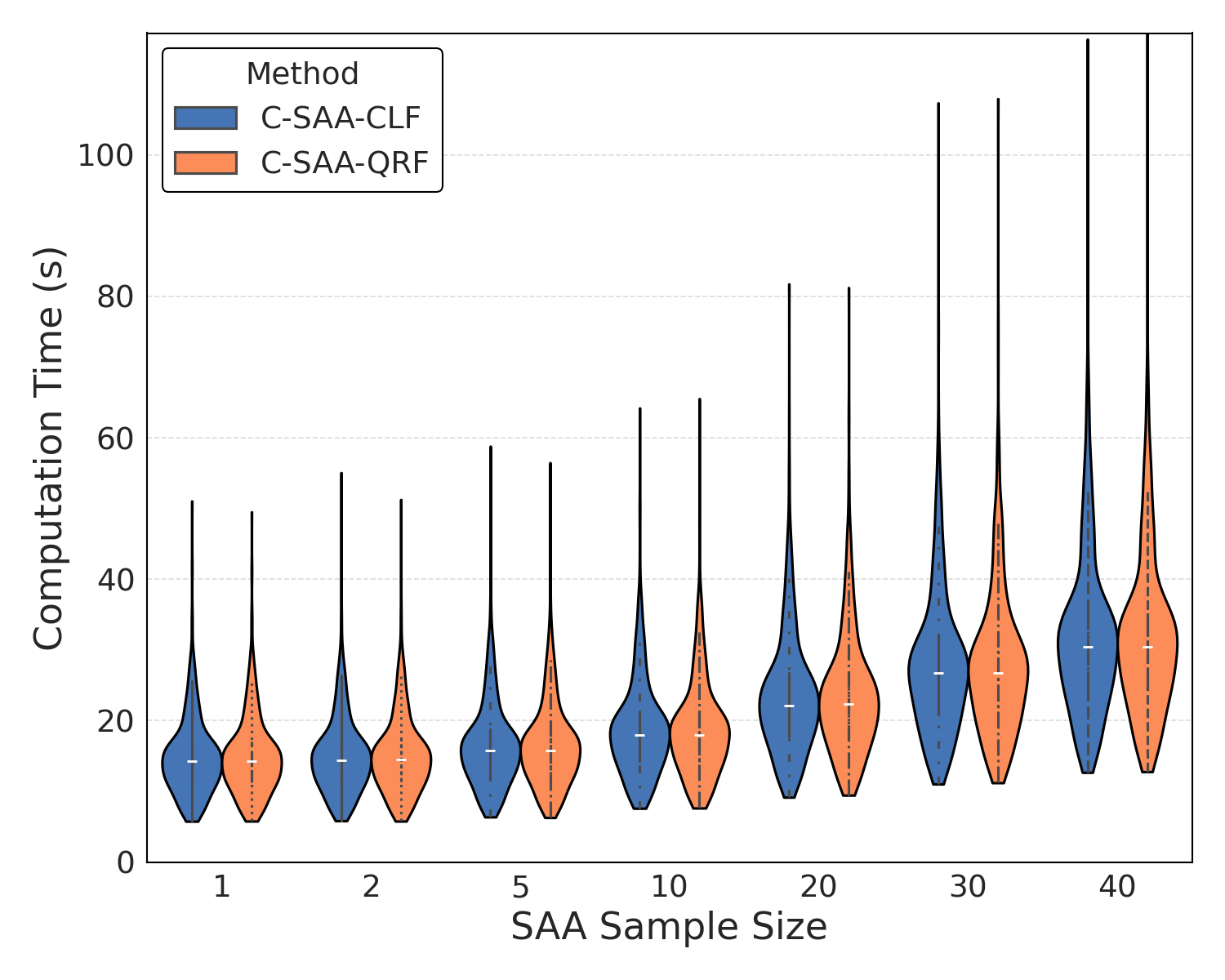} 
        \caption{Medium Instances.}
        \label{fig:scens_violin_solve_time_medium}
    \end{subfigure}
    \begin{subfigure}[b]{0.32\textwidth}
        \centering
        \includegraphics[width=\textwidth]{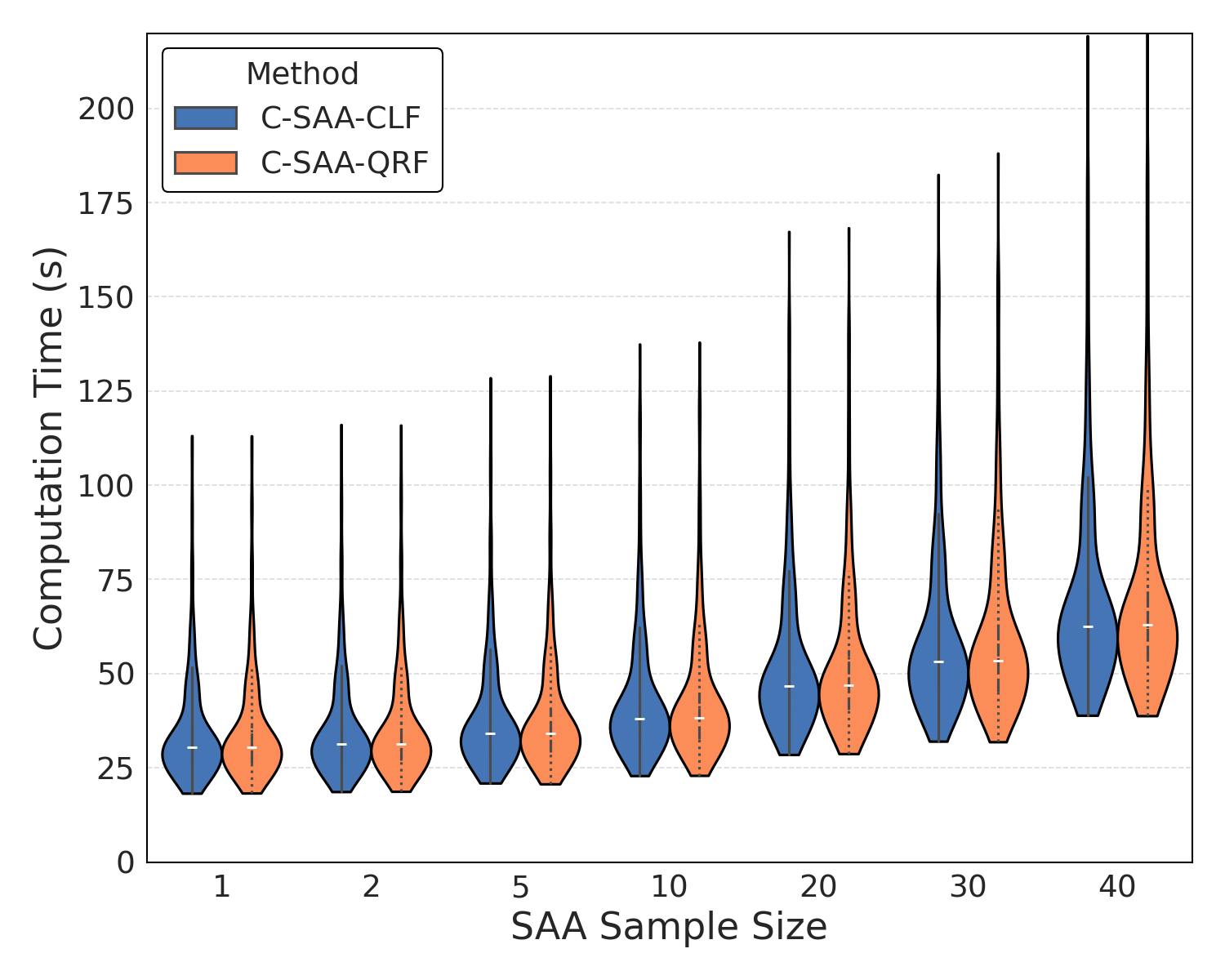} 
        \caption{Large Instances.}
        \label{fig:scens_violin_solve_time_large}
    \end{subfigure}
     \caption{Computation Time by C-SAA Sample Size.}
     \label{fig:saa_size_time}
\end{figure}

\begin{figure}[!ht]
\centering
    \begin{subfigure}[b]{0.32\textwidth}
        \centering
        \includegraphics[width=\textwidth]{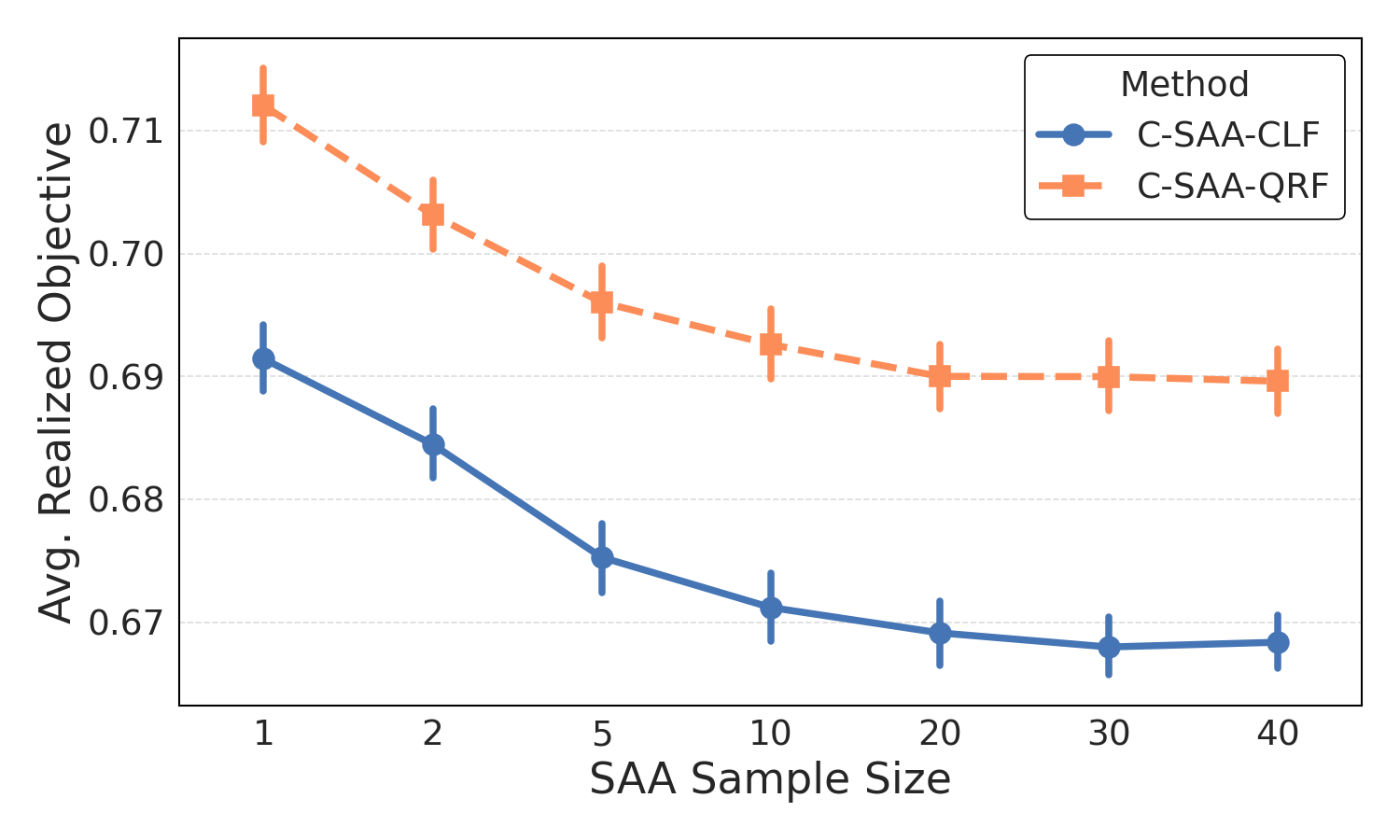} 
        \caption{Small Instances.}
    \end{subfigure}
    \begin{subfigure}[b]{0.32\textwidth}
        \centering
        \includegraphics[width=\textwidth]{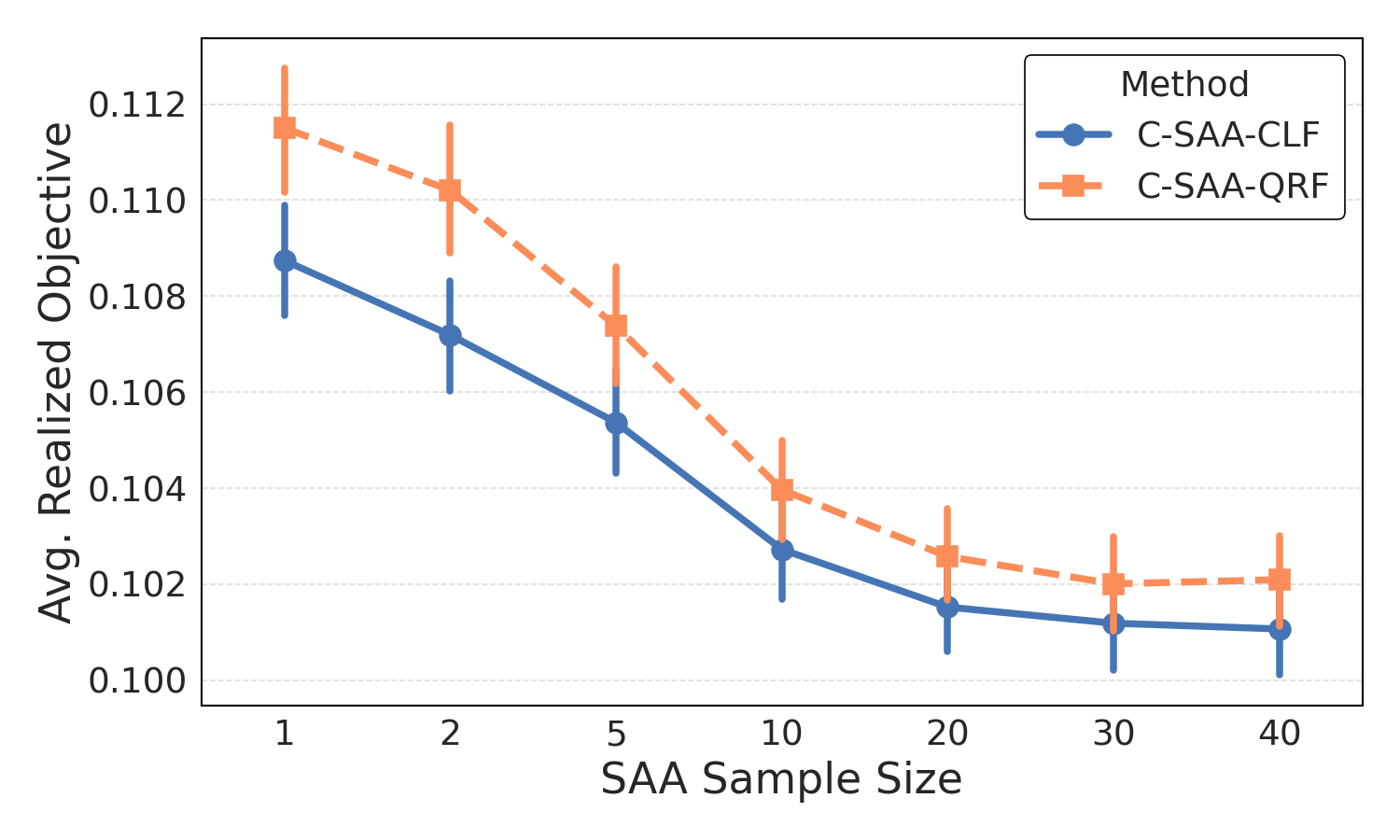} 
        \caption{Medium Instances.}
    \end{subfigure}
    \begin{subfigure}[b]{0.32\textwidth}
        \centering
        \includegraphics[width=\textwidth]{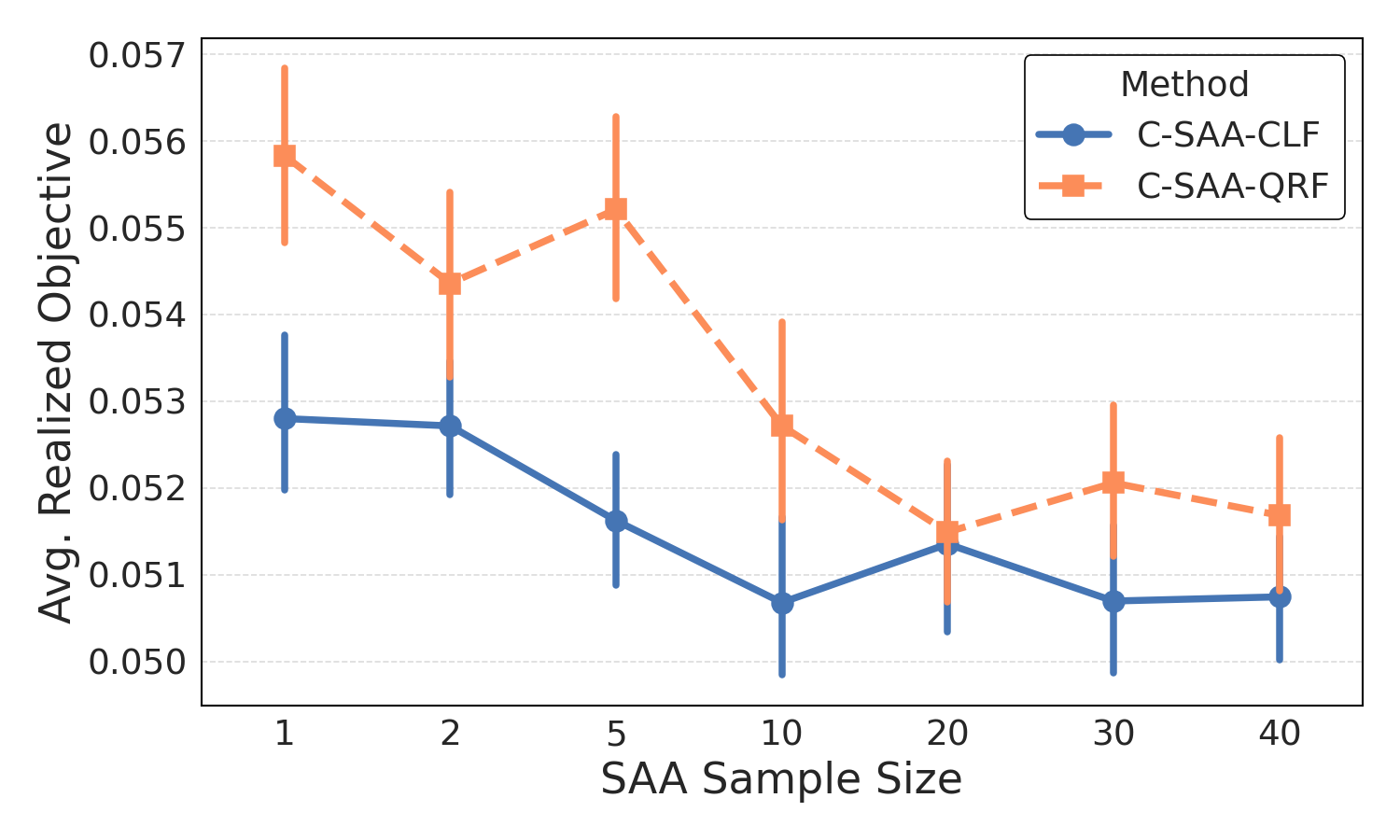} 
        \caption{Large Instances.}
    \end{subfigure}
     \caption{Average Realized Objective Value by C-SAA Sample Size.}
     \label{fig:saa_size_obj}
\end{figure}

\begin{figure}[!ht]
\centering
    \begin{subfigure}[b]{0.32\textwidth}
        \centering
        \includegraphics[width=\textwidth]{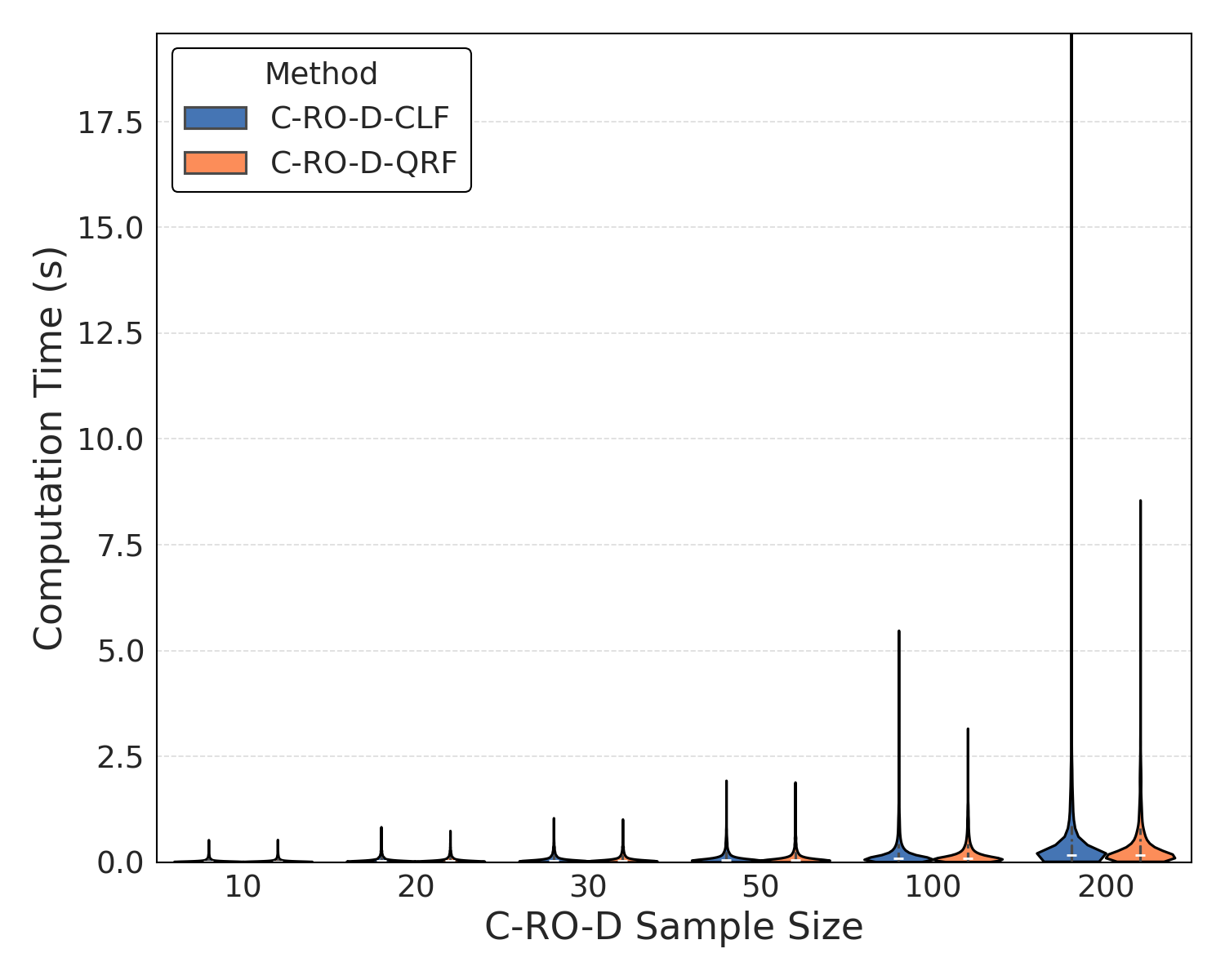}
        \caption{Small Instances.}
    \label{fig:ro_num_scenarios_violin_solve_time_small}
    \end{subfigure}
    \begin{subfigure}[b]{0.32\textwidth}
        \centering
        \includegraphics[width=\textwidth]{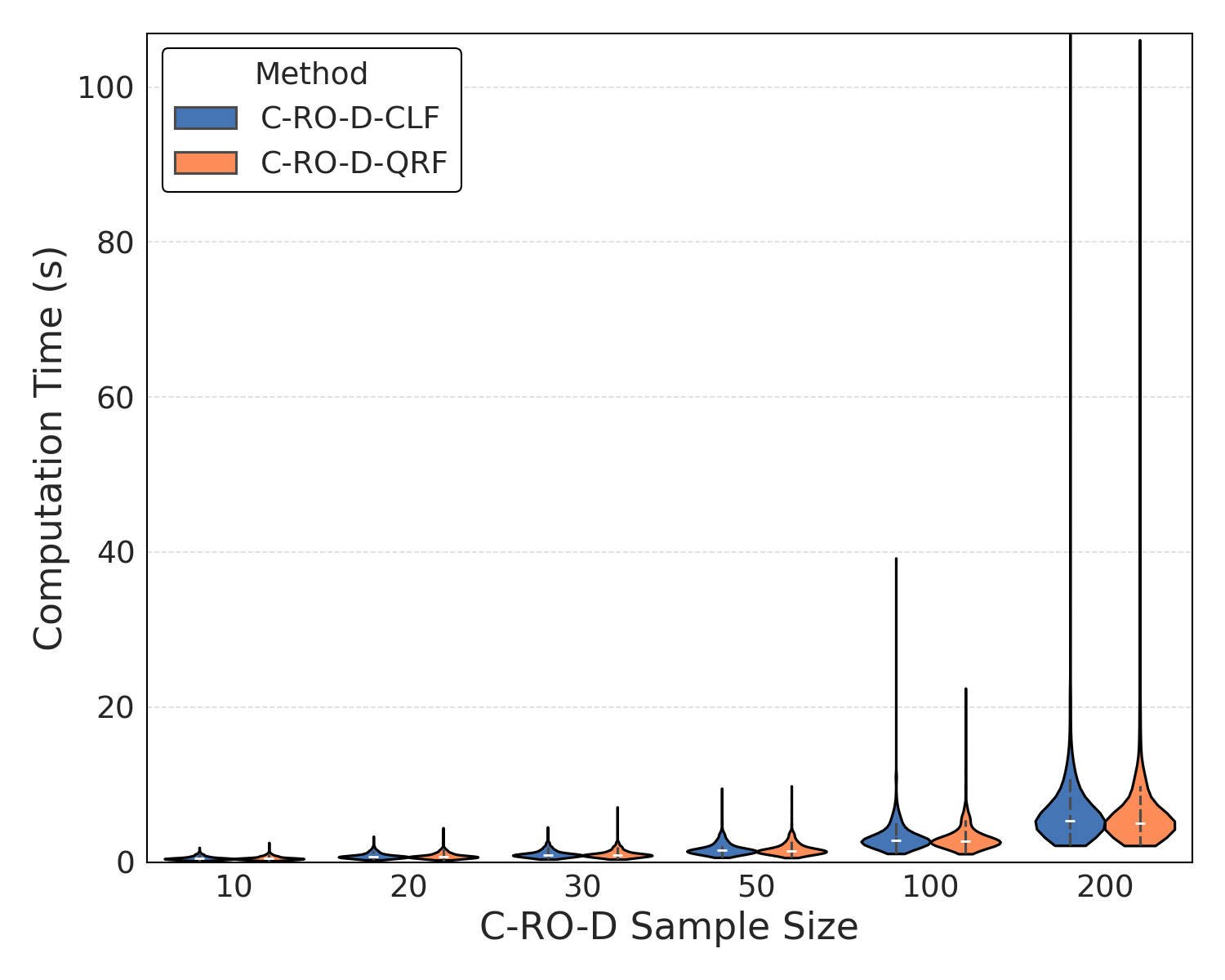} 
        \caption{Medium Instances.}
        \label{fig:ro_num_scenarios_violin_solve_time_medium}
    \end{subfigure}
    \begin{subfigure}[b]{0.32\textwidth}
        \centering
        \includegraphics[width=\textwidth]{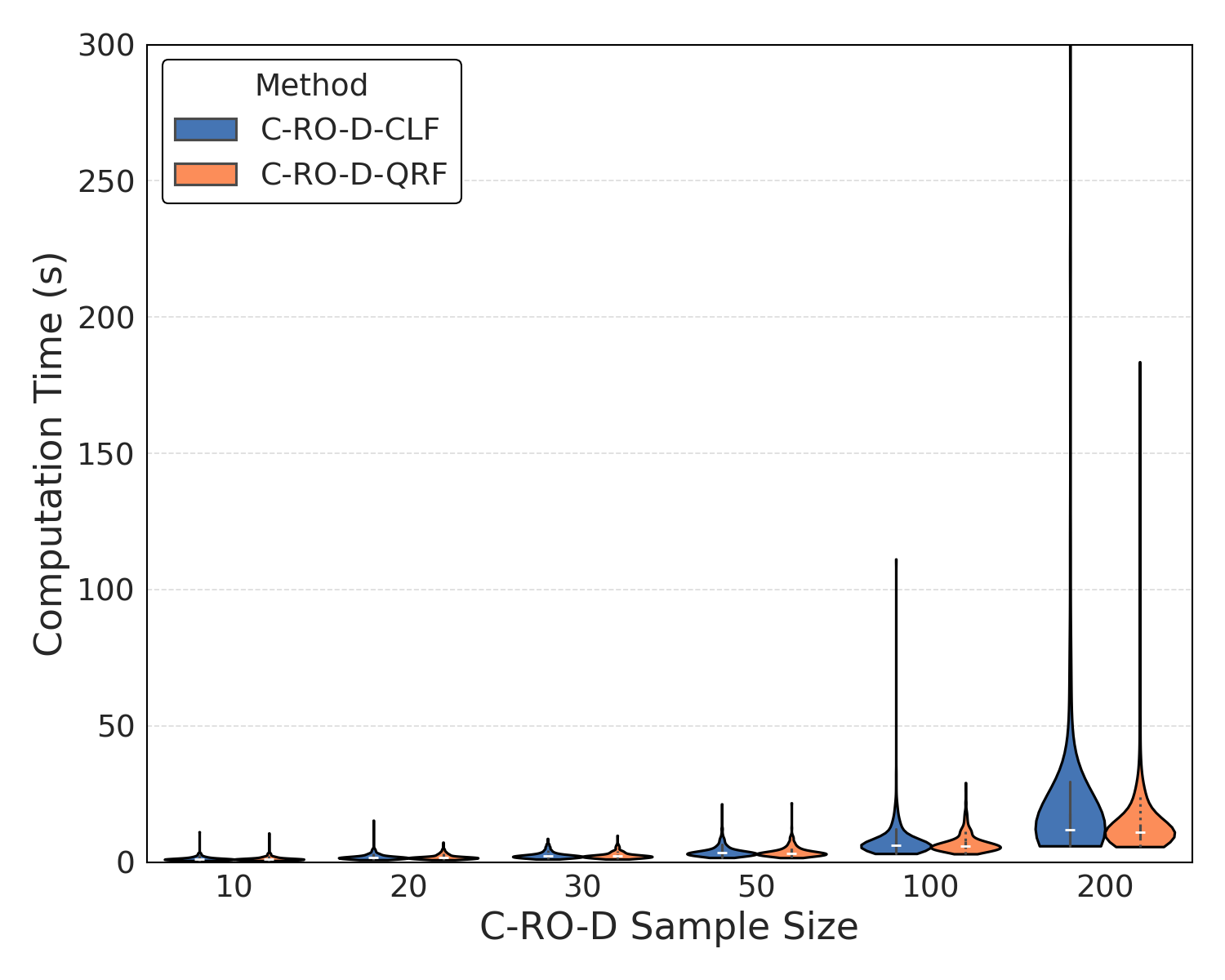} 
        \caption{Large Instances.}
        \label{fig:ro_num_scenarios_violin_solve_time_large}
    \end{subfigure}
     \caption{Computation Time by C-RO-D Sample Size.}
     \label{fig:cro_d_size_time}
\end{figure}

\begin{figure}[!ht]
\centering
    \begin{subfigure}[b]{0.32\textwidth}
        \centering
        \includegraphics[width=\textwidth]{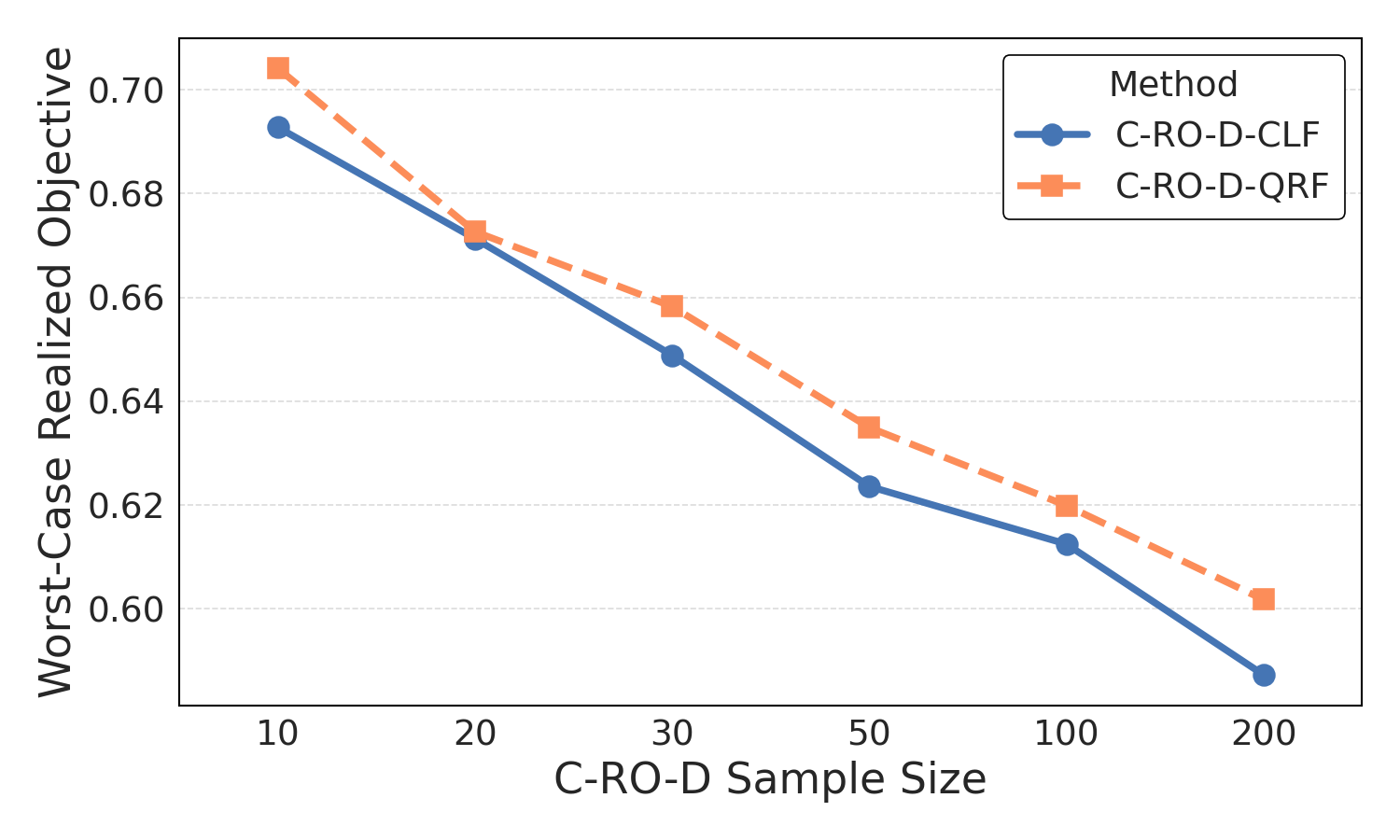} 
        \caption{Small Instances.}
    \end{subfigure}
    \begin{subfigure}[b]{0.32\textwidth}
        \centering
        \includegraphics[width=\textwidth]{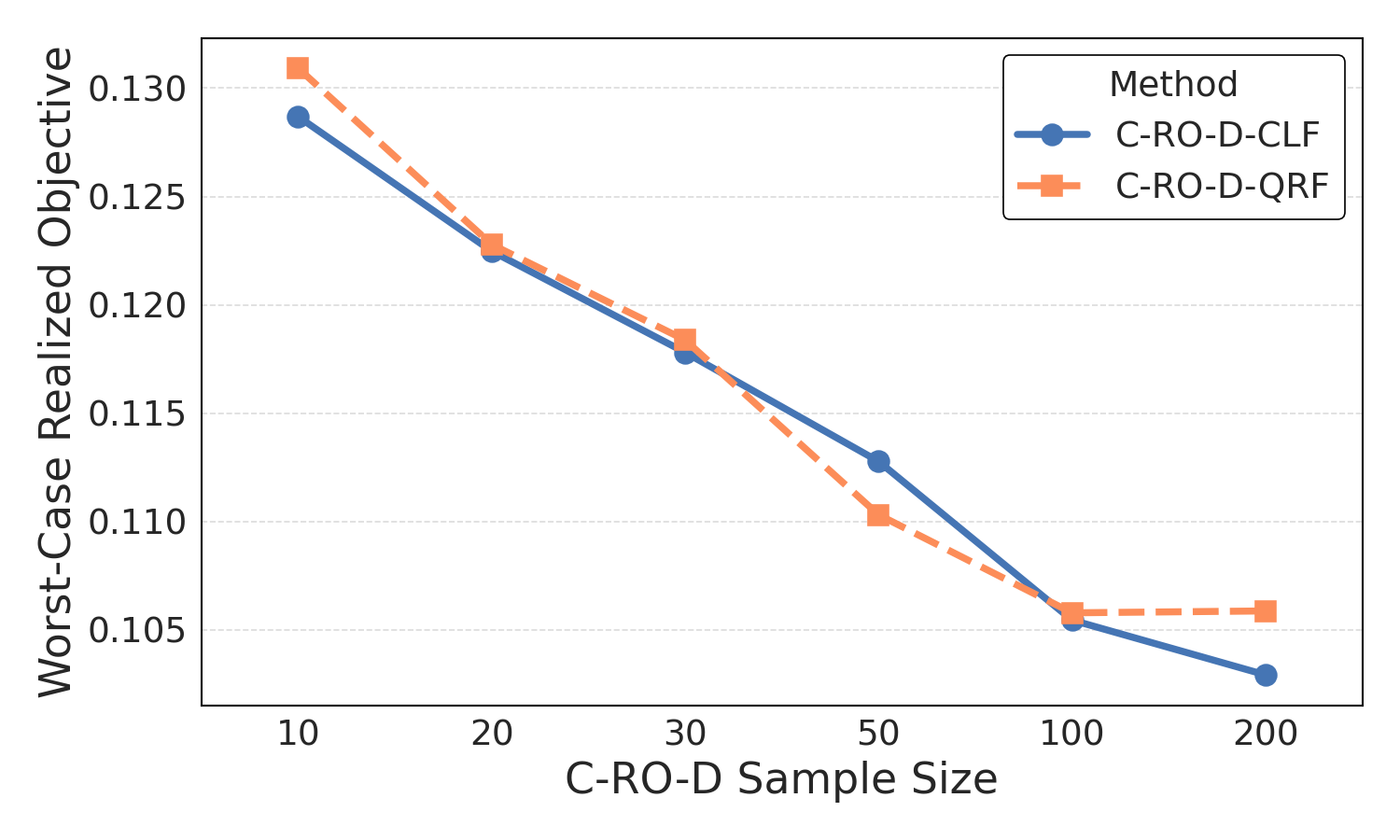} 
        \caption{Medium Instances.}
    \end{subfigure}
    \begin{subfigure}[b]{0.32\textwidth}
        \centering
        \includegraphics[width=\textwidth]{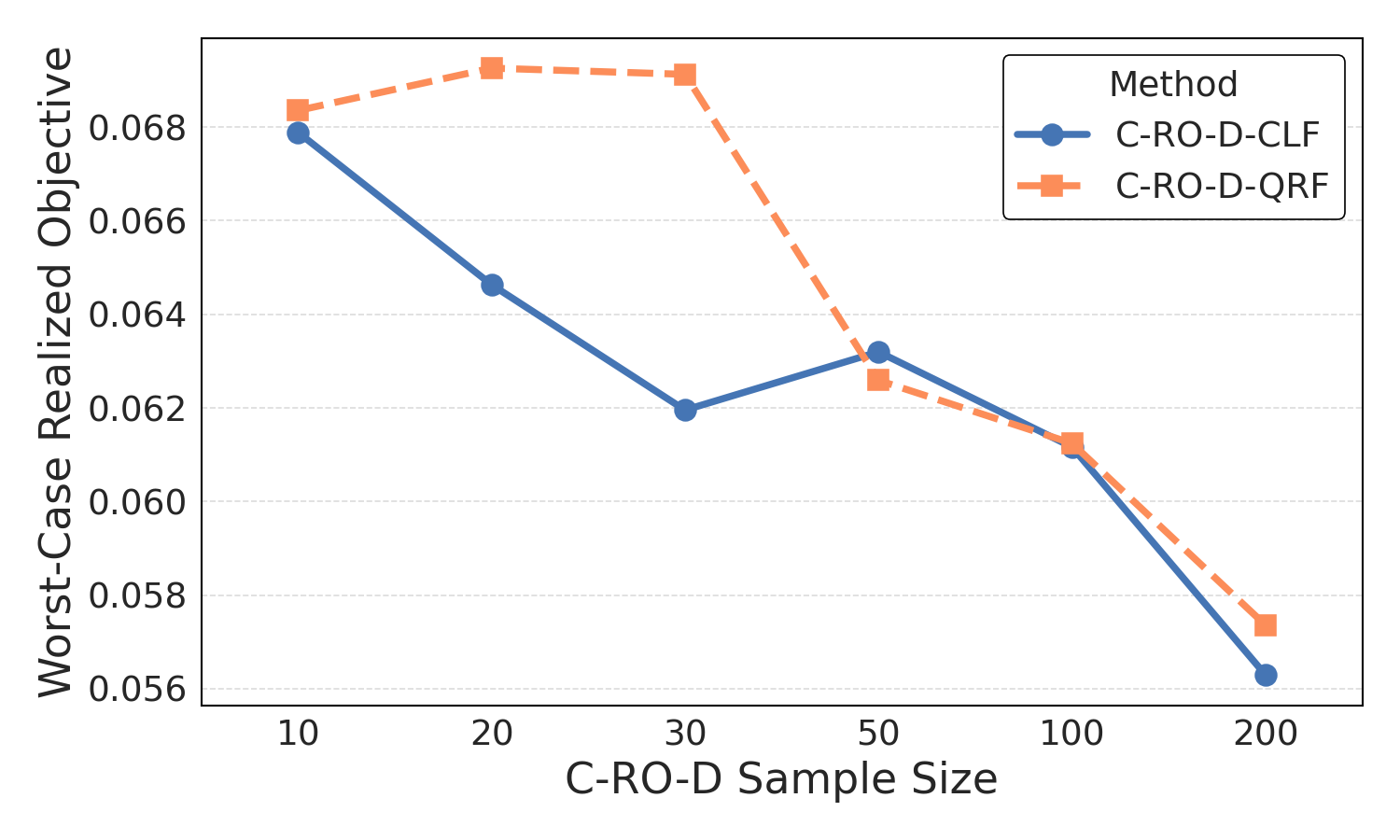} 
        \caption{Large Instances.}
    \end{subfigure}
     \caption{Worst-Case Realized Objective Value by C-RO-D Sample Size.}
     \label{fig:cro_d_size_obj}
\end{figure}

}

\end{APPENDIX}

\clearpage
\bibliographystyle{informs2014}
\bibliography{refs}

\end{document}